%% file: spcontv-arXiv.tex
\documentclass[11pt]{amsart}
\usepackage{amscd,amssymb,amsmath,latexsym,enumerate}
\usepackage{bbm}
\usepackage[mathscr]{euscript}
\usepackage{epsfig}
\usepackage{color}
\newcommand{\red}{\textcolor[rgb]{0.90,0.00,0.00}}

\usepackage[small,nohug,
]{diagrams}
\diagramstyle[labelstyle=\scriptstyle]
\newarrow{Corresponds}<--->
 
\newtheorem{theo}{Theorem}
\newtheorem{defini}{Definition}
\newtheorem{proposi}{Proposition}
\newtheorem{lemma}{Lemma}
\newtheorem{coro}{Corollary}
\newtheorem{rem}{Remark}
\newtheorem{exam}{Example}

\newcommand{\Aa}{{\mathcal A}}
\newcommand{\Bb}{{\mathcal B}}
\newcommand{\Cc}{{\mathcal C}}
\newcommand{\Dd}{{\mathcal D}}

\newcommand{\Ff}{{\mathcal F}}

\newcommand{\Hh}{{\mathcal H}}
\newcommand{\Ii}{{\mathcal I}}
\newcommand{\Jj}{{\mathcal J}}
\newcommand{\Kk}{{\mathcal K}}
\newcommand{\Ll}{{\mathcal L}}

\newcommand{\Pp}{{\mathcal P}}

\newcommand{\Uu}{{\mathcal U}}

\newcommand{\Ww}{{\mathcal W}}

\newcommand{\CM}{{\mathbb C}}

\newcommand{\NM}{{\mathbb N}}

\newcommand{\RM}{{\mathbb R}}

\newcommand{\SM}{{\mathbb S}}

\newcommand{\ZM}{{\mathbb Z}}

\newcommand{\CG}{{\mathfrak C}}

\newcommand{\FG}{{\mathfrak F}}

\newcommand{\hg}{{\widehat{g}}}
\newcommand{\hG}{{\widehat{G}}}

\newcommand{\hm}{{\widehat{m}}}

\newcommand{\hr}{{\widehat{r}}}

\newcommand{\hT}{{\widehat{T}}}

\newcommand{\hW}{{\widehat{W}}}

\newcommand{\as}{{\mathscr A}}

\newcommand{\cs}{{\mathscr C}}

\newcommand{\fs}{{\mathscr F}}
\newcommand{\gs}{{\mathscr G}}

\newcommand{\is}{{\mathscr I}}
\newcommand{\js}{{\mathscr J}}
\newcommand{\ks}{{\mathscr K}}

\newcommand{\ts}{{\mathscr T}}

\newcommand{\vs}{{\mathscr V}}

\newcommand{\tth}{{\tilde{h}}}

\newcommand{\tsig}{{\tilde{\sigma}}}

\newcommand{\Cs}{$C^{\ast}$-algebra }           
\newcommand{\CS}{$C^{\ast}$-algebra}            
\newcommand{\Css}{$C^{\ast}$-algebras }       	
\newcommand{\CsS}{$C^{\ast}$-algebras}          
\newcommand{\supp}{\textit{supp}}               
\newcommand{\oNM}{\overline{{\mathbb N}}}       
\newcommand{\GaO}{\Gamma^{(0)}}                 
\newcommand{\GaD}{\Gamma^{(2)}}                 
\newcommand{\UG}{\Uu\Gamma}                     
\newcommand{\isG}{\js(\Gamma)}                  
\newcommand{\tsG}{\ts(\Gamma)}                  
\newcommand{\tsGO}{\ts\Gamma^{(0)}}             
\newcommand{\tsGD}{\ts\Gamma^{(2)}}             
\newcommand{\tra}{\mbox{\sc t}}                 
\newcommand{\tri}{\mbox{\tiny\sc t}}            
\newcommand{\dist}{{\mbox{\rm dist}}}           
\newcommand{\Ima}{\mbox{\rm Im}}                
\newcommand{\Ker}{\mbox{\rm Ker}}               
\newcommand{\PGa}{\Gamma\times_r\Gamma}         
\newcommand{\AH}{L}                             
\newcommand{\spec}{{\rm Spec}}               

\begin{document}

\title{Spectral continuity for aperiodic quantum systems I.\\
        General Theory}
\thanks{Work supported in part by NSF Grant No. 0901514 and DMS-1160962 and SFB 878, M\"unster. GD's research is supported by the  FONDECYT grant \emph{Iniciaci\'{o}n en Investigaci\'{o}n 2015} - $\text{N}^{\text{o}}$ 11150143.}

\author{Siegfried Beckus, Jean Bellissard, Giuseppe De Nittis}

\address{Department of Mathematics\\
Technion - Israel Institute of Technology\\
Haifa, Israel}
\email{beckus.siegf@technion.ac.il}

\address{Georgia Institute of Technology\\
School of Mathematics\\
Atlanta GA, USA}
\address{Westf\"alische Wilhelms-Universit\"at\\
Fachbereich 10 Mathematik und Informatik\\
M\"unster, Germany}
\email{jeanbel@math.gatech.edu}

\address{
Facultad de Matem\'aticas \& Instituto de F\'{\i}sica\\
Pontificia Universidad Cat\'olica de Chile\\
Santiago, Chile}
\email{gidenittis@mat.uc.cl}

\begin{abstract}
How does the spectrum of a Schr\"odinger operator vary if the corresponding geometry and dynamics change? Is it possible to define approximations of the spectrum of such operators by defining approximations of the underlying structures? In this work a positive answer is provided using the rather general setting of groupoid $C^\ast$-algebras. A characterization of the convergence of the spectra by the convergence of the underlying structures is proved. In order to do so, the concept of continuous field of groupoids is slightly extended by adding continuous fields of cocycles. With this at hand, magnetic Schr\"odinger operators on dynamical systems or Delone systems fall into this unified setting. Various approximations used in computational physics, like the periodic or the finite cluster approximations, are expressed through the tautological groupoid, which provides a universal model for fields of groupoids. The use of the Hausdorff topology turns out to be fundamental in understanding why and how these approximations work.
\end{abstract}


\maketitle
\tableofcontents


\section{Introduction and main results}
\label{cspa.sect-MainResult}

\noindent The present paper delivers a characterization of the convergence of the spectrum of Schr\"odinger-type operators via the convergence of the underlying structures equipped with a suitable topology. This is the first paper of a series of articles aiming at providing a method to compute the spectrum of a self-adjoint operator on a Hilbert space. Such operators describe the quantum motion of a particle in a solid that might be aperiodic in its structure and/or submitted to a (not necessarily uniform) magnetic field. Special focus of this approach is to deal with aperiodic environments, in particular the ones with long range order, that have been mainly handled so far in one dimensional systems only. Our method delivers a general theory independent of the dimension by using $C^\ast$-algebraic techniques. Since this problem arises primarily in Solid State Physics, a special emphasis will be put on periodic approximations. This is because Physicists have at their disposal several numerical codes, based on Floquet-Bloch theory, permitting to compute the band spectrum of a periodic system.

\vspace{.1cm}

\noindent The present article is written for readers who have no expertise in the theory of groupoids and their \CsS. It is expected to serve as an introduction to these techniques for experts in Analysis and Spectral Theory. In order to do so, fundamental and known results in groupoid $C^\ast$-algebras are presented and proved for the convenience of the reader, a theory initiated in the classical work by Renault \cite{Re80}. Some of these results are slightly extended to cover all potential applications. In view of applications to physical problems, such as the electron motion in an aperiodic solid, submitted to a magnetic field, uniform or not, the case of $2$-cocycles depending on a parameter has been added. This addition does not change substantially the techniques used by experts in \Cs theory to prove the main results. But it has not been formally introduced in previous works. The central new concept introduced in this work, is the definition of the {\em tautological groupoid} of a groupoid $\Gamma$ (see Section~\ref{cspa.sect-taut}) providing the connection of the underlying structures with the spectra of the corresponding operators. It is a universal groupoid in classifying all continuous fields of groupoids with $\Gamma$ as their enveloping one 
(see Theorem~\ref{cspa.th-tautfield_2}).

\vspace{.1cm}

\noindent Within the tautological groupoid, the key concept is the use of the Hausdorff topology on the set of closed invariant subsets of the unit space, expressed in its various topological versions due to Vietoris \cite{Vi22}, Chabauty \cite{Ch50} and Fell \cite{Fe62}. The power of such a topology in dealing with various approximations has already been illustrated by various previous results obtained by the three authors over the last seven years and that can be found in the PhD Thesis of one of us \cite{Bec16} or in various talks that can be found online \cite{Be17}. It is expected that these results will be written in forthcoming publications.  Specifically, in the upcoming works it will be shown that Delone dynamical systems fall into this setting. These systems include  examples in Solid State Physics, such as the Sturmian sequences, the Penrose tiling, and, more generally, the model sets used to describe the structures of quasicrystals. As a warm up, it is shown here that this theory applies to a wide class of effective Hamiltonian and Schr\"odinger operators associated with dynamical systems. Within this future projects, the Hausdorff topology on the underlying structures is carefully analyzed explaining when Delone dynamical systems converge. Also the existence of, the construction of, and the obstruction to periodic approximations is on the schedule of this program. 

\vspace{.1cm}

\noindent {\bf Acknowledgements: } The three authors are indebted to Jean Renault for his invaluable email comments made to one of us (S.B.) about the validity of some arguments provided in his book \cite{Re80}. Many valuable comments were made during the last seven years during seminars, workshops and keynote presentations, by several participants, helping this work to be shaped in its present form. In particular, J.B. wants to thanks Pedro Resende about his input on the concept and the properties of open groupoids \cite{Re17}. He is indebted to Andr\'e Katz for a two-day conversation on defects for $2D$-quasicrystals \cite{Be92}. S.B. wants to express his deep gratitude to Daniel Lenz for his constant and fruitful discussions and support over the last years. Additionally, he wants to thank Anton Gorodetski for a stimulating discussion about the convergence of spectra, at BIRS, Oaxaca 2015. 
The authors wish to thank the anonymous referee for various fruitful observations that improved the work.

 \subsection{Proving the Continuity of the Spectrum: a Strategy}
 \label{cspa.ssect-strategy}

\noindent The first step of the strategy was proved in \cite{BB16} (see also the more detailed account in \cite{Aa16,Bec16}) and can be expressed as follows: let $T$ be a topological space. For each $t\in T$ let $\Hh_t$ be a Hilbert space and let $A_t$ be a bounded self-adjoint operator on $\Hh_t$. The family $A=(A_t)_{t\in T}$ is called a {\em field} of self-adjoint operators. It will be called $p2$-continuous whenever, for any polynomial $p\in \RM[X]$ of degree at most $2$, the map $t\in T\mapsto \|p(A_t)\|\in [0,\infty)$ is continuous.

\begin{theo}[\cite{BB16}]
\label{cspa.th-spcont}
A field $A=(A_t)_{t\in T}$ of self-adjoint operators is $p2$-continuous if and only if the spectral gap edges of the $A_t$'s are continuous in $t$, if and only if the spectrum $\spec(A_t)$ is continuous as a compact set {\em w.r.t.} the Hausdorff metric.   
\end{theo}

\noindent In \cite{BB16,Bec16}, the cases of fields of unitary or of unbounded self-adjoint operators are also treated. The results are similar in spirit.

\vspace{.1cm}

\noindent Proving the continuity of the norm of a field of operators might not always be easy though. So, it becomes convenient to enlarge the framework of this approximation method by considering continuous fields of \Css \cite{Ka51,TT61,To62,DD63,Di69}. If a field of \Css $\Aa=(\Aa_t)_{t\in T}$ is continuous then, from the definition, there is a dense algebra of continuous sections for which Theorem~\ref{cspa.th-spcont} applies. 
If it becomes possible to prove that the field $A$ can be seen as a continuous section of a continuous field of \CsS, then $A$ is automatically $p2$-continuous and the continuity of its spectrum is proved (see also \cite{Ka51}, Lemma 3.3). However, it is not necessarily easy to check that a field of \Css is actually continuous either.

\vspace{.1cm}

\noindent Thankfully there is a more accessible way to built continuous fields of \CsS. It is sufficient indeed to consider the convolution algebras of a field of groupoids. In \cite{LR01}, Landsman and Ramazan have investigated the definition of continuous field of groupoids and proved that the corresponding \Css make up a continuous field, at least if the groupoids are amenable. As it turns out, fields of groupoids are much easier to construct essentially because the construction boils down to building genuine topological spaces and dynamical systems. Examples of continuous fields of groupoids useful for Solid State Physics will be given in the second paper of this series \red{\cite{BBdN18}}. In the present paper, a general construction, valid for a large class of groupoids, will give rise to a {\em tautological} continuous field, indexed by the set of invariant subsets of the unit space. For example, a periodic orbit can be seen as an invariant subset of a groupoid. The main point will be to endow the set of invariant closed subsets with a suitable topology, which was first defined by Hausdorff for compact subsets of a metric spaces \cite{Ha14,Ha35}, then extended by Vietoris \cite{Vi22} and, more recently, by Chabauty \cite{Ch50} and Fell \cite{Fe62}. It will be called the {\em Hausdorff topology}. Then it will be shown that not every invariant set can be approximated by a family of periodic orbits. Characterizing invariant sets that can be approximated by periodic orbits  is exactly where the difficulty lies. In the second paper of this series criteria will be provided in various cases.

 \subsection{The Case of Dynamical Systems}
 \label{cspa.ssect-dynsyst}

\noindent Before going on, let the classical case of dynamical systems be considered. Let $(X,G,\alpha)$ be a topological dynamical system where $X$ is a second countable compact Hausdorff space, $G$ is a discrete countable group and $\alpha:G\times X\to X$ a continuous group action. In what follows $e$ will denote the neutral element of $G$. {\em For convenience $\alpha(x,g)=\alpha_g(x)$ will be denoted by $gx$ instead}. A subset $Y\subseteq X$ is called invariant if $gY\subseteq Y$ for all $g\in G$. The set of all non-empty closed invariant subsets of $X$ is denoted by $\js(G)$. The set $\js(G)$ equipped with the Hausdorff topology gets the structure of a second countable compact Hausdorff space.

\vspace{.1cm}

\noindent A generalized discrete Schr\"odinger operator is built as follows. Let $\ell^2(G)$ be the Hilbert space on which this operator acts. 
Then $G$ is represented by its left regular representation defined by
$$U_h\psi(g)=\psi(h^{-1} g)\,,
   \hspace{2cm}
    \psi\in\ell^2(G)\,,\, g,h\in G\,.
$$

\noindent Let $K\subseteq G$ be a {\em finite} subset such that $e\notin K$. Then for $k\in K$ let $t_k:X\to\CM$ be a continuous functions. Similarly let  $v:X\to\CM$ be continuous as well. For $x\in X$, let the operator $H_x:\ell^2(G)\to\ell^2(G)$ be defined by
%
\begin{equation}
\label{cspa.eq-schdis}
(H_x\psi)(g) \; := \; 
  \sum\limits_{h\in K} 
   t_h(g^{-1}x) \cdot \psi(gh)
    + v(g^{-1}x) \cdot \psi(g)\,,
\end{equation}

\noindent where $\psi\in\ell^2(G)$ and $g\in G$. The operator $H_x$ is linear bounded. In order to be self-adjoint the following are required:

(R1) $v$ is real valued,

(R2) $K$ is invariant by $k\to k^{-1}$,

(R3) the functions $t_k$ satisfy $t_{k^{-1}}(x)= \overline{t_k(k^{-1}x)}$.

\noindent The family of operators $H_X:=(H_x)_{x\in X}$ is obviously strongly continuous with respect to the variable $x\in X$. In addition, it is $G$-covariant, namely
$$U_h H_x U_h^{-1}= H_{hx}\,,
   \hspace{2cm}
    x\in X\,,\; h\in G\,.
$$
\noindent In particular, $\spec(H_x)=\spec(H_{hx})$ for all $h\in G$. Using strong continuity, it follows that if $y$ belongs to the closure $\Omega_x$ of the $G$-orbit of $x$, then $\spec(H_y)\subseteq \spec(H_x)$. It follows that if the $G$-action on $X$ is {\em minimal} (namely any orbit is dense), the previous inclusion holds both ways so that $\spec(H_y)= \spec(H_x)$ for any $x,y\in X$. It is worth noticing that such an equality between spectra for all $x,y\in X$ is equivalent to $(X,G,\alpha)$ being minimal (see \cite{Bec16}, Theorem 3.6.8, page 109, a result motivated by \cite{LS03}).

\vspace{.1cm}

\noindent It turns out that the family $H_X$ can be seen as an element (also denoted $H_X$) of the crossed product \Cs $\Aa=\Cc(X)\rtimes_\alpha G$ \cite{Be86}. The spectrum of this element in $\Aa$ can be shown to be \cite{Be86}
$$
\spec(H_X)=\overline{\bigcup_{x\in X}\spec(H_x)}\,.
$$
\noindent As will be seen, {\em amenability} is a crucial property required on the group $G$ (see \cite{Gr69} and Section~\ref{cspa.ssect-amen}) in order to get useful continuous fields. If $G$ is amenable, there is no need to take the closure so that $\spec(H_X)=\bigcup_{x\in X}\spec(H_x)$ holds \cite{Ex14,NP15}. 

\vspace{.1cm}

\noindent From the previous construction, if $Y\subseteq X$ is closed and $G$-invariant, then $(Y, G,\alpha)$ is also a dynamical system, so that $H_Y$ is well defined. The main question investigated in this work is: if $Y$ moves in the set of closed $G$-invariant subsets of $X$, how is the spectrum of $H_Y$ evolving~? The following result is a consequence of Theorem~\ref{cspa.theo-contSp}, in Section~\ref{cspa.ssect-groupoid}

\begin{theo}
\label{cspa.theo-ContSpectrGenSchrDyn}
Let $(X,G,\alpha)$ be a topological dynamical system where $G$ is an amenable discrete countable group. Let $\js(G)$ denotes the set of closed $G$-invariant subsets of $X$ equipped with its Hausdorff topology. Let $\ks(\RM)$ denotes the set of compact subsets of $\RM$, equipped with the Hausdorff metric. 

\noindent For each $Y\in \js(G)$, let $H_Y:=(H_x)_{x\in Y}$ be the generalized Schr\"odinger operator  defined in eq.~\eqref{cspa.eq-schdis} and satisfying (R1, R2, R3). Then the map 

$$\Sigma_H:\js(G)\to\ks(\RM)\,,
	\quad 
   Y\mapsto \spec(H_Y)\,,
$$

\noindent is continuous.
\end{theo}

\noindent It is worth noticing here that Theorem~\ref{cspa.theo-ContSpectrGenSchrDyn}, extends without further assumption to any self-adjoint element of the \CS, in particular allowing discrete Schr\"odinger operatos with infinite range. Furthermore, the convergence of non-empty invariant closed subsets is characterized by the convergence of the associated spectra for all generalized Schr\"odinger operators:


\begin{coro}
\label{cspa.cor-CharConvSpectrDyn}
Let $(X,G,\alpha)$ be a topological dynamical system where $G$ is an amenable discrete countable group. Then the following assertions are equivalent, whenever $Y_n\in\js(G)\,,\;n\in\oNM:=\NM\cup\{\infty\}$.

\begin{itemize}
\item[(i)] The sequence $(Y_n)_{n\in\NM}$ converges to $Y_\infty\in\js(G)$.
\item[(ii)] For all generalized Schr\"odinger operators $H_X:=(H_x)_{x\in X}$, the equation 

$$
\lim_{n\to\infty}\spec\big(H_{Y_n}\big)\;
	= \; \spec\big(H_{Y_\infty}\big)
$$ 

\noindent holds, where the limit is taken with respect to the Hausdorff metric on $\ks(\RM)$.
\end{itemize}
\end{coro}


\noindent {\bf Proof: } This Corollary is a consequence of
Corollary~\ref{cspa.cor-seqOid} in Section~\ref{cspa.ssect-groupoid}.
\hfill $\Box$

 \subsection{The Case of Groupoids}
 \label{cspa.ssect-groupoid}

\noindent The previous examples provide a large class of models liable to represent concrete situations in Condensed Matter Physics. However, a discrete translation group is not necessarily available in general. The most general situation in Condensed Matter Physics is the existence of a dynamical system $(X,G,\alpha)$ in which $G=\RM^d$ for some $d\in \NM$, called the {\em space dimension}. More precisely, in practical cases $X$ can be seen as a closed $G$-invariant subset of the space of Delone sets \cite{BHZ00,Be15}: each element of $X$ represents the instantaneous position of atomic nuclei and is a typical configuration of a solid or a liquid. Whenever $d=1$, discretization can be done using the concept of {\em Poincar\'e section}: namely it is a closed subset $\Xi\subseteq X$ such that for any $x\in \Xi$ the smallest $t\geq 0$ such that $\alpha_t(x)\in \Xi$ is strictly positive for any $x$. This concept can be generalized as follows: let $(X,G,\alpha)$ be a topological dynamical system, where $X$ is a compact Hausdorff second countable space, $G$ is a locally compact Hausdorff second countable group (where $G$ may not be discrete), and $\alpha$ is a continuous action of $G$ on $X$ by homeomorphisms. Then a closed subset $\Xi\subseteq X$ is called a {\em transversal} if for each $x\in \Xi$ there is an open neighborhood $U$ of the neutral element $e\in G$ such that for $g\in U$ then $gx\notin \Xi$. Then the set $\Gamma_\Xi$ of pairs $(x,g)\in \Xi\times G$ such that $gx\in \Xi$ becomes a groupoid (see below) \cite{Co79,Re80}. But in general it is not a topological dynamical system.

\vspace{.1cm}

\noindent To treat these cases, the concept of groupoid is more universal. A precise definition will be provided in Section~\ref{cspa.sect-grremind}. But a groupoid can be described as a 
small 
category with all morphisms invertible, for which the family of objects and of morphisms makes up a set. An object of this category is called a {\em unit} and $\GaO$ will denote the set of units. A morphism $\gamma:x\to y$ between the units $x$ and $y$ is called an {\em arrow} with {\em range} $r(\gamma)=y$ and {\em source} $s(\gamma)=x$. Then $\Gamma$ denotes the set of arrows. By definition, all arrows are invertible. In addition, any unit $x\in\GaO$ can be associated with a unique arrow $e_x$ such that $\gamma \circ e_x=\gamma = e_y\circ \gamma$ if $\gamma :x\to y$. Consequently $\GaO$ can be seen as a subset of $\Gamma$. An ordered pair of arrows $(\gamma,\gamma')\in \Gamma\times \Gamma'$ will be called {\em composable} if $s(\gamma)=r(\gamma')$ (morphism composition) in which case their composition is denoted by $\gamma\circ \gamma'$. The set of composable pairs is denoted by $\GaD$. At last, given $x\in \GaO$, then $\Gamma^x$ will denote the set of arrows with range $x$ and $\Gamma_x$ will denote the set of arrow with source $x$. 

\vspace{.1cm} 

\noindent In this work, all such groupoids will be endowed with a topology, making $\Gamma$ a {\em locally compact Hausdorff second countable space}. This topology will be such that the range and the source maps $r,s$, the inversion map $\gamma\to \gamma^{-1}$ and the composition map $(\gamma,\gamma')\in \GaD\mapsto \gamma\circ \gamma'\in \Gamma$ are continuous. It follows that $\GaO$ and $\GaD$ are closed. In addition, in this work,  all groupoids considered will be assumed to satisfy

\vspace{.1cm}

(i) their unit space is compact,

(ii) the two maps $r,s$ being open.

\vspace{.1cm}

\noindent Such groupoids will be called {\em handy} (see Section~\ref{cspa.ssect-defgr}), because they are convenient and practical. Moreover, they arise naturally in all known realistic problems motivated by Physics. The requirement that the unit space be compact is not necessary as will be seen in Section~\ref{cspa.ssect-haar}. But the space of closed invariant subset of the unit space has a unique Hausdorff topology only if the unit space is compact. Otherwise, the Vietoris and the Chabauty-Fell topologies may differ (see a discussion in \cite{BB16}).


\vspace{.1cm}

\noindent The analog of the case of a topological dynamical system with a discrete group action is provided by {\em \'etale groupoids}. $\Gamma$ will be called \'etale whenever any $\gamma\in \Gamma$ admits a neighborhood $U$ such that the map $\gamma'\in U\mapsto r(\gamma')\in \GaO$ is an homeomorphism onto its image.

\vspace{.1cm} 

\noindent A typical example of a groupoid is provided by a topological dynamical system as in Section~\ref{cspa.ssect-dynsyst}. In the latter case $\GaO=X$ and $\Gamma$ is the set of pairs $\gamma=(x,g)\in X\times G$ with 

(i) $r(x,g)=x$, $s(x,g)=g^{-1}x$, 

(ii) $(x,g)\circ(g^{-1}x, h)=(x,gh)$,

(iii) $e_x= (x,e)$ if $e\in G$ is the neutral element,

(iv) $(x,g)^{-1}=(g^{-1}x,g^{-1})$\,.

\noindent The topology on $\Gamma$ is defined by the product topology of $X\times G$. Such a groupoid is denoted by $X\rtimes_{\alpha} G$ and is called the {\em crossed product} of $X$ by $G$.

\vspace{.1cm} 

\noindent Another example is the {\em groupoid of a transversal}. Namely if $\Xi\subseteq X$ is a transversal, then the same definition applied to $\Gamma_\Xi$ gives a locally compact groupoid, which is {\em \'etale}.  

\vspace{.1cm} 

\noindent Any locally compact group admits a Haar measure. The convolution of functions allows to define an algebra of functions, and ultimately a \CS. For groupoid such a concept can be defined and will be called a {\em Haar system} (see Section~\ref{cspa.ssect-haar}) \cite{Re80}, or a {\em transverse function} \cite{Co79}. It is worth mentioning that not all locally compact groupoid admit a Haar system. However, a deep result by Blanchard \cite{Bl96} permits to prove that a groupoid is handy if and only if (i) its unit space is compact and (ii) it admits a Haar system (see Section~\ref{cspa.ssect-haar}, Theorem~\ref{cspa.th-haaropen}).  

\vspace{.1cm}

\noindent Correspondingly, the analog of the convolution algebra, defined with complex valued continuous functions on $\Gamma$, defines a \Cs after completion. However, this \Cs is not unique in general. The smallest such \Cs is called the {\em reduced} \Cs and is denoted by $\CG_{red}^\ast(\Gamma)$. 
In order to describe a Hamiltonian submitted to a magnetic field, a normalized $2$-cocycle $\sigma:\Gamma^{(2)}\to \SM^1$ is involved in the algebraic structure. In this case   $\CG_{red}^\ast(\Gamma,\sigma)$ will denote the corresponding \CS.  This algebra (with or without the twisting induced by the $2$-cocyle) is unique if $\Gamma$ is {\em amenable} \cite{AR00,AR01,BH14}.

\vspace{.1cm} 

\noindent Given a groupoid $\Gamma$ an {\em invariant set} is a subset $M\subseteq \GaO$ such that if $x\in M$ then any $\gamma\in \Gamma^x$ satisfies $s(\gamma)\in M$. If $\Gamma$ is a Hausdorff locally compact groupoid, let $\isG$ denotes the set of closed invariant subsets of $\GaO$, equipped with the Hausdorff topology. It will be shown that, whenever $\Gamma$ is handy, $\isG$ is compact. It follows that, for each $M\in\isG$ the set of arrows $\gamma$ with range in $M$ makes up a closed subgroupoid $\Gamma_M$. If $f: \Gamma\to \CM$ is a continuous function with compact support then its restriction to $\Gamma_M$ will be denoted by $f_M$. This restriction is continuous with compact support. It follows that the restriction map $f\to f_M$ extends as a $\ast$-homomorphism from $\CG_{red}^\ast(\Gamma,\sigma)$ onto $\CG_{red}^\ast(\Gamma_M,\sigma_M)$ where $\sigma_M$ is the restriction of the $2$-cocycle $\sigma$ to the subgroupoid $\Gamma_M$. It is worth reminding at this point that a {\em normal} element $f$, in a \Cs, is an element such that $f^\ast f=f\,f^\ast$.

\vspace{.1cm} 

\noindent The main result of this paper can be summarized as follows (see the proof in Section~\ref{cspa.sect-th3}):

\begin{theo}
\label{cspa.theo-contSp}
Let $\Gamma$ be an amenable handy groupoid and $\sigma$ a normalized $2$-cocycle. For every normal element $ f\in\CG^\ast_{red}(\Gamma,\sigma)$, the spectral map

$$\Sigma_ f:\isG\to\ks(\CM)\,, 
   \hspace{2cm}
    M\mapsto \spec(f_M)\,,
$$

\noindent is continuous.
\end{theo}

\noindent  {\bf Proof of Theorem~\ref{cspa.theo-ContSpectrGenSchrDyn}: } Theorem~\ref{cspa.theo-ContSpectrGenSchrDyn} immediately follows from Theorem~\ref{cspa.theo-contSp}. For the amenability of the groupoid $X\rtimes G$ follows from the amenability of the group $G$ \cite{AR00,AR01}.
\hfill $\Box$


\begin{rem}
\label{cspa.rem-amen}
{\em The amenability assumption is required only to make sure that the reduced and the full \Css coincide. In view of \cite{Wil15}, this coincidence may occur without the groupoid $\Gamma$ being amenable, in which case the previous Theorem applies. 
}
\hfill $\Box$
\end{rem}

\begin{coro}
\label{cspa.cor-seqOid}
Let $\Gamma$ be an amenable handy groupoid 
twisted by the normalized $2$-cocycle $\sigma$.
Consider a sequence $X_n\in\isG$ of invariant, closed sets for each $n\in\oNM=\NM\cup\{\infty\}$. Then the following assertions are equivalent.

(i) The sequence $(X_n)_{n\in\NM}$ converges to $X_\infty$ in the Hausdorff-topology of $\isG$.

(ii) For all self-adjoint (normal) $f\in\CG^\ast_{red}(\Gamma,\sigma)$, the equation

$$\lim_{n\to\infty} \spec\big(f_{X_n}\big)= \spec\big(f_{X_{\infty}}\big) \,,
$$ 

\noindent holds, where the limit is taken with respect to the Hausdorff metric $\ks(\RM)$.
\end{coro}

\noindent {\bf Proof: } (i)$\Rightarrow$(ii) follows by Theorem~3.

\vspace{.1cm}

\noindent (ii)$\Rightarrow$(i): Assume that $(X_n)_{n\in\NM}$ does not converge to $X_\infty\in\isG$. Since $\isG$ is compact, there is a subsequence $(X_{n_k})_{k\in\NM}$ such that $\lim_{k\to\infty} X_{n_k}=:Y$ exists and $Y\neq X_\infty$. There are two cases, namely (a) $X_\infty\setminus Y\neq\emptyset$ and (b) $Y\setminus X_\infty\neq\emptyset$. Both cases lead as follows to a contradiction. 

\vspace{.1cm}

\noindent (a) If $X_\infty\setminus Y$ is non-empty, consider an $x\in X_\infty\setminus Y$. Using the Lemma of Urysohn, there exists a continuous function 
$a:\Gamma\to[0,1]$ 
with compact support such that $a(x)=1$ and $a_{\Gamma_Y}\equiv 0$. Without loss of generality we can assume that $a$ is self-adjoint (otherwise consider $a^\ast a$). Then $\|a|_{X_\infty}\|>0$ while $\|a|_Y\|=0$. Thus, $\lim_{n\to\infty} \spec\big(a_{X_n}\big)= \spec\big(a_{X_{\infty}}\big)$ does not hold contradicting (ii).

\noindent (b) If $Y\setminus X_\infty\neq\emptyset$, a similar argument like in (a) leads to a contradiction.
\hfill $\Box$


 \subsection{How to Construct Approximations~?}
 \label{cspa.ssect-contapp}

\noindent What is the strategy to built some approximations, using Theorem~\ref{cspa.theo-contSp} or the Corollary~\ref{cspa.cor-seqOid}~? Here are some guesses. 

  \subsubsection{Periodic Approximations}
  \label{cspa.sssect-perapp}

\noindent Let $G$ be assumed to be $\ZM^d$ for instance (or more generally a free Abelian group), and let $(Y,\ZM^d,\beta)$ be a minimal dynamical system, with $Y$ a second countable compact Hausdorff space. Let assume that for each $n\in\NM$ there is a minimal periodic dynamical system $(Y_n,\ZM^d,\beta_n)$. By periodic, it is meant that there is a subgroup $H_n\subseteq \ZM^d$ of finite index fixing each point of $Y_n$. The minimality of $Y_n$ implies that it is finite then. Let $X$ be the disjoint union of the $Y_n$'s and of $Y$. Clearly $\ZM^d$ acts on $X$ through the action $\alpha$ such that $\alpha$ coincides with $\beta$ on $Y$ and with $\beta_n$ on $Y_n$. The next step will be to find a topology on $X$ making it a second countable compact Hausdorff space, so that (i) $Y$ and the $Y_n$ are closed $\ZM^d$-invariant subsets, and that (ii) $Y_n$ converges to $Y$ in the Hausdorff topology. In many cases such a strategy works (see for instance \cite{Pr13} for disordered systems and \cite{Bec16} for other examples). If so, then the previous Corollary~\ref{cspa.cor-CharConvSpectrDyn} applies and gives a sequence of periodic approximations to the generalized Schr\"odinger operator $H_Y$ permitting to approximate its spectrum. A full characterization of the existence of periodic approximations with explicit constructions in dimension $d=1$ was recently proved in \cite{BBdN18}.

  \subsubsection{Finite Cluster Approximations}
  \label{cspa.sssect-fincl}

\noindent Another approximation used in Physics, consists in cutting a finite piece of an infinite distribution of atoms, namely to choose a cluster (see for instance \cite{KS86}). Let $\Ll$ be a distribution of atoms: namely it will be a Delone set. In such a set, the Voronoi construction gives rise to a tiling by polyhedrons \cite{Be15} centered at each atom. Two atoms are nearest neighbors (n.n.) whenever their Voronoi tiles intersect along a facet (face of codimension one). If an edge is a pair of n.n., this gives a graph defining $\Ll$. In order to prove that this approximation is appropriate, let $C_1\subseteq \cdots \subseteq C_n\subseteq \cdots$ be an increasing family of such clusters such that the union is all the atomic sites. Each such cluster gives rise to a finite groupoid $\Gamma_n$, generated by the edges of the graph associated with it. In the infinite volume limit, there is a groupoid $\Gamma$ describing all such constructions. Then let $\Gamma_{tot}$ be obtained as the disjoint union of these groupoids, which is also a groupoid. If there is a topology on this union making it a continuous field over the set $\NM\cup\{\infty\}$, then the generalization to groupoids of the Corollary~\ref{cspa.cor-CharConvSpectrDyn}, will give convergence of the spectrum as well.

 \subsection{Methodology: a Review}
 \label{cspa.ssect-meth}

\noindent The problem of investigating the spectral properties of aperiodic quantum system has a long history, both in Physics and in Mathematics. 

  \subsubsection{Transfer Matrix Method}
  \label{cspa.sssect-trmat}

\noindent For one-dimensional models, the {\em transfer matrix} method has been quite successful. The earliest result obtained for Schr\"odinger's operators with non periodic coefficient is probably the one provided by Dinaburg and Sinai \cite{DS75}, using this method to provide results about the spectral properties of the $1D$-Schr\"odinger operator on the real line with a quasiperiodic potential. The method used for the KAM theorem turns out to be effective. This method was used successfully since the early eighties on similar problems on the discrete line. Eventually, it gave rise to the theory of cocycles. A spectacular achievement, using this method, has been the solution of the {\em Ten Martini Problem} \cite{AJ09}.

\vspace{.1cm}

\noindent However the transfer matrix method is limited in several ways. First it applies only to one dimensional systems. Secondly, it applies only to the case of finite range Hamiltonians. Nevertheless, this method was used successfully to investigate numerically the case of a random potential (Anderson model) using a transfer matrix in high dimension for systems of infinite length but with large finite transverse area \cite{KMcK93}. However, all attempt made to use this method for systems in higher dimension in order to get rigorous mathematical results, have failed so far, if one except models that can be decomposed into a direct sum or a product of commuting one-dimensional systems \cite{DGS15}.

  \subsubsection{Renormalization Group Method}
  \label{cspa.sssect-RG}

\noindent In the eighties and the early nineties the case of $1D$-discrete Schr\"odinger Operators with potential provided by a substitution was illustrated by several examples. The first example was the {\em Fibonacci Hamiltonian}, also called the {\em Kohmoto model} \cite{KKT83,OPRSS83,OK85,Ca86,Su87}. The spectrum could be computed through a Renormalization Group method induced by the substitution applied to the transfer matrix \cite{BIST89,Be90,BBG91}. Since then, the method has been refined and pushed to the limit, thanks to the two decades persistent work of D. Damanik and A. Gorodetski. A spectacular achievement, for instance can be the recent review on the Fibonacci Hamiltonian \cite{DGY16}.

\vspace{.1cm}

\noindent An attempt was made to use this approach in higher dimension \cite{SB90}. For indeed, most quasicrystals studied so far admits a special symmetry (mostly the icosahedral symmetry in $3D$) and, as a consequence admits also an inflation symmetry playing the role of a space Renormalization Map. However, the method turned out to be hopeless unless in the region of energy for which the kinetic energy can be treated as a very small perturbation. Example of Hamiltonian in high dimension that can be considered as Cartesian product of $1D$-operators, have shown that the qualitative results concerning the spectrum can be surprisingly complicate \cite{Si89}. In particular the RG-method is unlikely to give so much information.

  \subsubsection{Other Numerical Approaches}
  \label{cspa.sssect-numeric}

\noindent Soon after the discovery of quasicrystals \cite{SBGC84}, one of the earliest numerical calculations were made on a large cluster of the Penrose lattice \cite{Ch85,ON86,KS86}. The cluster size (varying from 251 to 3806) and the nature of the parameters in the model investigated, may change the results. The existence of inflation rules gives a fast algorithm to grow a cluster from a seed. This gives a sequences $C_1,\cdots,C_n$ of cluster with size increasing exponentially fast. However, this inflation rule can hardly been used to implement relations between the various Hamiltonians $H_1,\cdots,H_n$ in order to get more rigorous results on its spectral properties \cite{SB90}. The finite size effect were not systematically estimated, but they are visible. It is expected that the errors, as measured in Density of States (DOS), be of the order of the ratio $\{boundary\;area(C_n)\}/\{Volume(C_n)\}$ namely of the order of the inverse of the cluster diameter. This is a slow convergence.

\vspace{.1cm}

\noindent Another method, promoted in \cite{SB91}, consists in finding a periodic approximation and to use the Bloch Theory to compute the spectrum like in \cite{Pr13}. In most situations investigated so far, it can be proved that the error is much smaller and decay exponentially fast with the diameter of the unit cell of the periodic lattice.

  \subsubsection{Defect Creation}
  \label{cspa.sssect-def}

\noindent Periodic approximations may come with a prize. Depending upon how the sequence of approximations is built, defects might show up in the limit. This was demonstrated in \cite{OK85} and proved rigorously in \cite{BIT91} for the case of $1D$-quasicrystals (also called Sturmian sequences). As it turns out, such defects also occur, the so-called {\em worms}, in the $2D$-case for a nearest neighbor Hamiltonian on the quasicrystal with octagonal symmetry. This quasicrystal can be computed by mean of a projection from $\ZM^4$ onto $\RM^2$ respecting the $8$-fold symmetry. But the projection $2$-plane is not in generic position {\em w.r.t.} the lattice $\ZM^4$, precisely because of the symmetry. As a result, generically, a periodic approximation will produce defects \cite{Be92}. However, there is a way, in this example, to find a sequence of periodic approximation avoiding the occurrence of defects, as was proposed in \cite{DMO89}. In the second paper of this series, this problem will be addressed and solved from using the Anderson-Putnam complex \cite{AP98}.


\section{Groupoids: a Reminder}
\label{cspa.sect-grremind}

\noindent This Section offers a summary of the definitions and properties of groupoids that will be needed later. Most of this material can be found in \cite{Co79,Re80,Ka82}. 

 \subsection{Definitions}
 \label{cspa.ssect-defgr}

A {\em groupoid} $\Gamma$ is a small category in which all morphisms $\gamma: x\to y$ are invertible. By a small category it is meant that both the families of objects and morphisms are sets. Then $\Gamma$ will denote the set of morphisms. Let $\GaO$ denote the set of objects. If $\gamma:x\to y$, the notation will be $s(\gamma)=x$ which will be called the {\em source}, and $r(\gamma)=y$ which will be called the {\em range}. Similarly a morphism $\gamma$ will be called an {\em arrow}. Then $x$ can be identified with the morphism $\gamma^{-1}\circ \gamma=e_x$ and $y$ with $\gamma\circ \gamma^{-1}=e_y$. Then $x=s(e_x)=r(e_x)$ and $e_x\circ\gamma=\gamma$ if $r(\gamma)=x$, while $\gamma\circ e_x=\gamma$ when $s(\gamma)=x$. Hence $e_x$ is a unit and is the identity morphism for the object $x$, so that the map $x\to e_x$ gives a bijection identifying $\GaO$ with the {\em set of units}. 
Then $\Gamma^{(0)}$ is called 
the {\em unit space}, instead, and can then be seen as a subset of $\Gamma$. It will be convenient though to separate the concept of unit, the set of which will be sometimes denoted by $\UG$, from the concept of object in practice, even if the two sets are in bijection. If $ (\gamma_1, \gamma_2)\in\Gamma\times \Gamma$ the composition $\gamma_1\circ \gamma_2$ makes sense only if $s(\gamma_1)=r(\gamma_2)$. The notation will be
$$\GaD= 
   \{(\gamma_1,\gamma_2)\,;\,s(\gamma_1)=r(\gamma_2)\}\,.
$$

\noindent For $y\in \GaO$, its fibers are the $r$-fiber $\Gamma^y=\{\gamma\in \Gamma\,;\, r(\gamma)=y\}$, the $s$-fiber $\Gamma_y=\{\gamma\in \Gamma\,;\, s(\gamma)=y\}$ and the $rs$-fiber $\Gamma_y^y=\{\gamma\in \Gamma\,;\, r(\gamma)=s(\gamma)=y\}$ which is a group. The relation {\em ``$x$ is isomorphic to $y$''} , namely $\exists \gamma\in \Gamma\,,\, \gamma:x\to y$, is an equivalence relation which will be denoted by $x\sim y$. Given $F\subseteq \GaO$, its saturation $[F]$ denotes the set $\{x\in \GaO\,;\, \exists y\in F\,,\, x\sim y\}$. Given $F\subseteq \GaO$, $\Gamma_F$ will denote the subgroupoid $\Gamma_F=\{\gamma \in \Gamma\,;\, s(\gamma)\in F\,,\, r(\gamma) \in F\}$. A subset $F\subseteq \GaO$ will be called {\em invariant} if $F=[F]$. The concept of subgroupoid is straightforward and left to the reader.

\begin{exam}[Groups]
\label{cspa.exam-group}
{\em Let $G$ be a group with neutral element $e\in G$. Then $\Gamma:=G$ is a groupoid where $\Gamma^{(2)}=G\times G$. In addition, the range $r(g)$ and the source $s(g)$ are equal to $e$. Thus, $\GaO=\{e\}$.
}
\end{exam}

\begin{exam}[Sets]
\label{cspa.exam-set}
{\em Let $X$ be a set. Then $\Gamma:=X$ is a groupoid with $\GaD:=\{ (x,y)\in X\times X\;|\; x=y\}$, composition $x\circ x:=x$ and inverse $x^{-1}:= x$ for $x\in X$. Thus, the unit space $\GaO=X$.
}
\end{exam}

\begin{exam}[Equivalence Relation]
\label{cspa.exam-eqrel}
{\em Let $X$ be a set and let $\sim$ be an equivalence relation on $X$. Then the graph of the equivalence relation becomes a groupoid $\Gamma:=\{(x,y)\in X\times X\;|\; x\sim y\}$ with $\GaD:=\big\{\big((x_1,x_2),(x_3,x_4)\big)\;\big|\; x_2=x_3\big\}$, with the composition $(x,y)\circ(y,z):=(x,z)$ and inverse $(x,y)^{-1}:=(y,x)$. Then $\GaO=X$ with $e_x=(x,x)$.
}
\end{exam}

\begin{exam}[Dynamical Systems]
\label{cspa.exam-dynsys}
{\em Let $(X,G,\alpha)$ be a dynamical system  with $X$ a set, $G$ a group and $\alpha: (x,g)\in X\times G \to gx\in X$ a left action of $G$ on $X$ by bijections. Then $\Gamma= X\times G$ becomes a groupoid if $\GaO=X$, $r(x,g)=x$, $s(x,g)=g^{-1}x$, and the composition is defined by $(x,g)\circ(g^{-1}x,h)=(x,gh)$. Then $(x,g)^{-1}=(g^{-1}x,g^{-1})$. This groupoid will be denoted by $X\rtimes_\alpha G$ and will be called {\em crossed product}.
}
\hfill $\Box$
\end{exam}

\begin{exam}[Induced Subgroupoid]
\label{cspa.exam-trans}
{\em If $(X,G,\alpha)$ is a dynamical system, and if $Y\subseteq X$, then the set $\Gamma_Y=\{(y,g)\in Y\times G\,;\, g^{-1}y\in Y\}$ is a subgroupoid of $X\rtimes_\alpha G$.
}
\hfill $\Box$
\end{exam}

\begin{defini}
\label{cspa.def-topgr}
A groupoid $\Gamma$ is {\em topological} whenever the set of arrows $\Gamma$ and it unit space $\GaO$ admit a topology such that 

(i) the range map $r$, the source map $s$, are continuous, 

(ii) the composition and the inverse are continuous.
\end{defini}

\noindent A topological space $X$ is Hausdorff if given any two points $x,y\in X$ such that $x\neq y$, there are open sets $U\ni x$ and $V\ni y$ such that $U\cap V=\emptyset$. Equivalently this condition holds if and only if the diagonal $\Delta_X=\{(x,x)\in X\times X\,;\, x\in X\}$ is closed in $X\times X$. In such a case, for $x\in X$ the singleton $\{x\}$ is a closed subset.

\begin{lemma}
\label{cspa.lem-close}
(i) Let $\Gamma$ be a topological groupoid. If $\GaO$ is Hausdorff, then the $r,s$ fibers and $\GaD$ are closed. 

\noindent (ii) If, in addition $\Gamma$ is Hausdorff, then the set of units $\UG$ is closed inside $\Gamma$.   
\end{lemma}

\noindent  {\bf Proof: }(i) Since $\GaO$ is Hausdorff, any singleton $\{x\}$ in $\GaO$ is closed and the diagonal $\Delta$ in $\GaO\times\GaO$ is closed. Since $r,s$ are continuous, the preimages $r^{-1}(\{x\})=\Gamma^x$ and $s^{-1}(\{x\})=\Gamma_x$ are closed, implying that $\Gamma_y^x=\Gamma^x\cap \Gamma_y$ is closed as well. In addition let $\Phi:\Gamma\times \Gamma\to  \GaO\times \GaO$ be the map defined by $\Phi(\gamma_1,\gamma_2)=(s(\gamma_1),r(\gamma_2))$. This map is continuous by assumption and the preimage of the diagonal $\Delta$ is precisely $\GaD$, hence it is closed. 

\vspace{.1cm}

\noindent (ii) Let $\UG$ denotes the set of units inside $\Gamma$, namely the set of $e\in\Gamma$ such that there is $\gamma\in\Gamma$ for which $\gamma^{-1}\circ\gamma=e$. It follows that $r(e)=s(e)$ and that, from the associativity of the product, $e\circ e=e$ implying $e^{-1}\circ e= e\circ e^{-1}=e$. Consequently $e$ is a unit if and only if $r(e)=s(e)$ and $e\circ e=e$. The map $\xi : \Gamma \to \Gamma\times\Gamma$ given by $
\xi(\gamma)=  (\gamma,\gamma\circ\gamma^{-1})$ is continuous. The diagonal 
$\Delta\subseteq \Gamma\times\Gamma$ is closed if $\Gamma$ is Hausdorff. It follows that $\UG=\xi^{-1}(\Delta)$ hence proving that it is closed.
\hfill $\Box$

\vspace{.2cm}

\noindent Among topological groupoids, the {\em open} ones will be important. As a reminder, a continuous map $f:X\to Y$ between two topological spaces is called {\em open} whenever the image of any open set of $X$ by $f$ is open in $Y$.

\begin{defini}
\label{cspa.def-open}
A topological groupoid $\Gamma$ will be called open if its range and source maps are both open.
\end{defini}

\noindent It is worth remarking that because of the continuity of the map $i:\gamma\mapsto \gamma^{-1}$, the range map is open if and only if the source map is. In addition since $i$ is its own inverse, its continuity impies that $i$ is open as well. Being open is actually sufficient to define the {\em tautological groupoid} of a groupoid $\Gamma$ (see Section~\ref{cspa.sect-taut}). The following result is known in category theory \cite{Re17}, but it is worth giving a more topological proof.

\begin{lemma}
\label{cspa.lem-mopen}
Let $\Gamma$ be an open groupoid. Then, the multiplication map $m:\Gamma^{(2)}\to \Gamma$ is open.
\end{lemma}

\noindent  {\bf Proof: }(i) First let $\PGa\subseteq \Gamma\times \Gamma$ be the set of pairs $(\gamma_1,\gamma_2)$ such that $r(\gamma_1)=r(\gamma_2)$. Let $\pi_i\,,\, i=1,2$ be the projection onto the $i$-th coordinate defined by $\pi_i(\gamma_1,\gamma_2)=\gamma_i$. Then $\pi_i$ is open. For indeed, let $U_1,U_2$ be two open sets in $\Gamma$ and let $W$ be the trace of $U_1\times U_2$ on $\PGa$. Then $\pi_1(W)= U_1\cap r^{-1}(r(U_2))$. Since $r$ is both continuous and open, it follows that $r(U_2)$ is open, so that  $r^{-1}(r(U_2))$ is open too, showing that $\pi_1(W)$ is open. A similar argument holds for $\pi_2$. 

\vspace{.1cm}

\noindent (ii) Let now $\xi:\PGa\to \Gamma^{(2)}$ be the map defined by $\xi(\gamma_1,\gamma_2)=(\gamma_1,\gamma_1^{-1}\circ \gamma_2)$. This map is continuous. Moreover it is invertible with inverse $\xi^{-1}(\gamma,\gamma')= (\gamma, \gamma\circ \gamma')$ as can be checked immediately. In particular $\xi^{-1}$ is also continuous, so that $\xi$ is a homeomorphism. Thus it is clear that $m=\pi_2\circ \xi^{-1}$, which is open as composition of open maps.
\hfill $\Box$

\begin{coro}
\label{cspa.cor-unitopen}
Let $\Gamma$ be an Hausdorff open groupoid with Hausdorff unit space. Then the map $x\in\GaO\mapsto e_x\in \UG$ is a homeomorphism.
\end{coro}

\noindent  {\bf Proof: }Since $\Gamma$ and $\GaO$ are Hausdorff, it follows, from Lemma~\ref{cspa.lem-close} that $\UG$ and all fibers are closed subsets in $\Gamma$. In addition, thanks to Lemma~\ref{cspa.lem-mopen}, the product map is open. Let $U\subseteq \UG$ be an open set in the set of units for the induced topology. Let $\hr$ denote the restriction of the range map $r$ to $\UG$. Then let $\psi:\Gamma\to \UG$ be the continuous map defined by $\psi(\gamma)=\gamma\circ \gamma^{-1}$. Hence $W=\psi^{-1}(U)$ is open. In addition $r(W)=\hr(U)$, showing that $\hr(U)$ is open in $\GaO$. Thus the map $\hr:e\in\UG\mapsto r(e)\in\GaO$ is open. Since this map is open, continuous, one-to-one and onto, it follows that its inverse is continuous as well, so that this map is a homeomorphism.
\hfill $\Box$

\vspace{.2cm}

\noindent 

\noindent The groupoid analog of a discrete group is provided by

\begin{defini}
\label{cspa.def-etale}
A topological groupoid $\Gamma$ is \'etale whenever any $\gamma\in \Gamma$ admits an open neighborhood $U$ such that the map $\gamma'\in U\mapsto r(\gamma')\in \GaO$ is an homeomorphism onto its image.

\noindent Equivalently, $\Gamma$ is \'etale whenever its topology is generated by a basis of bisections, namely subsets $F\subseteq \Gamma$ such that the restriction of the range map to $F$ is a homeomorphism.
\end{defini}

\noindent Thanks to the continuity of the inverse map $i$, the range map $r$ can be replaced by the source map $s$ as well in the previous definition.

\vspace{.1cm}

\noindent The {\em canonical transversal} of the {\em Hull} of a {\em tiling} in $\RM^d$ with finite local complexity gives rise to an \'etale groupoid with Cantor unit space \cite{Kel97,BJS10,JS12}. \'Etale groupoid with Cantor unit space have recently been the focus of attention from the group theory community in order to produce groups of intermediate growth with remarkable properties, like simplicity, amenability, finite presentation or finitely generated \cite{Ne15}. This is because the set of bisections of an \'etale groupoid with range or source given by the whole unit space can be identified with a group of homeomorphisms of the unit space, called the {\em full group} of $\Gamma$. The concept of full group was introduced in \cite{Kr79,GPS99} for $\ZM^d$-actions on a Cantor set. The full group of a Cantor set was proved to be amenable in a remarkable paper \cite{JM13} .

\begin{lemma}[\cite{Re80}, Prop. 2.8]
\label{cspa.lem-discfib}
An \'etale groupoid $\Gamma$ has discrete fibers.
\end{lemma}

\noindent  {\bf Proof: }Let $x\in\GaO$. If $\gamma\in\Gamma^x$, let $U$ be an open neighborhood of $\gamma$ making the map $r_U:U\to r(U)\subseteq \GaO$ an homeomorphism. Since $r_U$ is one-to-one, it follows that $\Gamma^x\cap U$ is reduced to $\gamma$. Hence $\gamma$ is isolated inside $\Gamma^x$, showing that $\Gamma^x$ is discrete. The same argument works for the $s$-fiber.
\hfill $\Box$

\begin{lemma}
\label{cspa.lem-etale}
Any \'etale topological groupoid is open.
\end{lemma}

\noindent  {\bf Proof: } Let $\Gamma$ be an \'etale groupoid. Let $U$ be an open subset in $\Gamma$. If $x\in r(U)$, there is $\gamma\in U$ such that $x=r(\gamma)$. Since $\Gamma$ is \'etale, there is an open set $V_\gamma\subseteq U$, containing $\gamma$, on which the restriction of the source map is an homeomorphism on its image. By construction, $x\in r(V_\gamma)\subseteq r(U)$. Hence $r(U)$ is open. The same argument works for the source map.
\hfill $\Box$

\vspace{.2cm}

\noindent Most groupoids used to describe physical systems have additional properties that make them remarkable and more convenient. In this work they will be called {\em handy}. 

\begin{defini}
\label{cspa.def-handy}
A topological groupoid $\Gamma$ will be called handy whenever:

(i) it is locally compact, Hausdorff, second countable, 

(ii) its unit set is compact,

(iii) the source map and the range map are open.

\noindent It follows that $\GaO$ is Hausdorff.  
\end{defini}

\noindent Examples of handy groupoids are obtained from previous examples in an obvious way: (i) Example~\ref{cspa.exam-group} if the group is locally compact, Hausdorff, second countable; (ii) Example~\ref{cspa.exam-set} where $X$ is a compact Hausdorff second countable set; (iii) Example~\ref{cspa.exam-eqrel} if $X$ is a compact, Hausdorff, second countable space, and if the graph of the equivalence relation is closed.

\begin{exam}[Topological Dynamical System]
\label{cspa.exam-topdyn}
{\em Let $(X,G,\alpha)$ be a dynamical system. If (i) $X$ is compact, Hausdorff, second countable, (ii) $G$ is a locally compact, Hausdorff, second countable group, (iii) the action $\alpha: (g,x)\to gx$ is continuous and the maps $\alpha_g:x\in X\to gx\in X$ are homeomorphism, (iv) $X\rtimes_\alpha G$ is endowed with the topology of the Cartesian product $X\times G$. Then $X\rtimes_\alpha G$ is a handy groupoid. The proof is a consequence of the Theorem~\ref{cspa.th-haaropen} below.
}
\hfill $\Box$
\end{exam}


\begin{exam}
\label{cspa.exam-subgr}
{\em If $\Gamma$ is a handy groupoid, any closed sub-groupoid is handy. In particular, if $Y\subseteq \GaO$ is closed, then the induced groupoid $\Gamma_Y$ is handy.
}
\hfill $\Box$
\end{exam}

\begin{rem}
\label{cspa.rem-handy}
{\em The condition of second countability allows to use Urysohn's metrization theorem and the Tietze extension theorem. In practice it is not much of a restriction. The condition of local compactness is convenient to construct a \Cs by convolution. Without it, not much is known. The separation axiom is not absolutely needed \cite{Co82}, but it makes the groupoid much easier to handle. The compactness of the unit space is not necessary, but almost all examples used to describe physical systems so far, have this property. The openness of the source and range maps will turn out to be necessary when dealing with the tautological groupoid (see Section~\ref{cspa.sect-taut}).
}
\hfill $\Box$
\end{rem}

 \subsection{Continuous Haar Systems}
 \label{cspa.ssect-haar}

\noindent A locally compact group $G$ is known to admit a left invariant positive measure, called the {\em Haar measure} \cite{We40}, which is unique modulo a normalization factor. It follows that the space $L^1(G)$ can be endowed with an algebraic structure using the convolution product. Similarly an {\em adjoint} can also be defined. The algebra obtained in this way, denoted by $L^1(G)$, is called the {\em Wiener algebra}. It also gives rise to a \Cs denoted by $ \CG^\ast(G)$. If $G$ is Abelian, the Gelfand theorem and the Pontryagin duality give the Fourier transform as a $\ast$-isomorphism between $ \CG^\ast(G)$ and the space $\Cc_0(\hG)$ of continuous function of the Pontryagin dual $\hG$ vanishing at infinity.

\vspace{.1cm}

\noindent Such a construction can be extended to groupoids. The Haar measure is replaced by what has been called a {\em Haar system} \cite{We67,Se75,Se76,Re80} or a {\em transverse function} \cite{Co79}. The latter concept is actually helping in understanding the concept of {\em transverse measure} \cite{Co79}, which permits to express all positive weights of the \Cs of a groupoid in an analytic way. In particular it gives rise to a noncommutative analog of integration theory.

\begin{defini}
\label{cspa.def-HaarSyst}
A left-continuous Haar system on a topological groupoid $\Gamma$ is a family of positive Borel measures $(\mu^x)_{x\in\GaO}$ on $\Gamma$ such that the following assertions hold.

\begin{description}

\item[(H1)] For each $x\in\GaO$, the support of the measure $\mu^x$ is equal to the $r$-fiber $\Gamma^x$.

\item[(H2)] For each $f\in\Cc_c(\Gamma)$, the function $\GaO\to\CM$ defined by
$$x\mapsto \mu^x(f)=
   \int_{\Gamma^x} f(\gamma) d\mu^x(\gamma)\,.
$$

\noindent is continuous.

\item[(H3)] These measures are left-invariant, namely if $\gamma:x\to y$ then
$$
\int_{\Gamma^x} f(\gamma\circ\eta) d\mu^x(\eta)=
 \int_{\Gamma^y} f(\eta') d\mu^y(\eta')\,.
$$

\noindent for any $f\in\Cc_c(\Gamma)$.

\end{description}

\end{defini}

\noindent In this definition, the conditions (H1) and (H3) are similar to the case of a Haar measure on a group. However (H2) is new. It turns out to be necessary and sufficient so that the convolution defined on $\Cc_c(\Gamma)$ maps into $\Cc_c(\Gamma)$ see \cite{Se86} (see also Section~\ref{cspa.ssect-convol} below). 

\vspace{.1cm}

\noindent It is worth mentioning that not all locally compact groupoid admits a Haar system: a simple counter example can be founded in \cite{Gr09}, for which the axiom (H2) is violated. However every \'etale groupoid $\Gamma$ has one, using the counting measure
$$\mu^x(f) =\sum_{\gamma\in \Gamma^x} f(\gamma)\,,
   \hspace{2cm}
    f\in\Cc_c(\Gamma)\,.
$$

\noindent Indeed, since the fibers are discrete, if $f$ has compact support, the sum admits only a finite number of nonzero terms. By construction, the support of $\mu^x$ is $\Gamma^x$. By definition of ``\'etale'' the map $x\to \mu^x(f)$ is continuous. The left invariance comes from the remark that if $\gamma:x\to y$, then the map $\eta\in \Gamma^x\mapsto \gamma\circ\eta\in \Gamma^y$ is a bijection, since $\gamma$ is invertible. In particular this applies to Example~\ref{cspa.exam-set}, for the groupoid defined by a compact Hausdorff second countable set $X$, where $\mu^x(f)=f(x)$ is the Dirac measure. It is worth mentioning that any other Haar system on a handy \'etale groupoid is given by $\nu^x=g(x)\mu^x$ where $g:\GaO\to \RM_+$ is a non-negative continuous function such that $x\sim y\Rightarrow g(x)=g(y)$. In particular if $[\{x\}]$ is dense in $\GaO$, it follows that the counting Haar measure is the unique Haar system on $\Gamma$ modulo normalization. 

\vspace{.1cm}

\noindent Another important example is given by a topological dynamical system $(X,G,\alpha)$ (see Example~\ref{cspa.exam-topdyn}): since $G$ is locally compact it admits a Haar measure denoted by $dh$ here, leading to
$$\mu^x(f) =\int_G f(x,h)\,dh\,,
   \hspace{2cm}
    f\in\Cc_c(\Gamma)\,.
$$

\noindent It is easy to check that the axioms are satisfied. The previous two examples are the most important in view of applications to Physics. 

\vspace{.1cm}

\noindent The following characterization is a deep result

\begin{theo}[\cite{We67,Bl96}]
\label{cspa.th-haaropen}
A second countable locally compact Hausdorff topological groupoid admits a left-continuous Haar system if and only if it is open.
\end{theo}

\noindent The proof that such a groupoid admitting left-continuous Haar system is open is due to Westman \cite{We67}. The proof can be found below. The proof of the converse will be skipped here. But it can be found in \cite{Wi16} and is based upon a deep result by Blanchard \cite{Bl96}. 

\vspace{.1cm}

\noindent  {\bf Westman's Proof: } Here $\Gamma$ is a second countable locally compact Hausdorff topological groupoid admitting a left-continuous Haar system. Using the continuity of the inverse map, it is sufficient to show that the range map is open. Let $U\subseteq\Gamma$ be an open set and let $\gamma\in U$. Since $\Gamma$ is second countable, locally compact and Hausdorff, Urysohn's Lemma applies. Thus, there is a continuous function $f\in\Cc_c(\Gamma)$ satisfying $0\leq f \leq 1$, $f(\gamma)=1$ and $\supp(f)\subseteq U$. Let the interior of the support of $f$ be denoted by $V_\gamma$. It is an open subset of $U$ by construction. Hence $\mu^x(f)>0$ holds if and only if $x\in r(V_\gamma)$. Since the map $\Phi_f:x\in \GaO \mapsto \mu^x(f) \in \CM$ is continuous by (H2), the preimage $\Phi_f^{-1} \big((0,\infty)\big) \subseteq \GaO$ is open. By the previous considerations, $\Phi_f^{-1}\big((0,\infty)\big)$ is equal to $r(V_\gamma)$. Hence, $r(V_\gamma)$ is an open subset of $r(U)$ containing $r(\gamma)$. Altogether, $r(U)=\bigcup_{\gamma\in U} r(V_\gamma)$ is a union of open sets implying that $r(U)$ is open in $\Gamma^{(0)}$.
\hfill $\Box$

\begin{coro}
\label{cspa.cor-reghaar}
Any second countable locally compact topological groupoid, with compact unit space and with a left-continuous Haar system, is handy.
\end{coro}

 \subsection{Cocycles and Gauge Fields}
 \label{cspa.ssect-cocy}

\noindent A fair part of the book by J.~Renault \cite{Re80} is dedicated to the problem of groupoid cohomology. He introduced the concept of $n$-cocycles. For the purpose of the construction of the convolution algebra only 
the cases $n=2$ is really relevant. 
Here, $\Gamma$ will be a handy groupoid. It is convenient to introduce $\Gamma^{(n)}$ to be the subset of $\Gamma^{\times n}$ made of families $(\gamma_1,\cdots,\gamma_n)$ such that $s(\gamma_{k})=r(\gamma_{k+1})$ whenever $1\leq k <n$. Such a family defines a {\em path} within $\Gamma$. Since $\Gamma$ is Hausdorff, $\Gamma^{(n)}$ is a closed subset of $\Gamma^{\times n}$, thus a locally compact space. Let $\AH$ denote an Abelian topological group, which will be called the {\em group of coefficients}. The composition law in $\AH$ will be written as an addition and therefore the neutral element will be written as $0$. However, in application, $\AH$ will often be the unit circle $\SM^1$ seen as the multiplicative subgroup of $\CM_\ast=\CM\setminus \{0\}$ or the group of unitaries in an Abelian \CS. {\em The reader is invited to proceed to the change of notations between an additive and a multiplicative group law}. Then an {\em $n$-cochain} with coefficients in $\AH$ is a 
{continuous map $\sigma:\Gamma^{(n)}\to \AH$.
Here $\AH$ denotes any topological Abelian group. 
The set of such cochains will be denoted by $C^n(\Gamma,\AH)$. For the pointwise addition it is an Abelian group. The differential map $d$ is defined as a group homomorphism $d: C^n(\Gamma,\AH)\to C^{n+1}(\Gamma,\AH)$ such that
%
\begin{eqnarray}
\label{cspa.eq-dchain}
d\sigma(\gamma_0,\cdots,\gamma_{n})&=&
   \sigma(\gamma_1,\cdots,\gamma_{n})+ \\
&&    \sum_{k=1}^{n} (-1)^k
\sigma(\gamma_0,\cdots,\gamma_{k-2},\gamma_{k-1}\circ\gamma_k,\gamma_{k+1}, \cdots, \gamma_{n})\nonumber\\
  &&+(-1)^{n+1}\sigma(\gamma_0,\cdots,\gamma_{n-1})\nonumber\,.
\end{eqnarray}

\noindent It can be checked that $d^2\sigma=0$. An $n$-cocycle $\sigma$ is an 
$n$-cochain 
with $d\sigma =0$, namely $\sigma$ is {\em closed}. An $n$-coboundary $\sigma$ is an 
$n$-cochain for which there is an $(n-1)$-chain 
$\tau$ such that $\sigma=d\tau$. For the present purpose, only 
$1$- and $2$-cocycles matter. 

\vspace{.1cm}

\noindent {\bf The case $n=1$: }it can be checked that a $1$-cocycle is a groupoid homomorphism $\sigma: \Gamma\to \AH$, namely the condition $d\sigma=0$ becomes $\sigma(\gamma_0\circ \gamma_1)=\sigma(\gamma_0)+\sigma(\gamma_1)$. By definition, it is a $1$-coboundary if and only if it vanishes.

\begin{exam}
\label{cspa.exam-module}
{\em Let $X$ be a compact Hausdorff second countable space. Let $\tra$ be an homeomorphism of $X$. Then $\tra$ defines a $\ZM$-action. Let $\Gamma=X\rtimes_{\tri} \ZM$ be the corresponding handy (\'etale) groupoid. Then a {\em unitary module} $\delta$ is a $1$-cocycle with coefficients in $\AH=\SM^1$. Hence it can be checked that $\delta(x,0)=1$, while $\delta$ is entirely defined by the function $h:x\in X\mapsto \delta(x,1)\in \SM^1$, since, for any $n\in\NM$
$$\delta(x,n)\;=\;h(x)h(\tra^{-1}x)\cdots h(\tra^{-n+1}x)\,,
   \hspace{2cm}
    \delta(x,-n)\;=\;\overline{h(\tra x)h(\tra^2 x)\cdots h(\tra^n x)}\,.
$$
}
\hfill $\Box$
\end{exam}

\noindent {\bf The case $n=2$: }It will be convenient for the rest of the paper to denote by $\AH$ a compact Abelian group for which the group law will be denoted multiplicatively. Then  a $2$-cocycle is a continuous function on $\GaD$ with values in $\AH$ such that
%
\begin{equation}
\label{cspa.eq-coceq}
\frac{\sigma(\gamma_0,\gamma_1\circ \gamma_2)}
      {\sigma(\gamma_0\circ\gamma_1,\gamma_2)}\;=\;
 \frac{\sigma(\gamma_0,\gamma_1)}{\sigma(\gamma_1,\gamma_2)}\,.
\end{equation}

\noindent In addition $\sigma$ is a coboundary if and only if there is a continuous function $\tau:\Gamma\to \AH$ such that
$$\sigma(\gamma_0,\gamma_1)\;=\;
   \frac{\tau(\gamma_0)\tau(\gamma_1)}{\tau(\gamma_0\circ \gamma_1)}\,.
$$

\begin{exam}[Uniform magnetic fields: see \cite{Be86}]
\label{cspa.exam-magn}
{\em The main example relevant for Physics is the effect of a magnetic field. Let $X$ be a compact Hausdorff second countable space. Let $\tra$ be a $\ZM^2$-action on $X$ and let $\Gamma=X\rtimes\ZM^2$ be the corresponding handy (\'etale) groupoid. Then an arrow can be labeled by a pair $(x,n)$ with $x\in X$ and $n=(n_1,n_2)\in\ZM^2$. A unitary $2$-cocycle (namely with $\AH=\SM^1$) can be defined by 
$$\sigma\left((x,n),(\tra^{-n}x,m)\right)\;=\;
  e^{\imath B\; n\wedge m}\,,
   \hspace{2cm}
    B\in\RM\,,\;\; n\wedge m=n_1m_2-n_2m_1\,,
$$

\noindent as can be check easily. Here $B$ plays the role of a magnetic field perpendicular to the plane of the physical space in which the lattice $\ZM^2$ lies \cite{Be86}. Such a cocycle is not a coboundary, except for special values of $B$. 
}
\hfill $\Box$
\end{exam}

\vspace{.2cm}

\noindent A $2$-cocycle $\sigma$ is called {\em normalized} if $\sigma(e_x,e_x)=1$ for all $x\in\Gamma^{(0)}$. Throughout this work, normalized $2$-cocycles are considered.

\begin{lemma}
\label{cspa.lem-1coc}
Let $\sigma$ be a $2$-cocycle on $\Gamma$ with values in $\AH$. If $x\in\GaO$ then 
$$\sigma(e_x,\gamma)\;=\;
   \sigma(e_x,e_x)\;=\;
    \sigma(\gamma^{-1},e_x)\,,
     \hspace{1.5cm}
        \forall \gamma\in\Gamma^x\,.
$$

\noindent If, additionally, $\sigma$ is normalized, then $\sigma(\gamma,\gamma^{-1})=\sigma(\gamma^{-1},\gamma)$ holds for all $\gamma\in\Gamma$.
\end{lemma}

\noindent  {\bf Proof: }Use eq.~\eqref{cspa.eq-coceq} with $\gamma_0=e_x$ and $x=r(\gamma_1)$. It gives
$$\frac{\sigma(e_x,\gamma_1\circ \gamma_2)}
      {\sigma(\gamma_1,\gamma_2)} \;=\;
 \frac{\sigma(e_x,\gamma_1)}{\sigma(\gamma_1,\gamma_2)}\,.
$$

\noindent Consequently $\sigma(e_x,\gamma_1)=\sigma(e_x,\gamma_1\circ \gamma_2)$, for any choice of $\gamma_1,\gamma_2$ compatible with the definition. Choosing $\gamma_1=e_x$ and $\gamma_2=\gamma$ this gives the first set of identities. Now, choosing $\gamma_1=e_x$ with $x=s(\gamma_0)=r(\gamma_2)$, the same eq.~\eqref{cspa.eq-coceq} gives
$$
1\;=\;
	\frac{\sigma(\gamma_0,\gamma_2)}{\sigma(\gamma_0,\gamma_2)}
	\;=\; \frac{\sigma(\gamma_0,e_x)}{\sigma(e_x,\gamma_2)}\,,
	$$
\noindent 
which gives the second set of identities for $\gamma_0:=\gamma^{-1}$ and $\gamma_2:=\gamma$.

\vspace{.1cm}

\noindent Fianlly, $\sigma(\gamma,\gamma^{-1})=\sigma(\gamma^{-1},\gamma)$ follows from eq.~\eqref{cspa.eq-coceq} for $\gamma_0=\gamma_2=\gamma$ and $\gamma_1=\gamma^{-1}$ if $\sigma$ is normalized.
\hfill $\Box$

\begin{coro}
\label{cspa.cor-normcoc}
Let $\sigma$ be a $2$-cocycle on $\Gamma$ with values in $\AH$ and let $\tau:\Gamma\to \AH$ be a continuous function. Then the function $\tsig$ defined by
$$\tsig(\gamma_0,\gamma_1)\;=\;
   \sigma(\gamma_0,\gamma_1)\,
    \frac{\tau(\gamma_0)\tau(\gamma_1)}{\tau(\gamma_0\circ\gamma_1)}
$$

\noindent defines also a $2$-cocycle. If then $\tau$ is chosen so that $\tau(\gamma)=\sigma(e_x,e_x)^{-1}$ where $x=r(\gamma)$, then $\tsig$ is normalized.
\end{coro}

 \subsection{An Example: Non Uniform Magnetic Fields}
 \label{cspa.ssect-numagn}

\noindent In a series of articles \cite{MP02,MP05,MPR05}, M\u{a}ntoiu and his collaborators proposed a pseudo-differential calculus for studying the Schr\"odin\-ger operator describing the quantum motion of a charged particle submitted to both an electric potential and a non uniform magnetic field. To this end, they built a \Cs in which the product is defined through a cocycle taking values in an Abelian \CS. The same idea was also used by Rieffel in \cite{Rie89a} to build a more general theory of continuous fields of \CsS. In this Section, it will be shown that there is a way around such an extension to get an $\SM^1$-valued cocycle instead. In order to do so, the idea developed here is the same developed in \cite{Be86,Be93} to build the Hull. The general method is developed first and the specific for magnetic field will be explained afterwards. It is worth mentioning \cite{PR16} in which the case of a continuous field of cocycles corresponding to a family of magnetic fields, is proved to lead to continuity of the spectrum.

  \subsubsection{Hull of uniformly continuous functions}
  \label{cspa.sssect-hullUC}

\noindent Let $G$ be a locally compact Lindel\"of group. For simplicity it will be assumed to be Hausdorff and second countable. The neutral element of $G$ will be denoted by $e$ and the composition law, not necessarily commutative, will be denoted multiplicatively. Then a complex valued function $F:G\to \CM$ is called {\em uniformly continuous} whenever for any $\epsilon >0$, there is an open neighborhood $W$ of the neutral element $e$ such that $|F(g)-F(h)|<\epsilon$ for any pair $(g,h)\in G$ such that $hg^{-1}\in W$. Any such function is obviously continuous. Let $\Cc_u(G)$ denote the set of all complex valued, bounded, uniformly continuous function on $G$. It can be checked that the uniform limit of a uniformly convergent sequence of uniformly continuous functions is also unifomly continuous. Consequently, equipped with the pointwise addition, product, complex conjugacy and the uniform norm $\|F\|=\sup_{g\in G}|F(g)|$, $\Cc_u(G)$ becomes a unital, Abelian \CS. In addition, $G$ acts on $\Cc_u(G)$ through a group $\alpha$ of $\ast$-automorphisms defined by $\alpha_h(F)(g)=F(h^{-1}g)$. Thanks to the uniform continuity it follows that $\lim_{h\to e}\|\alpha_h(F)-F\|=0$ for any $F\in\Cc_u(G)$. Hence the group $\alpha$ is norm-pointwise continuous. Let now $\fs\subset \Cc_u(G)$ be a countable set. Then $A(\fs,G)$ will denote the sub-\Cs of $\Cc_u(G)$ generated by the family $\fs(G)=\{\alpha_g(F)\,;\, g\in G\,,\, F\in \fs\}$ of $G$-translated of elements of $\fs$. The \Cs $A(\fs,G)$ is unital, Abelian, so, by Gelfand's Theorem (see for instance \cite{Ar81,Co90}), it is isomorphic to $\Cc(\Omega_\fs)$ where $\Omega_\fs$ is a compact space. Since the family $\fs$ is countable and since $G$ is second countable, it follows that $A(\fs,G)$ is separable, so that $\Omega_\fs$ is also second countable. It is Hausdorff by construction. Actually, $\Omega_\fs$ is built as the {\em spectrum}, namely the set of {\em characters} defined as the $\ast$-homomorphisms $x:A(\fs,G)\to \CM$. It follows that $G$ acts by duality on $\Omega_\fs$ using $gx=x\circ \alpha_{g^{-1}}$. Consequently, $(\Omega_\fs,G)$ is a genuine topological dynamical system, where $G$ acts by homeomorphisms. Moreover, the Gelfand isomorphism is defined by the map $F\in A(\fs,G)\mapsto \gs F\in \Cc(\Omega_\fs)$ where the Gelfand transform $\gs$ is defined by $\gs F(x)=x(F)$. In $\Omega_\fs$ there is a remarkable point namely $x_e$ defined by $x_e(F)=F(e)$. It is straightforward to show that the $G$-orbit of $x_e$ is dense in $\Omega_\fs$, so that the Gelfand transform is entirely defined by the equation $\gs F(g^{-1}x_e)=F(g)$. Following \cite{Be86,Be93}

\begin{defini}
\label{cspa.def-hull}
The compact space $\Omega_\fs$, endowed with its canonical $G$-action and defined as the Gelfand spectrum of the \Cs generated by the $G$-translated elements of the family $\fs$, will be called the hull of $\fs$. 
\end{defini}

  \subsubsection{Non Uniform Magnetic Field: the continuum case}
  \label{cspa.sssect-contNU}

\noindent The idea is to use this hull to built a groupoid and a cocycle with values in $\SM^1$ to represent the algebra of observables behind a Schr\"odinger operator with magnetic field. To illustrate the method, let $\Ll$ be a Delone set in $\RM^d$ \cite{Be15}, supposed to represent the location of the atomic nuclei of an assembly of atoms making up a condensed material (solid or liquid). If there are several species of atoms, each atom is labeled by a letter $l$ in a finite alphabet $\as$. Then let $\Ll_l\subseteq \Ll$ denote the subset of positions of the atomic nuclei of species $l$. Let $V$ denote the potential energy seen by an individual valence electron, liable to travel through the material. With a very good approximation $V$ can be described as 
$$V(x)=\sum_{l\in\as} \sum_{q\in\Ll_l} v_l(x-q)\,,
$$

\noindent where $v_l(x)$ is the effective atomic Coulomb potential created by atoms of the $l$-species at the position $x$ relative to the position of it nucleus. It includes both the Coulomb attraction by the nucleus and the screening effect of the core electrons. A magnetic field is usually described as a closed two-form $B=\sum_{i,j=1}^d B_{ij}(x)dx_i\wedge dx_j$. In most cases the coefficients $B_{ij}$ are smooth functions of the variable $x\in\RM^d$. It will be assumed here that {\em these coefficients are uniformly continuous} over the whole $\RM^d$. Such a condition is satisfied for instance if the gradient of the coefficients are uniformly bounded\footnote{If the magnetic field is not bounded at infinity, the spectral properties of the Schr\"odinger operator are qualitatively different from the bounded case.} 
on $\RM^d$. Using the Poincar\'e construction, $B$ can be locally built from a vector-potential $A=\sum_{i=1}^d A_i(x)dx_i$, a one-form satisfying $dA=B$, that is $B_{ij}=\partial_i A_j-\partial_j A_i$. However, as can be seen from the case of a constant magnetic field, $A$ is not bounded in general. The Schr\"odinger operator describing the quantum motion of a particle of mass $m$ and charge $-e$ moving in the potential $V$, is the PDE
$$H= -\frac{\hbar^2}{2m}\,
   \sum_{i=1}^d(\partial_i-\imath e A_i)^2
    +V\,.
$$

\noindent It is a standard result to show that there is a dense domain $\Dd$ in $\Hh=L^2(\RM^d)$ such an operator, defined on $\Dd$ is essentially self-adjoint (see for instance \cite{RS75}). The kinetic term is the square of the generator of translations. But the presence of the magnetic field induces a change in the translation operators, leading to the concept of {\em magnetic translations} \cite{Zak64}. Namely the usual translation operators $T(a)$ by the vector $a\in\RM^d$, acting on the quantum Hilbert space of states $\Hh$, is replaced by (see for instance \cite{Zak64})

\begin{equation}
\label{cspa.eq-magtr}
U(a)\psi (x)=
   e^{\imath (e/\hbar)\,\oint_{x-a}^x A}\,\psi(x-a)\,,
  \hspace{2cm}
   \psi\in\Hh\,.
\end{equation}

\noindent where $\oint_{x-a}^x A$ denotes the line integral of the one-form $A$ along the straight line joining $x-a$ to $x$. This is the usual rule used in Physics when a gauge field is present, namely to change the momentum operator $p$, which is the infinitesimal generator of translations, into $p-eA$. The commutation rules for these unitary operators is given by
%
\begin{equation}
\label{cspa.eq-comrul}
U(a)U(b)U(a+b)^{-1}= e^{\imath (e/\hbar)\, \Phi(0,a,b)}\,,
\end{equation}

\noindent where $\Phi(0,a,b)$ represents the map $x\in\RM^d\mapsto \Phi_x(0,a,b)$, namely the magnetic flux through the oriented triangle defined by the ordered family of points $(x,x-a,x-a-b)$. If $(0,a,b)$ represents the oriented triangle with vertices $(0,a,a+b)$, 
and if $\partial(0,a,b)$ represents its oriented boundary, 
it can be written as
%
\begin{equation}
\label{cspa.eq-flux}
\Phi_x(0,a,b) =
 \int_{(0,a,b)}\sum_{i,j=1}^d B_{ij}(x-s) ds_i\wedge ds_j\,=\,
  \oint_{\partial(0,a,b)} A(x-s) ds\,.   
\end{equation}

\noindent seen as a multiplication operator acting on $\Hh$. To implement a good product on the set of observable, this rule must be used, leading to the appearance of a cocycle (see for instance \cite{Be86,Be93}).
In the case of a uniform magnetic field, the cocycle is described in Example~\ref{cspa.exam-magn}. But if $B$ is not uniform, the main difficulty at this point is that the magnetic translation is no longer independent of the position $x$ at which the state is evaluated. The solution to this problem proposed by M\u{a}ntoiu {\em et. al.}, consists in seeing the phase factor in the {\em r.h.s.} of eq.~\eqref{cspa.eq-comrul}, as a bounded continuous function of $x$, and therefore as an element of a ``good'' sub-algebra of $\Cc_b(\RM^n)$.

\begin{rem}
\label{cspa.rem-PotVec}
{\em Physicists know from Classical Hamiltonian Mechanics that the vector potential enters explicitly in the Hamiltonian, not the magnetic field. They wondered whether this is physical or not. Indeed, the vector potential is defined modulo an exact $1$-form, namely a gauge. However, the Aharonov-Bohm effect \cite{AB59}, confirmed by experiments, showed that indeed, in interference experiments, the quantum particle feels the magnetic field even in regions where it vanishes provided that, in this region, the vector potential does not vanish. This is why, {\em representing} the magnetic translation in the Hilbert space of quantum states require the use of the vector potential and cannot be defined otherwise. Changing gauge, gives rise to a unitary equivalent representation of these translations. However, when lifting the question to the level of the observable algebras this problem disappears, because only the magnetic flux is required to define the corresponding \CS. Indeed, as will be seen later in Section~\ref{cspa.ssect-convol}, the definition of the \Cs requires only the $2$-cocycle, which depends only upon the magnetic field itself.  
}%
\hfill $\Box$
\end{rem}

  \subsubsection{The combined Hull}
  \label{cspa.sssect-combHull}

\noindent This difficulty admits another solution though, using the groupoid approach of the present work. First, if the Delone set $\Ll$ is not periodic, it has a hull $\Omega_\Ll$, a Hausdorff, second countable compact space. This hull can be defined as the closure of the orbit of $\Ll$ under the translation group $G=\RM^n$, after a genuine topology is defined on the space of Delone sets (see for instance \cite{BHZ00}). 
It is also worth noticing that $\Omega_\Ll$ admits a canonical transversal, namely $\Xi_\Ll=\{\Ll'\in \Omega_\Ll\,;\, 0\in \Ll'\}$.
This Hull allows to see the electric potential $V$ as a function $V_\omega (x)= \vs(\tra^{-x}\omega)$ for a continuous function $\vs$ on the compact space $\Omega_\Ll$ and an $\RM^d$-action by homeomorphism $\tra$. 
Then, since $G=\RM^d$ is locally compact and Hausdorff and since the components $B_{i,j}$ of the magnetic field are assumed to be uniformly continuous, the family $\fs= \{ B_{i,j}\,;\, 1\leq i<j\leq d\}$ admits another hull, namely the Gelfand spectrum of the Abelian \Cs $A(\fs, G)$, which will be denoted by $\Omega_B$. The latter is another Hausdorff, second countable, compact space, and is endowed with an $\RM^d$-action which, for simplicity, will also be denoted $\tra$. This leads to a product hull $\tilde{\Omega}=\Omega_\Ll\times \Omega_B$ endowed with the diagonal $\RM^d$-action $\tra^x(\omega, \omega')=(\tra^x\omega,\tra^x\omega')$. The product hull $\tilde{\Omega}$ admits also a remarkable point, $\omega_0$, defined as the pair $(\Ll,B)$, so that both the actual magnetic field and the electric potential can be written in the form $V(x)=\vs(\tra^{-x}\omega_0)$ and $B_{ij}(x)=\Bb_{ij}(\tra^{-x}\omega_0)$. 
%
\begin{defini}
\label{cspa.def-hullLB}
Let $\Ll\subset \RM^d$ 
be a Delone set 
with Hull $\Omega_{\Ll}$. Let $B=(B_{ij})_{1\leq i<j\leq d}$ be a family of bounded, uniformly continuous functions defined on $\RM^d$ with Hull $\Omega_B$. Let $\omega_0\in \Omega_\Ll\times \Omega_B$ be the point representing the pair $(\Ll,B)$. Then\\
(i) the hull of this pair $\Omega$ is the closure, in the product space $\Omega_\Ll\times \Omega_B$, of the $\RM^d$-orbit of $\omega_0$.\\ 
(ii) the canonical transversal $\Xi$ of this pair is the closure of the set $\{\tra^{-y}\omega_0\,;\, y\in\Ll\}$.
\end{defini}
%

\begin{rem}
\label{cspa.rem-hullLB}
{\em It might happen that the Hulls of $\Ll$ and of $B$ ``coincide''. Namely that the two Hulls admit an homeomorphism conjugating their $\RM^d$-actions. In such a case, the Hull of $(\Ll,B)$ also coincides with each of these Hulls. For instance, when both $\Ll$ and $B$ 
    are periodic with the same period group $L$ (namely a discrete co-compact subgroup of $\RM^d$), their Hulls coincide both with $\RM^d/L$. As another example, if $B$ is weakly $\Ll$-pattern equivariant (see \cite{Kel03,Kel08}), the combined Hull coincides with $\Omega_\Ll$.  
}
\hfill $\Box$
\end{rem}

\begin{rem}
\label{cspa.rem-tranLB}
{\em That $\Xi$ is a transversal comes from the fact that $\Xi_\Ll$ is a transversal in $\Omega_\Ll$. Indeed, since $\Ll$ is a Delone set, there is $r>0$ such that any open ball of radius $r$ in $\RM^d$ contains at most one point of $\Ll$. The same property is shared by any element of $\Omega_\Ll$. Consequently if $\Ll'\in \Xi_{\Ll'}$, then $\Ll'$ contains the origin of $\RM^d$, which is the only point in $B(0;r)\cap \Ll$. Hence if $y\in\RM^d$ satisfies $y\neq 0$ and $\Ll'-y\in\Xi_{\Ll}$, namely if $y\in\Ll'$ then $|y|\geq r$. Consequently, the same property will occur for the combined transversal $\Xi$. 
}
\hfill $\Box$
\end{rem}

\begin{rem}
\label{cspa.rem-flux}
{\em Since the combined hull of the pair $(\Ll,B)$ involves the hull of $B$ itself, the magnetic flux through any loop in $\RM^d$ can also be viewed as a continuous function on $\Omega$. However, the vector potential is only defined on each $\RM^d$-orbit, through a Poincar\'e construction. Hence there will be an advantage in using the magnetic field instead of the vector potential, when constructing the \CS.  
}
\hfill $\Box$
\end{rem}

\noindent The underlying groupoid will then be $\Omega\rtimes_{\tri} \RM^d$ here. It ought to be reminded here that a typical arrow in this groupoid has the form of a pair $\gamma=(\omega,x)$ with $\omega\in \Omega$ and $x\in\RM^d$, and that $r(\gamma)=\omega$ while $s(\gamma)=\tra^{-x}\omega$. Then the magnetic translation operator becomes $\omega$-dependent, through the left regular representation of this groupoid. Similarly the cocycle defining the \Cs is given now by

$$\sigma((\omega, x),(\tra^{-x}\omega,y)) =
e^{\imath e\Phi_\omega(0,x,y)/\hbar}\,,
$$

\noindent where $\Phi_\omega(0,x,y)$ represent the magnetic flux through the oriented triangle $(0,x,y)$ with vertices $(0,x,x+y)$ defined by

$$\Phi_\omega(x,y) = \int_{(0,x,y)}
   \sum_{i,=1}^d \Bb_{ij}(\tra^{-s}\omega)ds_i\wedge ds_j\,.
$$

  \subsubsection{The discrete case}
  \label{cspa.sssect-discNU}
\noindent A similar situation occurs for the discrete version of the Schr\"odinger operator. Given the Delone set $\Ll$ and a point $\omega$ in the combined hull $\Omega$, there is $\Ll_\omega$ another Delone set such that $\Ll_\omega+x=\Ll_{\tra^x\omega}$. The canonical transversal is the closed subset $\Xi\subset \Omega$ of points such that $\Ll_\omega$ contains the origin $0\in\RM^d$. Let $\Gamma_\Xi$ be the groupoid associated with $\Xi$. Then, the Hilbert space of quantum states will be $\Hh_\xi=\ell^2(\Ll_\xi)$ for $\xi\in\Xi$. The translation allowed are those moving an element of the transversal onto another one. Such translations are represented by the arrow of the groupoid $\Gamma_\Xi$. Namely, as explained in \cite{Co79}, each arrow $(\xi,a)\in\Gamma_\Xi$ gives a unitary operators $T(\xi,a): \Hh_{\tri^{-a}\xi}\to \Hh_\xi$ defined by 
$$T(\xi,a)\psi(x)= \psi(x-a)\,,
   \hspace{2cm}
    \psi\in\Hh_{\tri^{-a}\xi}\,,\, x\in\Ll_\xi\,, \xi\in\Xi\,.
$$
\noindent As in the continuous case, and as it is explained in Remark~\ref{cspa.rem-PotVec}, the translation operators are modified by the magnetic field. Following the prescription proposed in eq.~\eqref{cspa.eq-magtr}, the corresponding magnetic translations gives a new unitary operator $U(\xi,a): \Hh_{\tri^{-a}\xi}\to \Hh_\xi$ defined by
$$U(\xi,a)\psi(x)= 
   e^{\imath(e/\hbar)\oint_{x-a}^x A_\xi}\,\psi(x-a)\,,
   \hspace{2cm}
    \psi\in\Hh_{\tri^{-a}\xi}\,,\, x\in\Ll_\xi\,, \xi\in\Xi\,.
$$
\noindent An elementary calculation shows that the eq.~\eqref{cspa.eq-comrul} becomes now
$$U(\xi,a)U(\tra^{-a}\xi,b)U(\tra^{-a-b}\xi, 0)=
   e^{\imath(e/\hbar)\Phi_\xi(0,a,b)}\,.
$$
\noindent A typical covariant Hamiltonian acting on $\Hh_\xi$ will be given by a combination of such magnetic translations with operators of multiplication as coefficients, namely an expression of the form
$$H_\xi \psi(a) = \sum_{b\in\Ll_\xi}
   h(\tra^{-a}\xi,b) 
    e^{\imath(e/\hbar)\oint_{b}^a A_\xi}\,
     \psi(b)\,,
      \hspace{2cm}
       a\in \Ll_\xi\,,\; \psi\in \ell^2(\Ll_\xi)\,,
$$
\noindent where (i) $h$ is a continuous function on $\Gamma_\Xi$ decaying fast enough at infinity to insure that the sum converges and (ii) $A_\xi$ is a vector potential associated with the magnetic field $B_\xi(x)=\Bb(\tra^{-x}\xi)$. So the treatment will be similar to the continuum case, provided one pays attention to use the groupoid $\Gamma_\Xi$ as a substitute for the translation group. In particular, the representation described in this section will follow the same prescription as the one provided in the next Sections.

 \subsection{Convolution and the Full $C^\ast$-algebra}
 \label{cspa.ssect-convol}

\noindent Given a locally compact space $Y$, let $\Cc_c(Y)$ denote the space of complex valued continuous functions defined on $Y$ with compact support. This is a topological complex vector space. The topology is defined so that a net $(f_\alpha)_{\alpha\in A}$ converges to $f$ if and only if (i) there is $K\subseteq Y$ compact and $\alpha\in A$, such that for $\beta\geq \alpha$ $\supp(f_\beta)\subseteq K$, and (ii) for any $\epsilon >0$ there is $\alpha\leq \alpha '\in A$ such that $\sup_{y\in K}|f_\beta(y)-f(y)|<\epsilon$ for $\beta \geq \alpha'$. 

\vspace{.1cm}

\noindent Let $\Gamma$ be a handy groupoid. Then $\Gamma$ admit a Haar system $\mu=(\mu^x)_{x\in\GaO}$. In addition, let $\sigma$ be a normalized $2$-cocycle with values in the unit circles $\SM^1$ seen as a multiplicative subgroup of $\CM_\ast$. Then a structure of $\ast$-algebra can be defined on the function space $\Cc_c(\Gamma)$ in complete analogy with the case of locally compact groups (see \cite{EH67,Di69}). If $f, g\in \Cc_c(\Gamma)$, their convolution is defined  by
\begin{equation}
\label{cspa.eq-Ccprod}
fg(\gamma) =
   \int_{\Gamma^x} \sigma(\eta, \eta^{-1}\circ \gamma)\;
    f(\eta)\, g(\eta^{-1}\circ\gamma) \,d\mu^x(\eta)\,,
    \hspace{2cm}
     r(\gamma)=x\,.
\end{equation}

\noindent Following \cite{Re80}, thanks to the cocycle property of $\sigma$, it is straightforward to check that this product is associative and bilinear. In addition, the function $fg$ is continuous thanks to (H2) and has compact support \cite{Se86}. An {\em adjoint} map is now defined as follows
\begin{equation}
\label{cspa.eq-Ccadj}
f^\ast(\gamma) =\overline{\sigma(\gamma,\gamma^{-1})}\;\;
   \overline{f(\gamma^{-1})}\,.
\end{equation}

\noindent A straightforward calculation shows that $f\to f^\ast$ is antilinear. Moreover, since $\sigma$ is normalized, Lemma~\ref{cspa.lem-1coc} implies that the adjoint map is involutive, namely $(f^\ast)^\ast=f$. It also satisfies 
$$
(f\,g)^\ast = g^\ast\,f^\ast\,.
$$

\noindent In addition, these two operations are continuous with respect to the topology of $\Cc_c(\Gamma)$ \cite{Re80}.

\begin{defini}
\label{cspa.def-algC-c}
The topological $\ast$-algebra obtained from $\Cc_c(\Gamma)$ by using the product (eq.~\eqref{cspa.eq-Ccprod}) and the adjoint (eq.~\eqref{cspa.eq-Ccadj}) will be denoted by $\Cc_c(\Gamma,\mu,\sigma)$.
\end{defini}

\noindent The topology of $\Cc_c(\Gamma,\mu,\sigma)$ is not coming from a norm. It is often considered as more  convenient to define algebraic norms and to complete this topological algebra in order to get a more practical tool for analysis. In this Section only two examples of such norms will be proposed. The first one is the extension of the Wiener algebra. For this purpose, let $\|f\|_{\infty,1}^\sim$ be defined by \cite{Ha78,Re80}
$$\|f\|_{\infty,1}^\sim=
   \sup_{x\in \GaO} \int_{\Gamma^x} |f(\gamma)|\,d\mu^x(\gamma)\,.
$$

\noindent It is straightforward to check that this defines a norm on $\Cc_c(\Gamma)$ such that
$$\|fg\|_{\infty,1}^\sim\leq 
   \|f\|_{\infty,1}^\sim\|g\|_{\infty,1}^\sim\,.
$$

\noindent In general though, this norm is not invariant under taking the adjoint. This leads to set
$$\|f\|_{\infty,1} =
   \max\{\|f\|_{\infty,1}^\sim,\|f^\ast\|_{\infty,1}^\sim\}
$$

\noindent The completion of $\Cc_c(\Gamma,\mu,\sigma)$ {\em w.r.t.} this norm will be denoted by $L^{\infty,1}(\Gamma, \mu)$. Using a $3\epsilon$-argument it is easy to check that any element of this normed $\ast$-algebra can be seen as a $\mu$-measurable function which is $L^1$ along each fiber, such that $x\in\GaO\to \mu^x(f)$ is continuous.

\vspace{.1cm}

\noindent In order to define a \CS, the usual procedure consists in looking at the representations of $\Cc_c(\Gamma,\mu,\sigma)$ by operators on a Hilbert space. Following \cite{Re80}, a {\em representation} of $\Cc_c(\Gamma,\mu,\sigma)$ on a Hilbert space $\Hh$ is a weakly continuous $\ast$-homomorphism $\rho: \Cc_c(\Gamma,\mu,\sigma)\to \Bb(\Hh)$, such that the linear span of the set $\{\rho(f)\psi\,;\, f\in \Cc_c(\Gamma,\mu,\sigma)\,,\, \psi\in \Hh\}$ is dense in $\Hh$. This representation is called {\em bounded} whenever $\|\rho(f)\|\leq \|f\|_{\infty,1}$. A $C^\ast$-norm can be defined by

\begin{equation}
\label{cspa.eq-fullnorm}
\|f\|=\sup_\rho \|\rho(f)\|\,,
   \hspace{2cm}
    \rho\;\mbox{\rm a bounded representation.}
\end{equation}

\begin{defini}
\label{cspa.def-full}
The full \Cs of $\Gamma$ is the completion of $\Cc_c(\Gamma,\mu,\sigma)$ under the norm $\|\cdot\|$ (eq.~\eqref{cspa.eq-fullnorm}). It will be denoted by $ \CG^\ast(\Gamma,\mu,\sigma)$ or $ \CG^\ast(\Gamma,\sigma)$ if there is no ambiguity on $\mu$.
\end{defini}


 \subsection{The Reduced $C^\ast$-algebra}
 \label{cspa.ssect-red}

\noindent Let $G$ be a locally compact group equipped with a left-invariant Haar measure $\lambda$. A {\em unitary representation} $U$ of $G$, in a Hilbert space $\Hh$, is a group homomorphism $U:G\to \Uu(\Hh)$, where $\Uu(\Hh)$ denotes the set of unitary operators on $\Hh$. Such a representation is {\em strongly continuous} whenever the map $g\in G\mapsto U(g)\in\Uu(\Hh)$ is strongly continuous. Since $U(g)$ is a unitary operator, strong continuity holds if and only if it is weakly continuous, that is, for any pair $\phi,\psi$ of elements of $\Hh$, the map $g\to \langle\phi|U(g)\psi\rangle$ is continuous. The {\em left-regular representation} is an important example. 
Let $\Hh=L^2(G)$ be the space of square integrable complex valued functions with respect to $\lambda$. Then $G$ acts on $\Hh$ by left multiplication, namely if $\psi\in L^2(G)$ and $g\in G$, then $L(g)\psi(h)=\psi(g^{-1} h)$ defines a strongly continuous unitary representation of $G$ in $\Hh$. Similarly it acts by right multiplication with $R(g)\psi(h)=\psi(hg)$, giving another strongly continuous unitary representation of $G$ in $\Hh$. It can be checked that $L(g)R(g')=R(g')L(g)$. If, in addition, $G$ is a compact group, the Peter-Weyl Theorem expresses that $L$ ({\em resp.} $R$) is a direct sum of irreducible unitary representations of $G$, each such representation occurring with multiplicity equal to its dimension. It is therefore legitimate to ask whether the left-regular ({\em resp.} the right-regular) representation is sufficient to recover entirely the group $G$. 

\vspace{.1cm}

\noindent To express this concept, let $ \CG^\ast_{red}(G)$ denotes the \Cs generated by the bounded operators acting on $\Hh$ of the form $R(f)=\int_G f(g)\, R(g) d\lambda(g)$ for $f\in \Cc_c(G)$, and let $\|\cdot\|_{red}$ be its norm. It is worth noticing that $L(g)R(f)L(g)^{-1}=R(f)$.
In particular, if $G$ is a compact group, $R(f)$ is the direct sum of representative in each irreducible component of $L$. If the right-regular representation suffices to generate all unitary representation, as suggested by the Peter-Weyl theorem, then any other strongly continuous unitary representation $V$ of $G$, on some Hilbert space $\Kk$, must satisfy $\|\int_G f(g)\,V(g)\,d\lambda(g)\|\leq \|f\|_{red}$. As it turns out, this is the case if and only if $G$ is {\em amenable} \cite{Hu66} (for the concept of amenability, see Section~\ref{cspa.ssect-amen} and \cite{Fo55,Gr69}). Any locally compact Abelian group, or any compact group, is amenable. On the other hand, the free group with $n\geq 2$ generators is not amenable and it does not satisfy the previous criterion. 

\vspace{.1cm}

\noindent A similar phenomenon occurs for groupoids, for which there is also a concept of amenability \cite{AR00,AR01}. To understand more precisely, the first concept to be defined is the definition of a strongly continuous unitary representation of a groupoid (analogous to the definition of  measurable representations \cite{Co79}). Since a groupoid $\Gamma$ is a {\em category}, a unitary representation will be a (covariant) functor from $\Gamma$ into the category {\bf Hilb} the objects of which are Hilbert spaces, with morphism given by unitary operators. Hence for each object (unit) $x\in\GaO$, a Hilbert space $\Hh_x$ should be given so that if $\gamma:x\to y$, then there is a unitary operator $U(\gamma):\Hh_x\to \Hh_y$. To be a covariant functor, the additional property should be that $U(\gamma_1\circ\gamma_2)= U(\gamma_1)U(\gamma_2)$. To express the strong continuity, it will be assumed that the field $(\Hh_x)_{x\in\GaO}$ of Hilbert space is {\em continuous} \cite{Di69}. Namely there are enough {\em continuous vector fields} $\psi=(\psi_x)_{x\in\GaO}$ to generate each of the Hilbert spaces $\Hh_x$. By definition of a continuous vector field, the map $x\in \GaO\mapsto \|\psi_x\|\in [0,\infty)$ is continuous. Hence, by analogy with the case of groups, strong continuity of $U$ can be expressed by demanding that for any pair $\phi,\psi$ of continuous vector fields, the map $\gamma\in \Gamma\mapsto \langle \phi_{r(\gamma)}|U(\gamma)\psi_{s(\gamma)}\rangle\in \CM$ is continuous.

\vspace{.1cm}

\noindent If $\sigma$ is a normalized $2$-cocycle, a $\sigma$-representation of $\Gamma$ will be defined as in the previous construction, with the following modifications
$$U(\gamma_1)U(\gamma_2)= 
   \sigma(\gamma_1,\gamma_2)\;U(\gamma_1\circ\gamma_2)\,,
    \hspace{2cm}
     U(\gamma)^\ast= 
      \overline{\sigma(\gamma,\gamma^{-1})}\;
        U(\gamma^{-1})\,.
$$

\vspace{.1cm}

\noindent The simplest example is provided by {\em left regular representation}. Namely, given a continuous Haar system on $\Gamma$, then $\Hh_x=L^2(\Gamma^x,\mu^x)$. If $\gamma:x\to y$, the unitary operator $U(\gamma)$ is now replaced by
$$L(\gamma)\psi(\eta)= \sigma(\gamma, \gamma^{-1}\circ \eta)\;
   \psi(\gamma^{-1}\circ \eta)\,,
    \hspace{2cm}
     \psi\in L^2(\Gamma^x,\mu^x)\,.
$$

\noindent Thanks to the Assumption (H3), it follows that $L(\gamma)$ is a unitary operator. It becomes straightforward to check that $L(\gamma_1) \,L(\gamma_2)= \sigma(\gamma_1,\gamma_2)\;L(\gamma_1\circ\gamma_2)$. It is worth noticing that the subspace of $L^2(\Gamma^x,\mu^x)$ generated by the restriction to $\Gamma^x$ of elements of $\Cc_c(\Gamma)$ is dense. 

\vspace{.1cm}

\noindent The analog of the operators $R(f)$ acting on $L^1(G)$ in the group case, can be also defined in the groupoid case. Indeed, if $f\in\Cc_c(\Gamma,\mu,\sigma)$, an operator can be defined on $L^2(\Gamma^x,\mu^x)$ as follows
\begin{equation}
\label{cspa.eq-redrep}
\pi_x(f)\psi(\gamma)=
   \int_{\Gamma^x} f(\gamma^{-1}\circ\eta)\,
    \sigma(\gamma, \gamma^{-1}\circ \eta)\;
     \psi(\eta)\,d\mu^x(\eta)\,,
    \hspace{2cm}
     \psi\in L^2(\Gamma^x,\mu^x)\,.
\end{equation}

\noindent Standard estimates imply 
$$\|\pi_x(f)\|\leq \|f\|_{\infty,1}\,.
$$

\noindent showing that $\pi_x(f)$ is a bounded operator. A tedious but simple computation shows that $\pi_x$ is linear in $f$, that
$$\pi_x(fg) =
   \pi_x(f)\;\pi_x(g)\,,
    \hspace{2cm}
     \pi_x(f)^\ast=\pi_x(f^\ast)\,,
      \hspace{2cm}
       f,g\in\Cc_c(\Gamma)\,.
$$

\noindent Therefore $\pi_x$ is a bounded representation of $\Cc_c(\Gamma,\mu,\sigma)$. In addition, using the density of $\Cc_c(\Gamma,\mu,\sigma)$ in $L^2(\Gamma^x,\mu^x)$, the Assumption (H2) and a $3\epsilon$-argument, it follows that $x\to \pi_x(f)$ is strongly continuous. At last, if $\gamma:x\to y$ then, the following {\em covariance} property holds
$$U(\gamma)\pi_x(f) U(\gamma)^{-1}=
   \pi_y(f)\,.
$$

\noindent It then follows that $x\sim y$ implies $\|\pi_x(f)\|=\|\pi_y(f)\|$. Moreover, whenever $f=f^\ast$ is self-adjoint, then $\pi_x(f)$ and $\pi_y(f)$ have the same spectral measure and the same spectrum. It is also straightforward to check that $\pi_x(f)=0$ if and only if $f$ vanishes on the fiber $\Gamma^x$. Consequently $\pi_x(f)=0$ for all $x\in\GaO$ if and only if $f=0$. This implies that
$$\|f\|_{red}=
   \sup_{x\in\GaO} \|\pi_x(f)\|\,,
$$

\noindent defines a norm. By construction this norm is algebraic, invariant by the adjoint and satisfies $\|f^\ast f\|_{red}=\|f\|_{red}^2$. Hence it is a $C^\ast$-norm. 

\begin{defini}
\label{cspa.def-red}
Let $\Gamma$ be a handy groupoid endowed with a Haar system $\mu$ and a normalized $2$-cocycle $\sigma$. Then the reduced \Cs of $\Gamma$ is the completion of $\Cc_c(\Gamma,\mu,\sigma)$ under the $C^\ast$-norm $\|\cdot\|_{red}$. This completion will be denoted by $\CG_{red}^\ast(\Gamma, \mu,\sigma)$.
\end{defini}

\noindent In general the left regular representation is not sufficient to encode all representations. A result of Anantharaman-Delaroche and Renault \cite{AR00} shows that, at least for handy \'etale groupoids, amenability implies $\CG^\ast_{red}(\Gamma,\mu)= \CG^\ast(\Gamma,\mu)$ (Theorem 6.1.8). In \cite{BH14}, Corollary~4.3, the analog result is proved with a normalized $2$-cocycle, namely $\CG^\ast_{red}(\Gamma,\mu,\sigma)= \CG^\ast(\Gamma,\mu,\sigma)$ if $\Gamma$ is amenable. However, even in the \'etale case, R.~Willett \cite{Wil15} gave an example of handy \'etale groupoid that is not amenable and for which $\CG^\ast_{red}(\Gamma,\mu)=\CG^\ast(\Gamma,\mu)$.

 \subsection{Amenability}
 \label{cspa.ssect-amen}

\noindent At this point it might help the reader to give some hint about the concept of amenability. The most elementary way to address this problem is as follows: consider a {\em time-dependent} function $f$. Here the time is represented by a variable $t\in\RM$. The average of $f$ over time is simply defined as 
$$\langle f\rangle =
   \lim_{T\to\infty}
    \frac{1}{2T}\int_{-T}^{+T} f(t)dt\,.
$$

\noindent If $f$ is unbounded such a limit might not even exist. But even if $f$ is bounded and continuous, {\em there might be several limits}. In this precise example, all possible limits are classified by ultrafilters, or, equivalently by the \v{C}ech compactification of $\RM$. Then any of these limits are {\em translation-invariant}, namely if $s\in \RM$ and if $f_s(t):=f(t-s)$, then $\langle f_s\rangle=\langle f\rangle$ for all $s$'s. This is because, if $\|f\|=\sup_t|f(t)|$,
$$\left|\frac{1}{2T}\int_{-T}^{+T} f(t-s)dt-
   \frac{1}{2T}\int_{-T}^{+T} f(t)dt\right| \leq 
    \frac{|s|\,\|f\|}{T}\;\;
     \stackrel{T\uparrow\infty}{\longrightarrow} \;\;0\,.
$$

\noindent A similar problem might occur in higher dimension, if the time axis $\RM$ is replaced by the position space $\RM^d$, in order to average over space. This problem occurs in Thermodynamics for which an observable can be defined as a position dependent function, $f(x)\,,\, x\in \RM^d$. Then its space-average can be defined by
$$\langle f\rangle =
   \lim_{\Lambda\uparrow \RM^d}
    \frac{1}{|\Lambda|}\int_{\Lambda} f(x)d^dx\,,
$$

\noindent where $|\Lambda|$ denotes the volume of $\Lambda$, namely its Lebesgue measure. Again in this case, at least if $f$ is continuous and bounded, the \v{C}ech compactification of $\RM^d$ classifies all limits and it is possible to show that, under suitable conditions on the choice of $\Lambda$, these limits are translation invariant. Such a limit defines what is called an {\em invariant mean}. In dealing with translation invariance, it becomes easy to realize that the ratio of the area $|\partial\Lambda|$ of the boundary of $\Lambda$ should become negligible {\em w.r.t.} the volume of $\Lambda$ as this volume grows to become infinite. This concept is crucial in Thermodynamics or in Statistical Mechanics, where the limit above is controlled by the concept of {\em Van Hove} limit (see \cite{VH49} Section 2 Axiom (e)). 

\vspace{.1cm}

\noindent If $\RM^d$ is replaced by a locally compact group $G$, the same problem can be addressed in similar terms. A mean is defined as a state on the Abelian \Cs $\Cc_b(G)$ of complex valued, continuous, bounded functions on $G$. That is to say

(i) $m:\Cc_b(G)\to \CM$ is linear,

(ii) $m$ is positive, namely if $f\in \Cc_b(G)$ is non-negative everywhere in $G$, then $m(f)\geq 0$,

(iii) $m$ is normalized, namely $m(1)=1$.

\noindent Given $s\in G$ and $f\in \Cc_b(G)$, let $f_s$ be defined by $f_s(t)=f(s^{-1}t)$ for $t\in G$. Then a mean $m$ is $G$ left-invariant whenever $m(f_s)=m(f)$ for all $s\in G$.

\vspace{.1cm}

\noindent However, there are examples of groups for which there is no invariant mean~! Such is the case for the free group with $n\geq 2$ generators. In general, a locally compact group $G$ will be called {\em topologically amenable} whenever there is a $G$-invariant mean. Characterizing amenable groups has been a challenge for a long time. The concept of F{\o}lner sequence became central \cite{Fo55,Gr69}: it expresses the idea that in taking the average over a subset $\Lambda\subseteq G$, the boundary of $\Lambda$ has a Haar-measure negligible compared to the Haar measure of $\Lambda$. More precisely, if $\lambda$ is a left-invariant Haar measure on $G$, a sequence $(\Lambda_n)_{n\in\NM}$ of $\lambda$-measurable subsets of $G$ is called a F{\o}lner-sequence whenever 

(i) $\lambda(\Lambda_n)<\infty$ for all $n\in\NM$,

(ii) $\lambda(\Lambda_n)\to \infty$ as $n\uparrow \infty$,

(iii) for any $s\in G$, the symmetric difference $s\Lambda_n\triangle \Lambda_n$ satisfies
$$\lim_{n\to\infty}
   \frac{\lambda(s\Lambda_n\triangle \Lambda_n)}{\lambda(\Lambda_n)}=0\,.
$$

\begin{rem}
\label{cspa.rem-vanHove}
{\em A F{\o}lner sequence might escape to infinity, far from the ``origin'' (namely the neutral element) of the group. In particular the construction of a F{\o}lner sequence is less constrained than the construction of a Van~Hove sequence that is commonly used in Physics for the purpose of controlling the infinite volume limit \cite{VH49,VH55}.
}
\hfill $\Box$
\end{rem}

\noindent Given a F{\o}lner sequence and given an ultrafilter $\FG$ on the set $\NM$ of natural integers, let $m$ be the mean defined by
$$m(f) =\lim_{\FG} 
   \frac{1}{\lambda(\Lambda_n)}
    \int_{\Lambda_n} f(t) d\lambda(t)\,,
   \hspace{2cm}
    f\in \Cc_b(G)\,.
$$

\noindent That such a limit is well defined comes from the fact that sequence on the {\em r.h.s.} is bounded by $\|f\|$. Moreover, 
$$|m(f_s)-m(f)|\leq 
   \lim_{\FG}\;\|f\|
    \frac{\lambda(s\Lambda_n\triangle \Lambda_n)}{\lambda(\Lambda_n)}
     =0\,,
$$ 

\noindent showing that this mean is $G$-invariant. Hence, the existence of a F{\o}lner sequence implies the existence of an invariant mean. It turns out that the converse is true \cite{Gr69}.

\vspace{.1cm}

\noindent Eventually, at least for locally compact Hausdorff and second countable groups, amenability could be described in terms of the equality $\CG^\ast_{red}(G)=\CG^\ast(G)$ \cite{Hu66}. For a handy groupoid $\Gamma$, the concept of amenability is much more involved. Anantharaman-Delaroche and Renault expressed it in terms of positive definite functions in \cite{AR00,AR01}. They distinguish between several concepts of amenability, topological, Borel or measurable. 

\vspace{.1cm}

\noindent As a consequence, a dynamical system $(X,G,\alpha)$ where $G$ is a locally compact, Hausdorff, second countable and amenable group and $X$ a compact Hausdorff, second countable space, gives rise to an amenable groupoid $\Gamma= X\rtimes_\alpha G$. If $Y\subseteq X$ is closed, any sub-groupoid $\Gamma_Y\subseteq X\rtimes_\alpha G$ is also amenable. And this is the only property that will be needed in the present work. However, there are dynamical systems for which $G$ is not amenable but the action on $X$ admits an invariant probability measure, leading to an amenable groupoid. The typical example is the action of the free group with 2-parameters acting on its boundary (see \cite{AR00}, Ex. 3.8). 

\section{The Tautological Groupoid}
\label{cspa.sect-taut}

\noindent In all this Section $\Gamma$ denotes a handy groupoid. The set of closed invariant sets $\isG$ of a handy groupoid $\Gamma$ is introduced. Among these subsets, some are minimals, like in the case of dynamical systems. In many examples of applications to Physics, minimal {\em finite} invariant subsets can be identified with periodic orbits. Using the Hausdorff topology (as defined either by Vietoris \cite{Vi22}, by Chabauty \cite{Ch50} or by Fell \cite{Fe62}), the space $\isG$ becomes a compact Hausdorff second countable space. This leads to the definition of the {\em tautological groupoid} of $\Gamma$, denoted by $\ts(\Gamma)$. Consequently, at least whenever $\Gamma$ is \'etale and amenable, thanks to a Theorem by Landsman and Ramazan \cite{LR01}, the \Cs 
(without a cocycle)
of the tautological groupoid appears as the enveloping algebra of a continuous field of \Css indexed by $\isG$. 
This machinery permits to control whether or not a closed invariant subset can be approximated by a sequence of periodic orbits. If so, a genuine Schr\"odinger operator $H$, defined on the closed invariant subset, can be approximated by periodic operators for which Bloch Theory applies, leading to a convergent sequence of approximations of the spectrum of $H$.

 \subsection{Invariant Sets}
 \label{cspa.ssect-inv}

\noindent If $M\subseteq \GaO$, its saturated is defined by $[M]=\{x\in\GaO\,;\, \exists y\in M,\, x\sim y\}$. $M$ is called invariant whenever $M=[M]$. Equivalently, $M$ is invariant if and only if, whenever $x\in M$ and $y\sim x$ then $y\in M$. Hence an invariant set is the union of the equivalence classes of its elements. If $x\in \GaO$, $[x]$ is its equivalence class, namely it is the saturated of $\{x\}$. If $\Gamma$ is provided by a dynamical system $(X,G,\alpha)$ (see Example~\ref{cspa.exam-topdyn}), and if $x\in X=\GaO$, then $[x]$ is nothing but the set $\{g^{-1}x\,;\, g\in G\}$, namely it is the {\em G-orbit} of $x$. By analogy, $[x]$ will be also called the $\Gamma$-orbit of $x$  (or simply ``orbit'' if there is no ambiguity) and the notation $[x]=Orb(x)$ will be used. The first result is the following

\begin{proposi}
\label{cspa.prop-clinv}
Let $M\subseteq \GaO$ be invariant. Then its closure $\overline{M}$ is also invariant.
\end{proposi}

\noindent {\bf Proof: }Let $x\in \overline{M}$. For $y\sim x$ let $\gamma:y\to x$. Let then $V$ be an open neighborhood of $y$. Since the source map is continuous, it follows that $s^{-1}(V)$ is an open set in $\Gamma$ containing $\gamma$. Since the range map is open, $U=r(s^{-1}(V))$ is also open in $\GaO$ and contains $r(\gamma)=x$. Since $x\in \overline{M}$ the intersection $U\cap M$ is not empty. Let $x'\in U\cap M$. Then there is $\gamma'\in s^{-1}(V)$ such that $r(\gamma')=x'$, so that $s(\gamma')\in V$. Since $x'\in M$ and $M$ is invariant $s(\gamma')=y'\in M$. Hence $V\cap M\neq \emptyset$. Since $V$ is arbitrary, it follows that $y\in\overline{M}$.
\hfill $\Box$

\begin{defini}
\label{cspa.def-clinv}
The set of closed invariant subset in $\GaO$ is denoted by $\isG$. 
\end{defini}

\noindent Using the analogy with topological dynamical systems \cite{GH55}, the following result holds

\begin{theo}
\label{cspa.th-min}
The set $\isG$, ordered by inclusion, has minimal elements. A closed invariant subset $M\in \isG$ is minimal if and only if any orbit is dense in $M$.
\end{theo}

\noindent  {\bf Proof: }(i) Any intersection of closed invariant subsets $(M_j)_{j\in J}$ of $\GaO$ is closed and invariant: (a) it is closed since any intersection of closed sets is closed, (b) it is invariant because if $N$ denotes the intersection, if $x\in N$ then, for all $j\in J$, $x\in M_j$. Hence if $y\sim x$, it follows that $y\in M_j$ as well for all $j$'s. Hence $y\in N$ and $N$ is closed and invariant. In particular, the intersection of a decreasing family of elements of $\isG$ belongs to $\isG$. By Zorn's Lemma, $\isG$ admits minimal elements.

\vspace{.1cm}

\noindent (ii) Let $M$ be minimal. If $x\in M$ then $[x]\subseteq M$ since $M$ is invariant. In addition the closure $\overline{[x]}$ is also included in $M$ since $M$ is closed, and is invariant by Proposition~\ref{cspa.prop-clinv}. Since $M$ is minimal $\overline{[x]}=M$.

\vspace{.1cm}

\noindent (iii) Conversely, let $M\in \isG$ be such that every of its orbit is dense. If $N\subseteq M$ is closed and invariant, then given $x\in N$ it follows that $M=\overline{[x]}\subseteq N$, so that $M=N$ and $M$ is minimal.
\hfill $\Box$

\vspace{.2cm}

\noindent If $(X,G,\alpha)$ is a topological dynamical system giving rise to a handy groupoid, then for $x\in X$, $[x]$ is the $G$-orbit of $x$. In particular, $G$ acts {\em transitively} on $[x]$. Let $H_x$ be the {\em stabilizer} of $x$ in $G$, namely $H_x=\{h\in G\,;\, h^{-1}x=x\}$. It follows that the $g^{-1}x$ depends only on the right coset $[g]=H_x g$, so that $H_x g\in H_x\backslash G\mapsto g^{-1}x\in [x]$ is a bijection. In addition, it is continuous. The following result, initially proved by Freudenthal \cite{Fr36} gives a condition under which this map is an homeomorphism.

\begin{theo}[Open Action \cite{Fr36}]
\label{cspa.th-open}
Let $G$ be a group acting on a space $X$ in a transitive way (the orbit of any point is $X$). If $G$ and $X$ are locally compact and if $G$ is a Lindel\"of space, then the map $\phi_x: G\to X$, defined by $\phi_x(g) = g^{-1}x$ is open.
\end{theo}

\noindent  It is worth reminding that a topological space is Lindel\"of whenever every open covers contains a countable sub-cover. In particular all compact space are Lindel\"of. The proof below is provided because it will be used later in the case of groupoid.

\vspace{.1cm}

\noindent {\bf Sketch of the Proof of Theorem~\ref{cspa.th-open}: } Let $K\subseteq G$ be any compact symmetric (namely $K$ is invariant by $g\mapsto g^{-1}$) neighborhood of the identity in $G$. Then $G$ can be seen as the countable union of sets of the form $g^{-1}K$. Let $\stackrel{\circ}{A}$ denote the interior of a subset $A\subseteq X$. By Baire's theorem, some of the $\phi_x(g^{-1}K)$ has a closure with non empty interior. By transitivity of the action, there is $h\in K$ such that $\phi_x(h)\in \stackrel{\circ}{\overline{\phi_x(K)}}$. In particular $x=\phi_x(e)\in \stackrel{\circ}{\overline{\phi_x(h^{-1}K)}}\subseteq \stackrel{\circ}{\overline{\phi_x(K^2)}}$. Let now $U\subseteq G$ be open and let $g\in U$. Then $K$ can be chosen small enough so that $gK^2\subseteq U$. Therefore $\phi_x(g)\in \stackrel{\circ}{\overline{\phi_x(gK^2)}} \subseteq \phi_x(U)$. Hence $\phi_x(U)$ is open. 
\hfill $\Box$

\begin{coro}
\label{cspa.cor-dynper}
Let $(X,G,\alpha)$ be a topological dynamical system  giving rise to a handy groupoid. Let $x\in X$ have a closed $G$-orbit. If $H_x\subseteq G$ is the stabilizer of $x$, then $[x]$ is homeomorphic to $H_x\backslash G$ and the groupoid $\Gamma_{[x]}$ is isomorphic to the dynamical system $(H_x\backslash G,G, \beta)$ where $\beta$ is the canonical right $G$-action on $H_x\backslash G$. 
\end{coro}

\noindent  {\bf Proof: }If $[x]$ is closed, then it is compact (since $X$ is compact). In particular the condition that the groupoid $X\rtimes_{\alpha} G$ is handy, implies that $G$ is Lindel\"of. Of course so is $[x]$ then. The map $g\in G\to g^{-1}x\in [x]$ being open, it follows from the definition of the quotient topology that the induced map $[g]\in H_x\backslash G\to [x]$ is also open. Hence the inverse of this map is also continuous, meaning that this map is actually an homeomorphism. Consequently, $H_x\backslash G$ is compact, namely $H_x$ is co-compact in $G$. The last statement is an exercise left to the reader.
\hfill $\Box$

\vspace{.2cm}

\noindent  In the specific case for which $G=\RM^d$, a co-compact subgroup is isomorphic to $H=\ZM^n\times \RM^{d-n}$ for some $0\leq n\leq d$. The special case $n=d$ corresponds to a non-singular {\em periodic dynamical system}. So, it is natural to study points $x\in\GaO$ with closed $\Gamma$-orbit. If so $[x]$ is automatically minimal in $\isG$. The following definition and the arguments after it are a generalization of Corollary~\ref{cspa.cor-dynper} to the groupoid case:

\begin{defini}
\label{cspa.def-qper}
A unit $x\in \GaO$ is called periodicoid if its $\Gamma$-orbit is closed. 
\end{defini}

\noindent Let $x$ be periodicoid unit of $\Gamma$. We want to know what is the structure of the sub-groupoid $\Gamma_{[x]}$. In order to do so, let $\Gamma^x$ be the $r$-fiber of $x$ and let $H_x=\Gamma_x^x$ its $rs$-fiber. Then $\Gamma^x$ and $H_x$ are both locally compact, Hausdorff and second countable. In particular they are Lindel\"of. $H_x$ acts on $\Gamma^x$ on the left, namely if $\eta\in H_x$ then $\gamma\to \eta\circ\gamma$ gives this action. Similarly, $H_x$ acts diagonally on the product space $\Gamma^x\times \Gamma^x$ namely through $(\gamma,\gamma')\mapsto (\eta\circ \gamma,\eta\circ \gamma')$. Let then $S\Gamma^x=H_x\backslash (\Gamma^x\times \Gamma^x)$ be the set of diagonal left $H_x$-cosets of $\Gamma^x\times \Gamma^x$. The equivalence class of a pair $(\gamma, \gamma')\in \Gamma^x\times \Gamma^x$ will be denoted by $[(\gamma, \gamma')]$.

\vspace{.1cm}

\noindent A groupoid structure can be defined on $S\Gamma^x$ as follows:

(i) {\em Range map:} $r[(\gamma,\gamma')]=s(\gamma)$,

(ii) {\em Source map:} $s[(\gamma,\gamma')]=s(\gamma')$

\noindent Since $s(\eta\circ\gamma)=s(\gamma)$ for all $\eta\in H_x$, these definitions make sense in the coset space. Then $[(\gamma_1,\gamma_1')]$ and $[(\gamma_2,\gamma_2')]$ are composable if and only if $s(\gamma_1')=s(\gamma_2)$, namely if and only if $\eta=\gamma_2\circ \gamma_1'^{-1}\in H_x$. Hence, $\gamma_2=\eta\circ \gamma_1'$. Therefore there is a unique $\gamma_3\in \Gamma^x$ such that $(\gamma_2,\gamma_2')=(\eta\circ \gamma_1',\eta\circ \gamma_3)$. This gives the product rule

(iii) {\em Product: } $[(\gamma_1,\gamma_1')]\circ [(\eta\circ \gamma_1',\eta\circ \gamma_3)]=[(\gamma_1,\gamma_3)]$,

(iv) {\em Inverse: } $[(\gamma,\gamma')]^{-1}=[(\gamma',\gamma)]$

(v) {\em Unit: } if $y\in [x]$ and $\gamma:y\to x$ then $e_y=[(\gamma,\gamma)]$

\noindent If $\gamma,\gamma'\in \Gamma^x$ it follows that $\gamma^{-1}\circ \gamma'$ is an element of $\Gamma_{[x]}$. This is because 
$r(\gamma^{-1}\circ \gamma')=s(\gamma)\sim x$ and $s(\gamma^{-1}\circ \gamma')=s(\gamma')\sim x$ as well.

\begin{proposi}
\label{cspa.prop-gammx}
Let $\Gamma$ be a handy groupoid. Let $x\in\GaO$ be a periodicoid element. Then there is a topological groupoid isomorphism $\phi: S\Gamma^x\to \Gamma_{[x]}$ defined by
$$\phi\left(
    [(\gamma,\gamma')]
      \right)= \gamma^{-1}\circ \gamma'\,,
$$ 
\end{proposi}

\noindent {\bf Proof: }(i) Let $\varphi: \Gamma^x\times \Gamma^x \to \Gamma_{[x]}$ be defined by $\varphi(\gamma,\gamma')= \gamma^{-1}\circ \gamma'$. Clearly $\varphi$ is continuous. In addition, $\varphi (\gamma,\gamma') = \varphi(\eta\circ\gamma,\eta\circ \gamma')$ for all $\eta \in H_x$. Hence $\varphi$ defines a map $\phi: S\Gamma^x\to \Gamma_{[x]}$. Using the quotient topology implies that $\phi$ is continuous as well. 

\vspace{.1cm}

\noindent (ii) The map $\phi$ is onto: if $\delta:z\to y$ belongs to $\Gamma_{[x]}$, then there are $\gamma:y\to x$ and $\gamma':z\to x$ so that $\eta= \gamma\circ \delta\circ \gamma'^{-1}\in H_x$. Therefore $\delta= \gamma^{-1}\circ \eta\circ \gamma'= \varphi(\gamma, \eta\circ \gamma')$.

\vspace{.1cm}

\noindent (iii) The map $\phi$ is one-to-one. This is because $\gamma_1^{-1}\circ \gamma_1'=\gamma_2^{-1}\circ \gamma_2'$, if and only if $\gamma_2'\circ \gamma_1'^{-1} =\gamma_2\circ \gamma_1^{-1} =\eta\in H_x$, namely if and only if $(\gamma_2,\gamma_2')=(\eta\circ \gamma_1,\eta\circ \gamma_1')$.

\vspace{.1cm}

\noindent (iv) Under the composition $\phi$ satisfies
%
\begin{eqnarray*}
\phi\left(
   [(\gamma_1,\gamma_3)]
      \right) &=& 
    \gamma_1^{-1}\circ \gamma_3\\
&=& \gamma_1^{-1}\circ \gamma_1'\circ 
      \gamma_1'^{-1}\circ \eta^{-1}\circ \eta\circ \gamma_3\\
&=& \phi[(\gamma_1,\gamma_1')]\circ 
        \phi[(\eta\circ \gamma_1',\eta\circ \gamma_3)]\,.
\end{eqnarray*}

\noindent Showing that $\phi$ is a groupoid isomorphism.

\vspace{.1cm}

\noindent (v) The map $\phi$ is open: if $U\subseteq S\Gamma^x$ is open, then, by the very definition of the quotient topology, there is $V\subseteq \Gamma^x\times \Gamma^x$ open and $H_x$-invariant such that $[V]=U$. Hence, because the topology on $\Gamma^x\times \Gamma^x$ is induced by the topology on $\Pp\Gamma=\PGa$, there is $W\subseteq\Pp\Gamma$ which is open and satisfies $W\cap \Gamma^x\times \Gamma^x=V$. Now, let $\hr$ be the map defined on $\Pp\Gamma$ by $\hr(\gamma_1,\gamma_2)= r(\gamma_1)=r(\gamma_2)$. It is clearly continuous. Thanks to Lemma~\ref{cspa.lem-mopen}, it follows that $\hr$ is also open. Consequently $\hW=\hr^{-1}\left(\hr(W)\right)$ is open. It follows also, from the proof of Lemma~\ref{cspa.lem-mopen}, that the map $\hm:\Pp\Gamma\to \Gamma$ defined by $\hm(\gamma_1,\gamma_2)=\gamma_1^{-1}\circ\gamma_2$ is continuous and open. Consequently $\hm(\hW)$ is open in $\Gamma$. It should be remarked, at this point, that $\varphi$ coincides with the restriction of $\hm$ on $\Gamma^x\times \Gamma^x$. Thus $\varphi(\hW\cap \Gamma^x\times \Gamma^x)$ is obviously included in $\hm(\hW)\cap \Gamma_{[x]}$. This inclusion is actually an equality: this is because, 
if $\gamma\in \hm(\hW)$, there is a pair $(\gamma_1,\gamma_2)\in \hW$, with $\gamma_1^{-1}\circ\gamma_2=\gamma$. Hence $u=\hr(\gamma_1,\gamma_2)\in [x]$. By the transitivity property, it follows that there is $\delta:u\to x$ such that replacing $\gamma_i$ by $\delta\circ\gamma_i$ $\hm(\delta\circ\gamma_1,\delta\circ\gamma_2)=\gamma$. Since $\hW$ is $\hr$-invariant, it follows that $(\delta\circ\gamma_1,\delta\circ\gamma_2)\in\hW\cap \Gamma^x\times \Gamma^x$ as well. Therefore $\phi(U)=\varphi(\hW\cap \Gamma^x\times \Gamma^x)=\hm(\hW)\cap \Gamma_{[x]}$, proving that $\phi(U)$ is open. Since this conclusion holds for any open set $U\subseteq S\Gamma^x$, $\phi$ is open. Hence, $\phi $ is a homeomorphism. 
\hfill $\Box$

\begin{proposi}
\label{cspa.prop-etalper}
Let $\Gamma$ be a handy \'etale groupoid. Then a periodicoid element  $x\in \GaO$ has a finite orbit. Conversely any finite invariant subset is a disjoint union of periodicoid orbits. 
\end{proposi}

\noindent  {\bf Proof: }(i) Let $X, Y$ be topological spaces and let $\phi:X\to Y$ be a surjective open map. If $X$ is discrete, then $\{x\}$ is open for any $x\in X$. It follows that $\phi(\{x\})=\{\phi(x)\}$ is also open. Since $\phi$ is onto, it implies that $Y$ is discrete as well.
On the other hand, if $X$ is also compact the family $\Uu= \{\{x\}\,;\, x\in X\}$ is an open cover. Since $X$ is compact one can extract a finite subcover. But $\Uu$ is also a partition, therefore it must be finite, namely $X$ is finite.

\vspace{.1cm}

\noindent (ii) Let $x\in\GaO$ be periodicoid, namely the equivalence class $[x]$ is closed. Since $\Gamma$ is handy, $\GaO$ is compact, and therefore $[x]$ is compact as well. Thanks to Proposition~\ref{cspa.prop-gammx}, the topological groupoid $\Gamma_{[x]}$ is isomorphic to $S\Gamma^x$. Since $\Gamma$ is \'etale, the $r$-fiber $\Gamma^x$ is discrete, so that $\Gamma^x\times \Gamma^x$ is discrete as well. The quotient map $\pi: \Gamma^x\times \Gamma^x\to S\Gamma^x$ is always open by definition of the quotient topology. It follows that $S\Gamma^x$ is discrete, and thus $\Gamma_{[x]}$ is discrete. Since $\Gamma$ is handy, the source and the range maps are open. Hence $[x]=r\left(\Gamma_{[x]}\right) = s\left(\Gamma_{[x]} \right)$ is also discrete. Since $[x]$ is compact, it must be finite.

\vspace{.1cm}

\noindent (iii) Let $M\subseteq \GaO$ be finite and invariant. Then it is closed, since $\GaO$ is Hausdorff. The equivalence classes $\{ [x]\,;\, x\in M\}$ make up a finite partition of $M$. By construction each $x\in M$ is periodicoid.
\hfill $\Box$

\begin{rem}
\label{cspa.rem-proper}
{\em If the groupoid $\Gamma$ is proper, in addition to being handy, namely if the source and range maps are proper, then every unit is periodicoid (see \cite{Cl07}).
}
\hfill $\Box$
\end{rem}

 \subsection{The Hausdorff Topology}
 \label{cspa.ssect-haus}

\noindent In his famous 1914 memoir \cite{Ha14,Ha35,Ha62}, Hausdorff defined the foundations of topology. He also defined a topology on the set of closed subsets of a complete metric space $(X,d)$ using the so-called {\em Hausdorff metric}, defined as follows: if $A,B$ are two subset of $X$, then 

$$\delta(A,B)=\sup_{x\in A} \dist(x,B)\,,
   \hspace{2cm}
    d_H(A,B)= \max\{\delta(A,B),\delta(B,A)\}\,.
$$

\noindent By definition $d_H(A,B)=d_H(B,A)$. In addition, it is straightforward to check that $d_H(A,B)\leq d_H(A,C)+d_H(C,B)$ for any $C$. In general $d_H(A,B)$ may not be finite. In addition $d_H(A,B) =0$ if and only if $A$ and $B$ have the same closure. But if both $A,B$ are compact, $d_H$ is finite and defines a metric on the space $\ks(X)$ of non empty compact subsets of $X$. If $X$ is complete, so is $\ks(X)$. If $X$ is compact, so is $\ks(X)$ (see \cite{Mi51,Mu75,Ba88} for instance).

\vspace{.1cm}

\noindent In 1922, Vietoris \cite{Vi22} published a paper in which he defined a topology on the set $\cs(X)$ of nonempty closed subsets of $X$, whenever $X$ is only a topological space and showed that whenever $X$ is metric and complete, his definition coincides, on $\ks(X)$ with the topology defined by the Hausdorff metric \cite{CV77}. In 1950, Chabauty \cite{Ch50} revisited the problem in the framework of locally compact Abelian groups, and in 1962, Fell \cite{Fe62} did the same in the context of \CsS. Both defined a modification of the Vietoris topology leading to $\cs(X)$ being compact, but not necessarily Hausdorff. However, both topology coincide if $X$ is compact.

\vspace{.1cm}

\noindent A basis for the Vietoris topology on a topological space $X$ is defined by the family $\Uu(F, \Ff)$ defined below, where $F\subseteq X$ is closed and $\Ff$ is a finite family of open sets:
$$\Uu(F,\Ff)=
   \{A\in \cs(X)\,;\, A\cap F=\emptyset\;\,\&\;\, 
    \forall O\in \Ff\,,\, A\cap O\neq \emptyset 
   \}\,.
$$

\noindent Such a topology is also known under the name of {\em hit-and-miss} \cite{LP94}, since elements of $\Uu(F,\Ff)$ must hit each open set in $\Ff$ and miss $F$. It is tedious but elementary to check that this family defines indeed a basis for a topology. The Fell topology is defined in a similar way by restricting $F$ to be compact. If $X$ is compact the Fell topology and the Vietoris topology coincide and will be called the {\em Hausdorff topology}. The following result is a consequence of the results obtained by Vietoris, Michael and Fell \cite{Vi22,Mi51,Fe62}

\begin{theo}
\label{cspa.th-viet}
Let $X$ be a compact Hausdorff second countable space. Then, equipped with the Hausdorff topology, $\cs(X)$ is Hausdorff, compact, and second countable. 
\end{theo}

\noindent If now $\Gamma$ is a handy groupoid, the following holds

\begin{proposi}
\label{cspa.prop-isGamma}
Let $\Gamma$ be a handy groupoid. Then, if $\cs(\GaO)$ is equipped with the Hausdorff topology, the subset $\isG$ is closed. Consequently it is compact, Hausdorff and second countable.
\end{proposi}

\noindent {\bf Proof: } It is sufficient to prove that $\isG$ is closed, since the rest follows from Theorem~\ref{cspa.th-viet}. Let $M$ be a closed subset of $\GaO$ belonging to the Vietoris closure of $\isG$. If, by contradiction, $M$ is not an element of $\isG$, then there is a unit $x\in M$, and an arrow $\gamma\in\Gamma$, such that $r(\gamma)=x$, $s(\gamma)=y$, and $y\notin M$. For any open neighborhood $O_x$ of $x$ the preimage $V=r^{-1}(O_x)$ is an open neighborhood of $\gamma$. Similarly, given any open neighborhood $O_y$ of $y$ the preimage $W=s^{-1}(O_y)$ is an open neighborhood of $\gamma$. Let then $U$ be an open neighborhood of $\gamma$ contained in $V\cap W$. Since $\Gamma$ is handy, $x\in r(U)\subseteq O_x$ is open and $y\in s(U)\subseteq O_y$ is also open. The separability of $\Gamma$ allows to choose $O_y$ such that $M\cap \overline{O_y}=\emptyset$. Let $F,\Ff$ be defined by $F=\overline{s(U)}$ and $\Ff=\{r(U)\}$. It follows that $F$ is closed and does not intersect $M$. Then $\Uu(F,\Ff)$ is an open neighborhood of $M$ in the Hausdorff topology. Therefore it intersects $\isG$ by assumption. In particular there is a closed invariant subset $N$ of $\GaO$ such that $N\cap F=\emptyset$ and $N\cap r(U)\neq\emptyset$. Therefore, if $x'\in N\cap r(U)$, there is $\gamma'\in U$ such that $x'=r(\gamma')\in N$. As $N$ is invariant, it follows that $s(\gamma')\in s(U)\cap N\subseteq F\cap N=\emptyset$, a contradiction. Hence $M\in\isG$.
\hfill $\Box$

\vspace{.2cm}

\noindent The question addressed in this paper is to control the approximation of an invariant subset by a sequence of periodic orbits. In the general framework we chose to use here, periodic orbits are, in general only periodicoid elements. In the case of \'etale groupoids, these are finite invariant subsets. So the question is {\em ``Can any invariant subset be approximate by a sequence of periodicoid orbits~?''}. The following example shows that the answer is: NO.

\begin{exam}[Periodic Orbits are not Hausdorff-dense]
\label{cspa.exam-noapp}
{\em Let $\as$ be a finite set, called an {\em alphabet}. $\as$ will be seen as equipped with its discrete topology. Let $X=\as^{\ZM}$ be equipped  with the product topology. By Tychonov's Theorem, $X$ is compact. It is also clearly Hausdorff and second countable. The {\em shift} defines a $\ZM$-action by 
$$(Sx)_n=x_{n-1}\,,
   \hspace{2cm}
    x= (x_n)_{n\in\ZM}\in X\,.
$$

\noindent A periodic point, of period $p\in\NM$, is a sequence $x$ such that $x_{n+p}=x_n$ for all $n\in\ZM$. Hence if $w=x_0x_1\cdots x_{p-1}$, then $x=w^\infty$ is the concatenation of the word $w$ infinitely many times on both sides. Given $a, b\in \as$ with $a\neq b$, let $y^{ab}=(y_n)_{n\in\NM}$ be the sequence defined by $y_n=b$ if $n\geq 0$ and $y_n=a$ if $n<0$. Then $y^{ab}=a^\infty\cdot b^\infty$ where the $\cdot$ indicates the position of the separation between negative and nonnegative integers. Then,

(i) The set of periodic point is dense in $X$, a property of the product topology,

(ii) The closure $Y$ of the orbit of $y^{ab}$ cannot be approximated, in the Hausdorff topology by a sequence of periodic orbits.

\noindent To see this it is sufficient to equip $X$ with a compatible metric defining the topology, and to check that the Hausdorff distance between $Y$ to any periodic orbit is bounded from below by a constant independent of the periodic orbit (see \cite{Bec16} Example 5.2.8).
}
\hfill $\Box$
\end{exam}

 \subsection{The Tautological Groupoid: Definition}
 \label{cspa.ssect-taut}

\noindent Let $\Gamma$ be a handy groupoid. Let then $\tsG\subseteq \isG\times \Gamma$ be the set of pairs $(M,\gamma)$ such that $r(\gamma)\in M$. Since $M$ is a closed invariant subset of $\GaO$, it follows that $s(\gamma)\in M$ as well. Equivalently $\gamma\in \Gamma_M$. Then $\tsG$, equipped with the topology induced by the product topology of $\isG\times \Gamma$ is closed, so that it is a locally compact Hausdorff, second countable space. In addition, it becomes a groupoid if endowed with the following structure

\begin{description}

 \item[(T1)] the set of units $\tsGO$ is the set of pairs $(M,x)\in\isG\times \GaO$ such that $x\in M$;

 \item[(T2)] the range and the source maps are defined by $r(M,\gamma)= (M,r(\gamma))$ and $s(M,\gamma)= (M,s(\gamma))$;

 \item[(T3)] the set $\tsGD$ of composable pairs of arrows is the set of pairs $\big((M,\gamma),(M,\eta)\big)$ such that $s(\gamma)=r(\eta)$;

 \item[(T4)] The composition is then given by $(M,\gamma)\circ(M,\eta)=(M,\gamma\circ\eta)$; consequently, the inverse map is given by $(M,\gamma)^{-1}=(M,\gamma^{-1})$. 

\end{description}

\begin{proposi}
\label{cspa.prop-tautgr}
Let $\Gamma$ be a handy groupoid. Then, if $\isG$ is equipped with the Hausdorff topology, the set $\tsG$, endowed with the product topology of $\isG\times \Gamma$ and with the groupoid structure defined by (T1)-(T4), is a handy groupoid as well, which will called the {\em tautological groupoid of $\Gamma$}.  
\end{proposi}

\noindent The term {\em tautological} reflects the fact that it is a ``universal'' object for all continuous fields of groupoids modeled on $\Gamma$, as will be seen in  Section~\ref{cspa.sect-cf} (see in particular   Theorem~\ref{cspa.th-tautfield} and Theorem~\ref{cspa.th-tautfield_2}).

\vspace{.1cm}

\noindent {\bf Proof: } (i) It is straightforward to check that (T1)-(T4) defines a groupoid structure. In addition, the range map is continuous: any open set in $\tsGO$ is containing an open set of the form $\Uu\times V\cap \tsGO$, the preimage of which by $r$ is $\Uu\times r^{-1}(V)\cap \tsG$, which is an open set in $\tsG$. Similarly, the source map is continuous. Since $\Gamma$ is handy, $\GaO$ is Hausdorff, so is $\isG$. Hence, since $\Gamma$ is also Hausdorff, it follows that $\tsG$ is Hausdorff as well. Hence $\tsGD$ is closed in $\tsG\times \tsG$. By the same argument, the groupoid product is continuous and the inverse map also. Since $\tsGO\subseteq \isG\times \Gamma^{(0)}$ is closed, it is compact. At last the range map is open because if $\Uu$ is an open set in $\isG$ and if $U\subseteq \Gamma$ is open, then $\Uu\times U\cap \tsG$ is open in $\tsG$ and such sets generate the topology. Therefore $r(\Uu\times U\cap \tsG)= \Uu\times r(U)\cap \tsGO$ as can be checked easily, which is open. Similarly, the source map $s$ is open.
\hfill $\Box$

\section{Continuous Fields}
\label{cspa.sect-cf}

 \subsection{Continuous Fields of Groupoids}
 \label{cspa.ssect-cfgr}

\noindent The definition below was proposed and studied in \cite{LR01} as a way to construct continuous fields of \CsS, in connection with the concept of deformation quantization. The latter concept has a long history (see \cite{LR01}), but Rieffel deserves the credit for providing a rigorous satisfactory definition for deformations of \Css \cite{Rie89a,Rie89b,Rie90a,Rie90b,Rie93}. 

\begin{defini}[\cite{LR01}]
\label{cspa.def-contfieldgr}
A  field of groupoids is a triple $(\Gamma,T, p)$ where $\Gamma$ is a groupoid, $T$ a set and $p:\Gamma\to T$ is a 
surjective 
map such that, if $p_0$ is the restriction of $p$ to the unit space $\GaO$, then $p=p_0\circ r=p_0\circ s$.

\noindent If, in addition, $\Gamma$ is a topological groupoid, $T$ a topological space and if $p$ is continuous and open, then $(\Gamma,T, p)$ is called a continuous field of topological groupoids.

\noindent If, moreover, $\Gamma$ is locally compact ({\em resp.} handy) and if $T$ is Hausdorff, and if $p$ is continuous and open, then $(\Gamma,T, p)$ is called a continuous field of locally compact ({\em resp.} handy) groupoids.
\end{defini}

\begin{rem}
\label{cspa.rem-onto}
{\em If $p$ is not onto, then $T$ can be replaced by $\hT= p(\Gamma)$. 
Since $p$ is assumed to be open, $\hT$ is open in $T$ and $p$ defines a continuous surjective open map
between $\Gamma$ and $\hT$.
}
\hfill $\Box$
\end{rem}

\noindent The following result shows that the structure of continuous field of groupoids is coded in the restriction $p_0:\Gamma^{(0)}\to T$ of the map $p$.

\begin{lemma}
\label{cspa.lem-p}
Let $(\Gamma,T, p)$ be a continuous field of groupoids, Then $p_0:\Gamma^{(0)}\to T$ is a continuous, open, surjective map. Moreover if 
$(\Gamma,T, p)$ is a continuous field of handy groupoids 
then $T$ is necessarily compact and second countable space.
\end{lemma}

\noindent {\bf Proof: }
Since the restriction of a continuous map is continuous, it suffices to show that $p_0$ is surjective and open. By definition of a continuous field of groupoids, $p(\gamma)=p(r(\gamma))$ follows implying $p(\Gamma)=p(r(\Gamma))$. Thus, $p_0$ is surjective as $r(\Gamma)=\Gamma^{(0)}$. 
Let $U\subseteq \Gamma^{(0)}$ be open. Since the range map is continuous, $r^{-1}(U)\subseteq\Gamma$ is open. Thus, $p_0(U)=p(r^{-1}(U))$ implies that $p_0(U)$ is an open subset of $T$ since $p$ is an open map. Hence, $p_0$ is an open map. If the field of groupoids is handy, $\Gamma^{(0)}$ is compact and the topology of $\Gamma$ is  is second countable by definition. Since $(\Gamma,T, p)$ is a continuous field of handy groupoids, $T$ is Hausdorff. Hence, $T$  
is necessarily compact as a continuous image of the compact set $\Gamma^{(0)}$ under the map $p_0$. 
Moreover,
 $T$ is second countable as the image of $\Gamma$ by the continuous, open, surjective map $p$, see e.g. \cite[Proposition~A.2.1]{Bec16}.
\hfill $\Box$

\medskip

\noindent The next result result contains the main construction of this work. It describes how to associate a continuous field of groupoids to a given initial groupoid $\Gamma$.

\begin{theo}[The tautological field of a groupoid]
\label{cspa.th-tautfield}
Let $\Gamma$ be a handy groupoid.
Let $p_\Gamma:\tsG\to \isG$ be the map defined by $p_\Gamma(M,\gamma)=M$ for $M\in \isG$ and $\gamma\in \Gamma_M$. Then $(\tsG,\isG,p_\Gamma)$ is a continuous field of handy groupoids, called the tautological field of $\Gamma$.
\end{theo}

\noindent {\bf Proof: }
(i) Since, by Proposition~\ref{cspa.prop-isGamma}, $\isG$ is compact, Hausdorff and second countable, it suffices to prove that $p_\Gamma$ is continuous, open and surjective. By definition $p_\Gamma(\tsG)=\isG$ holds implying that $p_\Gamma$ is surjective. Since the topology on $\tsG$ is induced by the product topology of $\isG\times \Gamma$ and since $p$ is the first projection, $p_\Gamma$ is clearly a continuous and open map.
\hfill $\Box$

\medskip

\noindent 
The construction of the tautological field of a groupoid is compatible with the structure given by the presence of a continuous Haar system
on the initial
 groupoid. This is shown in the following result.

\begin{lemma}
\label{cspa.lem-haarrest}
Let $\Gamma$ be a handy groupoid equipped with a left-continuous Haar system $\mu$. Then the restriction $\mu_M=x\in M\mapsto \mu^x$ is a well-defined left-continuous Haar system on $\Gamma_M$. In addition, if $(M,x)\in \tsGO$ and if $\delta_M$ denotes the Dirac measure supported by $M$ on $\isG$, then $\mu^{(M,x)}=\delta_M\times \mu^x$ on $\{M\}\times \Gamma^x$, defines a left-continuous Haar system on $\tsG$.
\end{lemma}

\noindent  {\bf Proof: }For each $x\in \GaO$, $\mu^x$ is a positive Borel measure on the $r$-fiber $\Gamma^x$ satisfying the conditions described in Definition~\ref{cspa.def-HaarSyst}. If $M\in \isG$ it follows that for any $x\in M$ and any $\gamma\in\Gamma^x$, then $s(\gamma)\in M$ as well. In particular $\Gamma^x\subseteq \Gamma_M$, namely $\Gamma_M^x=\Gamma^x$. Therefore the restriction $\mu_M=x\in M\mapsto \mu^x$ can be seen as a left-continuous Haar system on $\Gamma_M$. The rest is left to the reader.
\hfill $\Box$

\medskip

\noindent 
By construction the tautological field of $\Gamma$ is a ``fibered'' groupoid  $\tsG$
over the ``base space'' $\isG$ through the ``projection map'' $p_\Gamma$. In many situations it can be useful to have the freedom of changing the base space. For that one can generalize the pull-back construction typical of the theory of fiber bundles. 
Here the argument will be only sketched leaving the details for the interested reader.
Let $Y$ be a ``nice'' topological space (Hausdorff, second countable, (locally-) compact, etc.) and $f:Y\to \isG$ be a continuous map. Consider the space
$$
f^*\tsG:=
	\big\{(y,(M,\gamma))\in Y\times \tsG\; :\; f(y)=p_\Gamma(M,\gamma) \big\}
	\;\subseteq\; Y\times \tsG
$$
\noindent and the projection $p_f:f^*\tsG\to Y$ given by
$p_f(y,(M,\gamma))=y$. Then the triple $(f^*\tsG, Y, p_f)$ turns out to be  a continuous field of  groupoids over $Y$ called the \emph{pull-back field induced by $f$}. Moreover, the following diagram turns out to be commutative
$$
\begin{diagram}
f^*\tsG        &   \rTo^{\hat{f}}                &   \tsG \\
       \dTo^{p_f}     &   &   \dTo_{p_\Gamma} \\
Y&       \rTo^f           & \isG  \;
\end{diagram}
$$
\noindent where the map $\hat{f}$ is the projection on the second factor defined by $\hat{f}(y,(M,\gamma))=(M,\gamma)$. The commutativity of the 
diagram above implies that the two fibers  $f^*\tsG|_y:= p_f^{-1}(y)$ and $\tsG|_{f(y)}:=p_\Gamma^{-1}(f(y))$ are related by the projection $\hat{f}$ and this provides the (natural) groupoid isomorphism $f^*\tsG|_y\simeq\tsG|_{f(y)}$.

\medskip

\noindent 
A final question that it is worth to answer is the following: Given a continuous field of  groupoids $(\Gamma,T,p)$ what is the relation with the tautological field $(\tsG,\isG,p_\Gamma)$ associated to $\Gamma$?
In the case of a field of handy groupoids the answer is described in Theorem~\ref{cspa.th-tautfield_2} and is the following: if $(\Gamma,T,p)$ is a continuous field of handy groupoids, then it can be embedded as a
sub-groupoid of the tautological field $(\tsG,\isG,p_\Gamma)$ associated to $\Gamma$. This result is relevant since it allows to extend results valid for the tautological field to generic continuous field of handy groupoids.
Before proving Theorem~\ref{cspa.th-tautfield_2} a preliminary result is needed.

\begin{lemma}
\label{cspa.lem-isT}
Let $(\Gamma,T, p)$ be a continuous field of handy groupoids, with $T$ Hausdorff and denote by $p_0:\Gamma^{(0)}\to T$ be the restriction to the unit space of the map $p$. Then
$$\is_T
	:=\{M\subseteq \Gamma^{(0)}\,:\, M=p_0^{-1}(\{t\}) \text{ for some } t\in T\}
$$
is a closed subset of $\isG$ in the Hausdorff topology.
\end{lemma}

\noindent {\bf Proof: }{
According to Lemma~\ref{cspa.lem-p}, $p_0$ is an open, continuous and surjective map. Since $p_0$ is continuous and $p=p_0\circ r=p_0\circ s$, the preimage $p_0^{-1}(\{t\})\subseteq\Gamma^{(0)}$ is closed and invariant, namely $\is_T\subseteq\isG$. Let $M\in\overline{\is_T}$ and $M_n\in\is_T$ for $n\in\NM$ such that $M_n\to M$ in the Hausdorff topology. Thus, for each $x\in M$ there is an $x_n\in M_n$ such that $x_n\to x$ in $\Gamma^{(0)}$. Set $t:=p_0(x)\in T$ and $t_n:=p_0(x_n)$ for $n\in\NM$. Due to continuity of $p_0$, $(t_n)_n$ converges to $t$ in $T$. Since $t_n=p_0(y)$ holds for all $y\in M_n$ and $x\in M$ was arbitrary, the inclusion $M\subseteq p_0^{-1}(\{t\})$ follows.

\vspace{.1cm}

\noindent 
In order to show the converse inclusion $M\supseteq p_0^{-1}(\{t\})$, assume by contradiction that it does not hold. Thus, there is an $y\in \Gamma^{(0)}\setminus M$ such that $p(y)=t$. Let $U$ be an open neighborhood of $y$ such that $\overline{U}\cap M=\emptyset$ and $\overline{U}$ is compact, which is possible since $\Gamma$ is a handy groupoid. Hence, $y\in U$ and $p(y)=t$ imply that $p_0(U)\subseteq T$ is an open neighborhood of $t$ since $p_0$ is an open map. By convergence $t_n\to t$, there is an $n_U\in\NM$ satisfying $t_n\in p_0(U)$ for $n\geq n_U$. Consequently, $p_0^{-1}(t_n)\cap U\neq \emptyset$ follows. Since $M_n=p_0^{-1}(t_n)$, this yields $M_n\cap U\neq \emptyset$ for $n\geq n_U$, a contradiction as $M_n\to M$ in the Hausdorff topology and $M\cap \overline{U}=\emptyset$. Consequently, $M=p_0^{-1}(\{t\})$ follows implying $M\in\is_T$, namely $\is_T\subseteq\isG$ is a closed subset.
\hfill $\Box$

\medskip

\begin{theo}
\label{cspa.th-tautfield_2}
Let $(\Gamma,T, p)$ be a continuous field of handy groupoids, with $T$ Hausdorff and $\is_T$ be the closed subset of $\isG$ defined in Lemma~\ref{cspa.lem-isT}. Define $p_T:\is_T\to T\,,\, M\mapsto p_0(x)\,,$ where $x$ is any element in $M$. Then $p_T$ is a homeomorphism. Furthermore, there is an open, continuous, surjective map $\iota_T:\Gamma\to p_\Gamma^{-1}(\is_T)$, where $p_\Gamma^{-1}(\is_T)$ is a closed subgroupoid of $\tsG$, satisfying $p=p_T\circ p_\Gamma \circ \iota_T$.
\end{theo}

\noindent {\bf Proof: }{
The map $p_T:\is_T\to T\,,\, M\mapsto p_0(x)\,$ where $x\in M$ is well-defined since $p_0(x)=p_0(y)$ for all $x,y\in M\in\is_T$. Due to Lemma~\ref{cspa.lem-p}, the restriction $p_0:\Gamma^{(0)}\to T$ of $p$ is an open, continuous and surjective map. Thus, $p_T$ is surjective. If $t:=p(M)=p(N)$ for $M,N\in\is_T$, then $M=p_0^{-1}(\{t\})$ and $N=p_0^{-1}(\{t\})$ follow by definition of $\is_T$. Thus, $M=N$ follows implying that $p_T$ is injective. In order to show that $p_T$ is continuous, it suffices to restrict to sequences as all involved spaces are second countable, c.f. Proposition~\ref{cspa.prop-isGamma} and Lemma~\ref{cspa.lem-p}. Let $M_N\to M$ in $\is_T$. Let $x_n\in M_n$ and $x\in M$ such that $x_n\to x$ in $\Gamma^{(0)}$. Thus the continuity of $p_0$ implies
$$
p_T(M_n)=p_0(x_n)\to p_0(x)=p_T(M)\,,
$$
\noindent proving the continuity of $p_T$. According to Lemma~\ref{cspa.lem-isT} and Proposition~\ref{cspa.prop-isGamma}, $\is_T$ is a compact, Hausdorff space in the Hausdorff topology. Since $T$ is also compact and Hausdorff by Lemma~\ref{cspa.lem-p}, the previous considerations yield that $p_T$ is a homeomorphism.

\vspace{.1cm}

\noindent
If $M,N\in\is_T$ and $M\cap N\neq\emptyset$, then there is an $x\in M\cap N$ and $p_T(M)=p_0(x)=p_T(N)$. Thus, $M=N$ follows as $p_T$ is injective. Consequently, $\bigsqcup_{M\in\is_T} M=\Gamma^{(0)}$ is a disjoint union. With this at hand, define $i_T:\Gamma^{(0)}\to\is_T\,,\, x\mapsto M(x)\,,$ where $M(x)$ is the unique element in $\is_T$ containing $x$. By the previous considerations $i_T$ is well-defined and surjective. Furthermore, $i_T$ is an open map: Let $U\subseteq \Gamma^{(0)}$ open, then $i_T(U)=\{M\in\is_T\,:\, M\cap U\neq\emptyset\}$ which by definition is an open set in the Hausdorff topology on $\is_T$. In addition, $i_T$ is continuous: If $x_n\to x$ in $\Gamma^{(0)}$, then 
$$
p_T(i_T(x_n))=p_0(x_n)\to p_0(x)=p_T(i_T(x))
$$ 
\noindent by continuity of $p_0$. Since $p_T$ is a homeomorphism, this yields $i_T(x_n)\to i_T(x)$ in the Hausdorff topology on $\is_T$, namely, $i_T$ is continuous.

\vspace{.1cm}

\noindent 
Consider the closed subgroupoid $p_\Gamma^{-1}(\is_T)$ of the tautological groupoid $\tsG$. Define $\iota_T:\Gamma\to p_\Gamma^{-1}(\is_T)$ by $\gamma\mapsto (i_T(r(\gamma)),\gamma)$. Since all involved maps are continuous, open and surjective, the image $\iota_T(\Gamma)$ equals to $p_\Gamma^{-1}(\is_T)$ and $\iota_T:\Gamma\to p_\Gamma^{-1}(\is_T)$ is an open, continuous map. By construction, we have
$$p_T\big(p_\Gamma\big(\iota_T(\gamma)\big)\big) 
	= p_T\big(p_\Gamma\big( (i_T(r(\gamma)),\gamma)\big)\big) 
	= p_T\big(i_T(r(\gamma))\big)
	= p_0(r(\gamma))=p(\gamma)
$$
\noindent since $r(\gamma)\in i_T(r(\gamma))$.
\hfill $\Box$

 \subsection{Continuous Field of $C^\ast$-Algebras: a Reminder}
 \label{cspa.ssect-cfcst}

\noindent The concept of continuous field of \Css has been investigated at least since Kaplansky in \cite{Ka51}. In particular he gave the first proof that a self-adjoint continuous section of such a field has a spectrum varying continuously with the parameter. This concept was developed and used by Fell \cite{Fe60,Fe61}, Tomiyama \cite{TT61,To62} and later studied by Dixmier and Douady \cite{DD63,Di69}, who gave a complete definition. It is worth noticing that continuous fields are sometimes called {\em ``bundle''} \cite{KW95}. But this name does not corresponds to the concept of bundle in Differential Geometry, as continuous fields might not be locally trivial.

\vspace{.1cm}

\noindent Let $T$ be a topological space. For each $t\in T$ let $\Aa_t$ be a \CS. The family $\Aa=(\Aa_t)_{t\in T}$ is called a field. A {\em section} or a {\em field of operators} is a family $A=(A_t)_{t\in T}$ for which $A_t\in\Aa_t$ for all $t\in T$. Given a set $\Upsilon$ of sections,  $\Upsilon_t$ will denote the set of $A_t\in \Aa_t$ such that $A\in \Upsilon$.

\begin{defini}
\label{cspa.def-cfC*}
A field $\Aa$ of \Css will be called {\em continuous} whenever there is a set $\Upsilon$ of sections with the following property

\begin{description}
  \item[(CF1)] for each $t\in T$, then $\Upsilon_t$ is a dense $\ast$-subalgebra of $\Aa_t$,

  \item[(CF2)] for each $A\in\Upsilon$, the map $t\in T\mapsto \|A_t\|$ is continuous,

  \item[(CF3)] a section $B$ belongs to $\Upsilon$ if and only if for each $t\in T$, for each $\epsilon >0$, there is an open neighborhood $U$ of $t$ in $T$ and a section $A\in \Upsilon$, such that $\|A_s-B_s\|<\epsilon $ for $s\in U$.
\end{description}

\noindent Then an element of $\Upsilon$ is called a continuous section. 
\end{defini}

\noindent The set $\Upsilon$, endowed with the pointwise addition, scalar multiplication, product and adjoint, becomes a $\ast$-algebra. However, in general, this is not a normed algebra. A continuous section $A\in \Upsilon$ will be called bounded whenever $\|A\|=\sup_{t\in T}\|A_t\|$ is finite. In this case, $\|\cdot\|$ is a $C^\ast$-norm and the set of bounded continuous sections becomes a \CS, which will be denoted by $\Cc_b(T,\Aa)$. Then the pointwise multiplication by a complex valued bounded continuous function on $T$ makes it a $\Cc_b(T)$-$C^\ast$-module \cite{Bl96}. If, in addition, $T$ is locally compact, then the set of continuous sections vanishing at infinity is also a \Cs when endowed with that norm, which will be denoted by $\Cc_0(T,\Aa)$. Again, it can be seen as a $\Cc_0(T)$-$C^\ast$-module. 

\vspace{.1cm}

\noindent It is not difficult in practice to check that a field of \Cs satisfies (CF1). The axiom (CF3) looks hard to check, but, as can be seen in \cite{Di69}, $\Upsilon$ can be replaced by a smaller set generating it, and in most situations, it is not difficult to built such a generating set. The tricky part is the axiom (CF2). Proving the continuity of the norm is always a problem. For this reason, Rieffel defined the concept of upper and lower semi-continuous field, by replacing the continuity property (CF2) by the corresponding semi-continuity. 

\vspace{.1cm}

\noindent Following Rieffel \cite{Rie89a}, if $T$ is locally compact and if $\Bb$ is a \CS, it will be called a  $\Cc_0(T)$-$C^\ast$-module whenever there is an injective $\ast$-homomorphism defined on $\Cc_0(T)$ with values in the center of the multiplier algebra $M(\Bb)$ (see Section~\ref{cspa.ssect-algB} for a definition of the multiplier algebra). Then $\Bb$ can be seen as the set of continuous sections, vanishing at infinity, of an upper semi-continuous field of \Css \cite{Rie89a}. The lower semi-continuity requires a specific property, for instance the existence of a suitable representation of $\Bb$. In the case of groupoid $C^\ast$-algebras, the left regular representation does the job. But then, the difference between the reduced and the full algebra may prevent those fields to be more than only semi-continuous, one way or the other.

\vspace{.1cm}

\noindent It is also worth noticing that the tensor product of two continuous fields may not be continuous. This is due to the non uniqueness of the tensor product of \CsS, more precisely, on the algebraic tensor product of two \Css there are, in general, several $C^\ast$-norms, one being minimal and one maximal. The point of view of $\Cc(T)$-modules turns out to be very useful to investigate the tensor product \cite{Bl96,KW95}.

\vspace{.1cm}

\noindent Building upon such a concept, Landsman and Ramazan \cite{LR01} showed that continuous field of groupoids lead to semi-continuous fields of \Css depending upon whether the full or the reduced \Cs are used. For convenience, this theorem will be reformulated for handy groupoids, since the purpose of this paper is to consider applications rather than a general mathematical result. In \cite{LR01}, only the case of a trivial cocycle is considered

\begin{theo}[see \cite{LR01}, Section 5]
\label{cspa.th-LR01}
Let $(\Gamma, T,p)$ be a continuous field of handy groupoids. Let $\mu$ be a Haar system on $\Gamma$. Then,  the field $\left(\CG^\ast(\Gamma_t,\mu)\right)_{t\in T}$ is upper semi-continuous and the field  $\left(\CG_{red}^\ast(\Gamma_t,\mu)\right)_{t\in T}$ is lower semi-continuous.
\end{theo}

\noindent As a result, if the field has the {\em ``Amenability Property''} (see \cite{Rie89a}, Hypothesis 2.6), namely if each fiber satisfies $\CG_{red}^\ast(\Gamma_t,\mu)=\CG^\ast(\Gamma_t,\mu)$, then the field obtained this way is continuous. Here some remark is in order: since the fiber $\Gamma_t$ is a sub-groupoid of $\Gamma$, it follows that $\mu$ restricts to $\Gamma_t$ is a left-continuous Haar system as well.

 \subsection{Continuous Fields of $2$-cocycles}
 \label{cspa.ssect-cfcoc}

\noindent In the previous result, the \Css are defined without reference to a cocycle, meaning that $\sigma=1$. Adding a cocycle on a handy groupoid $\Gamma$, defines a field of cocycles $\sigma_M$ on each fiber of the tautological groupoid.
As shown in Theorem~\ref{cspa.th-tautfield_2}, this is not a restriction: any continuous field of cocycles can be seen as a continuous field over the tautological groupoid. In Physics, such fields occur whenever a uniformly continuous magnetic field is added, which might be uniform or not, as shown in Section~\ref{cspa.ssect-numagn}.
In the study of the Heisenberg group \cite{Rie89b}, which corresponds to the Weyl quantization procedure, a similar problem occurs where the magnetic field is replaced by a {\em semiclassical} parameter, generalizing the concept of {\em Planck constant}, leading to a continuous field of cocycles. It becomes important, in view of potential applications, to include the possibility of changing the cocycle along the way and to get results about the continuity of the field of \CsS. 

\vspace{.1cm}

\noindent In the case for which $\Gamma$ is actually a discrete group, the space of all possible $2$-cocycles can be equipped with a topology making it a compact space (\cite{Rie89a}, Corollary 2.8). But in the case of a handy groupoid, even an \'etale one, this is not the case. To see the difficulty, let the case of $1$-cocycles, called modules,  be considered as in Example~\ref{cspa.exam-module}: a module over the groupoid $X\rtimes \ZM$ is entirely defined by a continuous function $h:X\to \SM^1$. To put a topology on the set of such modules, it is sufficient to define a topology on the space $\Omega= \Cc(X,\SM^1)$ of such functions. The problem is that compactness can be achieved using the pointwise convergence, but then $\Omega$ is not closed, since discontinuous function will show up in the closure. If then the uniform convergence is defined, the space $\Omega$ is closed, but not compact, not even locally compact. In the definition of Dixmier-Douady \cite{Di69}, the space of parameters needs not be locally compact in general, so that it is not an obstacle. However, the argument provided by Rieffel \cite{Rie89a} to prove the upper semi-continuity fails if $\Omega$ is not locally compact, so that the definition of Dixmier-Douady cannot be used since the upper semi-continuity of the norm cannot be proved. As a compromise, Rieffel propose to use a {\em continuous field of $2$-cocycles} instead. Adjusting his definition to the groupoid case is provided below. The notations are the following: let $\tsGD\subseteq \GaD\times T$ denote the set of $(\gamma,\eta, t)$, such that $(\gamma,\eta)\in \Gamma_t^{(2)}$. Clearly, since the map $p$ is continuous, it follows that  $\tsGD$ is closed in $\GaD\times T$. Then it will be endowed with the induced topology. 

\begin{defini}
\label{cspa.def-cfcoc}
Let $(\Gamma, T,p)$ be a continuous field of handy groupoids. Then a continuous field $\sigma=(\sigma_t)_{t\in T}$ of $2$-cocycles is a continuous function $\sigma: \tsGD \to \SM^1$ such that, for each $t\in T$, the map $\sigma_t: (\gamma, \eta)\in \Gamma_t^{(2)}\mapsto \sigma(\gamma, \eta,t)\in \SM^1$ is a normalized $2$-cocycle.
\end{defini}

\begin{proposi}
\label{cspa.prop-cfcoc}
Let $\Gamma$ be a handy groupoid with normalized $2$-cocycle $\sigma$ and consider the continuous field of groupoids $(\tsG,\isG),p_\Gamma)$ defined by the tautological groupoid. Then $\hat{\sigma}:\tsGD\to\SM^1$ defined by $\hat{\sigma}\big((M,\gamma),(M,\eta)\big):=\sigma(\gamma,\eta)$ is a continuous normalized $2$-cocycle.
\end{proposi}
\noindent  {\bf Proof: } This is straightforward.
\hfill $\Box$

\section{Proof of Theorem~\ref{cspa.theo-contSp}}
\label{cspa.sect-th3}

\noindent In this Section all the previous concepts and results will be used to prove Theorem~\ref{cspa.theo-contSp}. The assumptions made in the statement of this theorem imply that the groupoid $\Gamma$ 
admits a Haar system,
 thanks to Theorem~\ref{cspa.th-haaropen}. The strategy for the proof is not new. It follows the ideas first used by Renault \cite{Re80}, by Elliott \cite{El82}, by Rieffel \cite{Rie89a}, by Landsman and Ramazan \cite{LR01} and by Latr\'{e}moli\`{e}re \cite{La04}. Since the aim of this work is to open this approach to a general audience that is not familiar with groupoid \CsS, the proof is provided.

\vspace{.1cm}

\noindent The proof is organized in two parts. The first, entitled {\em ``Upper Semi-Continuity''} follows the strategy proposed by Rieffel: if $A=(A_t)_{t\in T}$ is a continuous field of \CsS, then $J_t$ denotes the closed two sided ideal generated by continuous sections vanishing at $t$. If it can be proved that $A_t$ is $\ast$-isomorphic to the quotient algebra $A/J_t$, then the norm on $A_t$ is a quotient norm, namely the infimum over continuous semi-norms, leading to the upper semi-continuity (see \cite{Rie89a}, Proposition 1.2). 
The second part, entitled {\em ``Lower Semi-Continuity''}, uses the left regular representation to show that the norm on the reduced \Cs is the supremum of continuous functions, leading to the lower semi-continuity. It will be shown that adding a continuous field of $2$-cocycles does not change the strategy (the critical technicality is the proof of Lemma~\ref{cspa.lem-multh} below). 
This generalizes Theorem~\ref{cspa.th-LR01} to continuous fields of handy groupoids with continuous field of $2$-cocycle.

\vspace{.1cm}

\noindent Section~\ref{cspa.ssect-algB} contains preparatory technical results intended to prove that the norm on $A_t$ is the quotient norm. First the dense sub-algebra provided by continuous functions with compact support on the various groupoids are considered. In Section~\ref{cspa.ssect-full}, it is shown, by using general results on \CsS, that the conclusions extend to the full \CS, leading to the upper semi-continuity of the norm. In Section~\ref{cspa.ssect-reduced} the left regular representation is used to show that the norm on the reduced \Cs is the supremum of continuous functions, leading to the lower semi-continuity. Adding a continuous field of $2$-cocycles changes the formulas defining the left regular representation. But the result is the same as without it.
Finally, Section~\ref{cspa.ssect-Prtheo} provides the proof of Theorem~\ref{cspa.theo-contSp}.

\vspace{.1cm}

\noindent
Instead of studying a general continuous field of handy groupoids, the considerations are restricted to the tautological field $(\tsG,\isG,p_\Gamma)$. Thanks to Theorem~\ref{cspa.th-tautfield} (along with Lemma \ref{cspa.lem-haarrest}) and Theorem~\ref{cspa.th-tautfield_2}, this is not a restriction in general (see also Remark~\ref{cspa.rem-tautfield2} and Remark~\ref{cspa.rem-tautfield3} below). In light of this, throughout the section the tautological field $(\tsG,\isG,p_\Gamma)$ is considered equipped with a continuous field of $2$-cocycle $\sigma:\tsGD\to\SM^1$.

 \subsection{The Algebra $\Aa$ and its Ideals}
 \label{cspa.ssect-algB}

\noindent In this Section $\Aa$ will denote the $\ast$-algebra $\Aa=\Cc_c(\tsG,\mu,\sigma)$. It is reminded that its topology is defined as the inductive limit topology, namely a net $(f_\alpha)_{\alpha\in A}$ in $\Aa$ converges to $f\in \Aa$ if there is a compact subset $K\subseteq \tsG$ supporting each $f_\alpha$ and $f$ and $\sup_{(M,\gamma)\in K}|f(M,\gamma)-f_\alpha(M,\gamma)|\to 0$. With this topology the product and the adjoint are continuous.

\vspace{.1cm}

\noindent Let $\Cc$ denote the space $\Cc(\isG)$ of complex valued continuous functions on the space of closed invariant subsets of the groupoid $\Gamma$. Since $\isG$ is compact, it follows that $\Cc$ is a unital Abelian \Cs when endowed with the uniform norm

$$\|h\|=\sup_{M\in\isG}|h(M)|\,.
$$

\noindent For $F\in\isG$, let $\Ii_F$ denotes the set of functions $h\in\Cc$ such that $h(F)=0$. This is a maximal (closed, two-sided) ideal of $\Cc$ and any maximal ideal is of this type.

\vspace{.1cm}

\noindent As a reminder, a multiplier of $\Aa$ is a pair $(L,R)$ of linear continuous operators on $\Aa$ such that, for any $f, g\in\Aa$
\begin{equation}
\label{cspa.eq-mult}
L(fg)=L(f)g\,,
   \hspace{1cm}
    R(fg)=f R(g)\,,
     \hspace{1cm}
      R(f)g=f L(g)\hspace{1cm}\mbox{\rm  (compatibility)}\,.
\end{equation}
%
\noindent If $M=(L,R)$, it is common to write $M f=L(f)$ and $f M=R(f)$. The relation above are expressing the associativity of such a product. If $g\in \Aa$ then $g$ can be seen as a multiplier through 
$$g=(L_g,R_g)\,,
   \hspace{2cm}
    L_g(f)=gf\,,\;\;R_g(f) =fg\,.
$$

\noindent If $(L,R)$ and $(L',R')$ are two multipliers then $(L+L',R+R')$ is a multiplier, called their sum. If $\lambda \in \CM$ is a scalar then $\lambda(L,R)=(\lambda L, \lambda R)$ is also a multiplier. The product is defined by $(L,R)(L',R')=(LL',R'R)$, and it can be checked that this is a multiplier. At last, if $S:\Aa\to \Aa$ is linear and continuous, the operator $S^\ast$ is defined by $S^\ast(f)=\left(S(f^\ast)\right)^\ast$. Then the adjoint of a multiplier is defined by $(L,R)^\ast= (R^\ast,L^\ast)$. It can be checked that the adjoint is also a multiplier. Consequently, endowed with these algebraic operations, the set of multipliers is a $\ast$ algebra. An elementary exercise shows that given any two multipliers $(L,R)$ and $(L',R')$, then $LR'=R'L$ and, similarly $L'R=RL'$. In particular if $L=R$ then $(L,R)(L',R')= (LL',R'R)= (RL',R'L)= (L'R,LR')= (L'L,RR')=(L',R')(L,R)$.

\begin{lemma}
\label{cspa.lem-multh}
For $h\in\Cc$ and $a\in\Aa$ let $m_h:\Aa\to \Aa$ be the linear operator defined by 
$$(m_h f)(M,\gamma)=h(M)\,f(M,\gamma)\,,
    \hspace{2cm}
     M\in \isG\,,\, \gamma\in\Gamma_M\,.
$$

\noindent Then, $m_h$ defines a multiplier of $\Aa$ commuting with all multipliers of $\Aa$. In addition $m_h^\ast=m_{\overline{h}}$.
\end{lemma}

\noindent  {\bf Proof: } It is immediate to check that the pair $\hm_h=(m_h,m_h)$ of linear operators on $\Aa$ satisfies the algebraic relations~\eqref{cspa.eq-mult} defining a multiplier $(L,R)$. For indeed, since $h(M)\in\CM$, $m_h(fg)(M,\gamma)$ equals to
$$
h(M) fg(M,\gamma) 
  =
   \int_{\Gamma^{r(M,\gamma)}} 
    \sigma_M(\eta,\eta^{-1}\circ \gamma)\,
     h(M) f(M,\eta) g(M,\eta^{-1}\circ \gamma)
      d\mu^{r(M,\gamma)}(\eta)
  = 
      f(m_h g)(M,\gamma)\,.
$$
\noindent The relation $L=R$ gives the two other relations. In addition, it implies that this multiplier commutes with all elements of $\Aa$. Moreover, $m_h^\ast(f)=\left(m_h (f^\ast)\right)^\ast$ holds by definition, so that
\begin{eqnarray*}
m_h^\ast(f)(M,\gamma) 
  =
    \overline{\sigma_M(\gamma,\gamma^{-1})}\;
     \overline{h(M)}\;
      \overline{f^\ast(M,\gamma^{-1})}
  =
     \overline{h(M)}\,f(M,\gamma)
      =(m_{\overline{h}}f)(M,\gamma)\,.
\end{eqnarray*}
\noindent At last, for $f\in\Aa$ and $h\in\Cc$
$$\supp(m_h f)\subseteq \supp(f)\,,
   \hspace{2cm}
    |m_h f(M,\gamma)|\leq \|h\||f(M,\gamma)|\,,
$$
\noindent as can be checked by inspection. It follows that the map $f\in\Aa\mapsto m_h f\in\Aa$ is continuous.
\hfill $\Box$

\begin{defini}[See for instance \cite{Bl96}]
\label{cspa.def-Cmodule}
The multiplier $m_h$ will be written as $m_h(f)=hf=fh$. Then $\Aa$ becomes a $\ast$-$\Cc$-module.
\end{defini}

\noindent Let $n\in\NM$ be an integer. Then $M_n(\Aa)$ denotes the $n\times n$ matrices with elements in $\Aa$. The algebraic operations, sum, scalar multiplication, product and adjoint are defined in the usual way. An element $T\in M_n(\Aa)$ will be called {\em positive} whenever there is $S\in M_n(\Aa)$ such that $T=S^\ast S$. Then the relation $T\leq T'$ means $T'-T$ is positive. If $h\in \Cc$ then the multiplier $m_h$ extend to $M_n(\Aa)$ by multiplying $h$ to each matrix element.

\begin{proposi}
\label{cspa.prop-posestim}
Given $h\in\Cc$, then, for any $S\in M_n(\Aa)$, 
$$S^\ast h^\ast h S\leq \|h\|^2 S^\ast S\,.
$$

\noindent In particular, given $\{f_1,\cdots,f_n\}$ in $\Aa$, the matrix
$ T^h_{k,l}=f_k^\ast h^\ast h f_l$ satisfies 
$$T^h\leq \|h\|^2 T^1\,,
   \hspace{1cm} \mbox{\rm\em where}\hspace{1cm} 
    1(M)=1\,\;\forall M\in\isG\,.
$$
\end{proposi}

\noindent {\bf Proof: }Let $g$ be the function defined by
$$g(M)=\sqrt{\|h\|^2-|h(M)|^2}\,,
   \hspace{2cm}
    |h(M)|^2= (h^\ast h)(M)\,.
$$

\noindent Then as a multiplier $g^\ast=g\geq 0$. Moreover, $g^\ast g+h^\ast h=\|h\|^2 \,1$ where $1(M)=1$ for all $M\in\isG$. Thus 
$$\|h\|^2 S^\ast S=
   S^\ast h^\ast h S+
    S^\ast g^\ast g S\,,
$$

\noindent proving the inequality. If $f=\{f_1,\cdots,f_n\}$ in $\Aa$, let $S_f$ be the matrix defined by
$$S_f=\left[
\begin{array}{ccc}
f_1 & \cdots & f_n\\
0 & \cdots & 0\\
\vdots & \ddots & \vdots\\
0 & \cdots & 0
\end{array}
\right]\;\;\in M_n(A)\,.
$$

\noindent Then $T^h=S_f^\ast h^\ast h S_f$ and the second inequality follows from the first.
\hfill $\Box$

\begin{coro}
\label{cspa.cor-repAa}
Let $\rho$ be a representation of $\Aa$ in the Hilbert space $\Hh$. Then for each $h\in\Cc$ there is a bounded operator $\rho(h)$, commuting with the operators $\rho(f)$ for $f\in\Aa$ such that $\rho$ becomes a $\Cc$-module $\ast$-homomorphism.

\noindent Consequently, the $\Cc$-action extends by continuity to both the full and the reduced \Css of $(\tsG,\mu,\sigma)$.
\end{coro}

\noindent  {\bf Proof: }If $\rho$ is a representation, by definition the map $f\to \rho(f)$, is a weakly continuous $\ast$-homomorphism such that the linear span of the vectors of the from $\rho(f)\xi$ with $f\in \Aa$ and $\xi\in \Hh$ is dense in $\Hh$. Then, if $\xi_1,\cdots,\xi_n$ are vectors in $\Hh$ and if $f_1,\cdots,f_n$ are elements of $\Aa$ it follows from Proposition~\ref{cspa.prop-posestim} that 
$$\sum_{k,l} \langle \rho(h\,f_k)\xi_k|\rho(h\,f_l)\xi_l\rangle
   \leq \|h\|^2\;
    \sum_{k,l} \langle \rho(f_k)\xi_k|\rho(f_l)\xi_l\rangle\,.
$$

\noindent Hence if $\psi=\sum_k \rho(f_k)\xi_k$ and if $T\psi=\sum_k \rho(h f_k)\xi_k$, then $\|T\psi\|\leq \|h\|\|\psi\|$.
It follows that if $\psi=0$ then $\|T\psi\|\leq \|h\| \|\psi\|=0$, namely if $\psi=\phi\in\Hh$, then $T\psi=T\phi$. Thus the map $\psi\in\Hh\to T\psi\in\Hh$ is well defined, linear by construction and bounded.
This map is denoted by $\rho(h)$ and it is an exercise that this defines $\rho$ as a weakly continuous $\Cc$-module $\ast$-homomorphism. The second part of the statement is just a consequence of the definition of the reduced and the full \CsS.
\hfill $\Box$

\vspace{.2cm}

\noindent The following Lemma will be used in the proof of the next result

\begin{lemma}
\label{cspa.lem-techn}
Let $X$ be a compact Hausdorff space. Let $h:X\to [0,1]$ be a function for which there is $x_0\in X$ such that $h(x_0)=0$ and $h$ is continuous at $x=x_0$. 
\noindent Then there exists a continuous function $\tth:X\to [0,1]$ such that $\tth(x_0)=0$ and $h(x)\leq \tth(x)$ for $x\in X$.
\end{lemma}

\noindent  {\bf Proof: }(i) Since $h$ is continuous at $x_0$, for any natural integer $n$ there is an open set $U_n$ containing $x_0$ such that if $x\in U_n$ then $0\leq h(x)<1/n$. Let then $V_n$ be defined by 
$$V_n=\bigcup_{m\geq n} U_m\,.
$$

\noindent By construction $V_n$ contains $x_0$, it is open and $V_{n+1}\subseteq V_n$ for all $n$'s. Similarly let $W_n$ be the open set defined by 
$$W_n=V_n\setminus \overline{V_{n+2}}\,.
$$

\noindent Then, for each $n$, $x_0\notin W_n$. Let $W$ denote the union of all the $W_n$'s (which coincides with $V_1\setminus \{x_0\}$). Its complement $K=X\setminus W$ is a compact set, containing $x_0$. By compactness it is possible to extract a finite subset $\{x_0, x_1,\cdots,x_l\}\subseteq K$ and for each $0\leq k\leq l$ an open neighborhood $O_k$ of $x_k$ such that $K\subseteq \bigcup_{k=0}^l O_k$. It is even possible to choose this open cover so that $x_0\notin O_k$ for $k\geq 1$. Then let $W_0$ denotes the open set $W_0=O_0\cup O_1\cup\cdots\cup O_l\setminus \{x_0\}$. In addition the family $\Ww=\{W_n\,;\; n\in \ZM_+\}$ is an open cover of the locally compact normal space $X\setminus \{x_0\}$. By construction, it is locally finite. Using the Urysohn's Lemma, let $(\phi_n)_{n\in\ZM_+}$ be a partition of unity subordinate to the open cover $\Ww$ of $X\setminus \{x_0\}$: namely for each $n\geq 0$, $\phi_n$ is continuous, $0\leq \phi_n(x)\leq 1$, the support of $\phi_n$ is included in $W_n$, and for any $x\in X\setminus \{x_0\}$, $\sum_{n=0}^\infty \phi_n(x)=1$. It is easy to check, using the definition of the $W_n$'s, that if $x\neq x_0$ is contained in $W_n$ but not in $W_{n-1}$, then it is not contained in $W_k$ if $k\neq n,n+1$. Therefore the previous infinite sum has at most two nonzero terms, so it is finite.

\vspace{.1cm}

\noindent (ii) For $j\in \NM$ let $\tth_j$ be the function defined by
$$\tth_j(x) = \phi_0(x)+ \sum_{n=1}^j \frac{1}{n}\;\phi_n(x)\,.
$$

\noindent By construction, $\tth_j$ is continuous, it vanishes at $x=x_0$, and $0\leq \tth_j(x)\leq \tth_{j+1}(x)$. In addition
$$\tth_j(x)\leq \sum_{n=0}^j \phi_n(x)\leq 1\,.
$$

\noindent Consequently the sequence $(\tth_j)_{j\in\NM}$ converges monotonically pointwise to a function $\tth$, such that $\tth(x_0)=0$ and $0\leq \tth(x)\leq 1$. However, for $1\leq i<j$ two integers
$$0\leq \tth_j(x)-\tth_i(x)=
   \sum_{n=i+1}^j \frac{1}{n}\;\phi_n(x)\leq 
    \frac{1}{i+1}\; \sum_{n=i+1}^j \phi_n(x)\leq 
     \frac{1}{i+1}\,.
$$ 

\noindent Consequently the sequence $(\tth_j)_{j\in\NM}$ is Cauchy for the uniform topology, so that $\sup_{x\in X}|\tth(x)-\tth_j(x)|\to 0$ and therefore $\tth$ is continuous. By construction, if $x\in W_n$, $0\leq h(x)\phi_n(x)< \phi_n(x)/n$ for $n\geq 1$, while for $n=0$, $0\leq h(x)\phi_0(x)\leq \phi_0(x)$. Consequently, summing over $n$ and using the pointwise convergence, this gives $0\leq h(x)\leq \tth(x)$.
\hfill $\Box$

\vspace{.2cm}

\noindent Let now $M\in\isG$ be a closed invariant subset. Then $\Ii_M$ will denote the set of functions $h\in\Cc$ such that $h(M)=0$. Thus $\Ii_M$ is a maximal (closed two-sided) ideal of $\Cc$. Similarly $\Jj_M$ will denote the set of $f\in\Aa$ such that $f(m,\gamma)=0$ for $\gamma\in \Gamma_F$. The following result holds

\begin{proposi}
\label{cspa.prop-quotient}
(a) The set $\Jj_M$ is a closed two-sided ideal of $\Aa$ which coincides with the closed two-sided ideal generated by $\Ii_M \Aa$.

\noindent (b) The quotient algebra $\Aa_M=\Aa/\Jj_M$ is $\ast$-isomorphic to the algebra 
$\Cc_c(\Gamma_M,\mu,\sigma)$.
\end{proposi}

\noindent {\bf Proof: }(i) That $\Jj_M$ is a two-sided ideal can be checked by inspection: namely if $f\in\Jj_M$ and if $g\in\Aa$ then $fg(M,\gamma)=0=gf(M,\gamma)$ for $\gamma\in \Gamma_M$. That it is closed comes from the fact that if $f$ is the limit of a net $f_\alpha$ of elements of $\Jj_M$ then, in particular, $f(M,\gamma)=\lim_{\alpha} f_\alpha(M,\gamma)=0$ for any $\gamma\in \Gamma_M$.

\vspace{.1cm}

\noindent (ii) Given $f\in\Jj_M$ there are $f_0,f_1,f_2,f_3$ such that $f=f_0-f_2+\imath(f_1-f_3)$ and $f_l(N,\gamma)\geq 0$ for all $l$'s and $(N,\gamma)$. This is because $f$ can be decomposed into real and imaginary part. In addition, each continuous real valued function $h$ on a topological space $X$ can be written as $h=h_+-h_-$, with $h_\pm(x) = \max\{\pm h(x),0\}$ for $x\in X$. Hence $h_\pm$ are continuous with support  contained in the support of $h$. In particular each $f_l$ is continuous and supported in the support of $f$. In addition, each $f_l$ vanishes at $M$ so that $f_l\in\Jj_M$. Thus, there is no loss of generality in assuming that $f(N,\gamma)\geq 0$ for all $(N,\gamma)\in \tsG$. Similarly, changing $f$ into $\lambda f$ if necessary, where $\lambda \in \RM_+$ there is no loss of generality in assuming that $0\leq f(N,\gamma)\leq 1$. 
Let $h:\isG\to [0,1]$ be defined by $h(N)=\sup_{\gamma\in \Gamma_N}f(N,\gamma)$. Then $h$ is lower semi-continuous, $h(M)=0$ and $0\leq h(N)\leq 1$. 

\vspace{.1cm}

\noindent (iii) Let $K$ denote the support of $f$. By assumption, $K$ is compact in $\tsG$. It follows that $h$ is continuous at $N=M$. Indeed, thanks to the Bolzano-We\"{\i}erstrass theorem , given $N\in\isG$, there is $\gamma_N\in\Gamma_N$ such that $f(N,\gamma)$ reaches its supremum at $(N,\gamma_N)$. Hence $h(N)=f(N,\gamma_N)$. Let $(N_n)_{n\in\NM}$ be a sequence converging to $M$ in $\isG$. Thanks to the compactness of $K$, by extracting a sub-sequence if necessary, the corresponding $\gamma_{N_n}$ converge to a point $\eta\in\Gamma_M$. Since $f$ is continuous it follows that $\lim_{n\to\infty} h(N_n)=\lim_{n\to\infty} f(N_n,\gamma_{N_n})=f(M,\eta)=0$. Since $\Gamma$ is second countable, it follows, from the properties of the Hausdorff topology (see Proposition~\ref{cspa.prop-isGamma}), that $\isG$ is second countable as well, so that this convergence is sufficient to prove the continuity of $h$ at $N=M$.

\vspace{.1cm}

\noindent (iv) Thanks to Lemma~\ref{cspa.lem-techn}, let then $\tth\in \Ii_M$ be such that $\tth(N)\geq h(N)$. Then given $0<\epsilon<1$, let $f_\epsilon$ be defined by
$$f_\epsilon(N,\gamma)= \tth(N)^\epsilon\; f(N,\gamma)^{1-\epsilon}\,.
$$

\noindent Hence 
$f_\epsilon\in\Ii_M\Aa$. 
 It follows that if $R$ denotes the ratio $h(N)/f(N,\gamma)$ then $R\geq 1$ (it can be infinite) and  
$$0\leq (f_\epsilon-f)(N,\gamma)=
   (R^\epsilon-1)f(N,\gamma)=
    \int_0^\epsilon 
     \ln(R) R^s f(N,\gamma)\,ds\leq
      \epsilon R^{\epsilon-1}\ln(R) \tth(n)\,.
$$

\noindent where $f(N,\gamma)= R^{-1} \tth(N)$ is used. The support of $f_\epsilon$ is contained in the support of $f$. Moreover, $R^{\epsilon -1}\ln(R)\leq \left(e(1-\epsilon)\right)^{-1}$. Therefore $\lim_{\epsilon\to 0} f_\epsilon=f$ in the topology of $\Aa$. Hence $f\in\Ii_F\Aa$.

\vspace{.1cm}

\noindent (v) For $f\in \Aa$ let $f_M$ denote the function $f_M(\gamma)=f(M,\gamma)$ for $\gamma\in \Gamma_M$. The map $p_M:f\to f_M$ is a continuous $\ast$-homomorphism with kernel $\Jj_M$ by construction. Hence it defines an injective $\ast$-homomorphism from $\Aa_M$ into $\Cc_c(\Gamma_M,\mu,\sigma)$. 
Conversely, if $g\in \Cc_c(\Gamma_M)$, with support in the compact subset $K\subseteq \Gamma_M$ it can be seen as a function defined on the compact subset $\{M\}\times K\subseteq \tsG$. Since $\Gamma$ is handy, the Tietze Extension Theorem applies, so that there is a function $\hg\in\Aa$ such that $\hg_M=g$. Hence $p_M$ is onto.
\hfill $\Box$

 \subsection{The upper semi-continuity for the full $C^\ast$-algebra}
 \label{cspa.ssect-full}

\noindent The first part of this sections is devoted to prove the exactness of a sequence of groupoid \Css defined with a cocycle, see Corollary~\ref{cspa.cor-exact}. With this at hand, the upper semi-continuity of the full norm is shown in Proposition~\ref{cspa.prop-fullusc}.

\vspace{.1cm}

\noindent The definition of the full \Cs was provided in Section~\ref{cspa.ssect-convol}. As a preparation, the following result is elementary but is needed.

\begin{lemma}
\label{cspa.lem-ideal}
Let $0\to J\stackrel{i}{\to} A\stackrel{p}{\to} B\to 0$ be a sequence of \Css with $i,p$ being $\ast$-homomorphisms such that $i$ is injective, $p$ is surjective and $p\circ i=0$. If, for any $\ast$-representation $\pi$ of $A$ vanishing on $\Ima(i)$, there is a $\ast$-representation $\pi_B$ of $B$ such that $\pi_B\circ p=\pi$, then the previous sequence is exact at $A$, namely $\Ker(p)=\Ima(i)$. 
\end{lemma}

\noindent  {\bf Proof: } From the assumptions it follows easily that $\Ima(i)\subseteq \Ker(p)$. If $i(J)\neq \Ker(p)$, there is $a\in \Ker(p)$ such that $a\notin \Ima(i)$. Consequently there is a $\ast$-representation $\pi$ of $A$ vanishing on $\Ima(i)$ such that $\pi(a)\neq 0$. By the assumption, there is a representation $\pi_B$ of $B$ such that $\pi(a)=\pi_B\circ p(a)=0$, leading to a contradiction.
\hfill $\Box$

\vspace{.2cm}

\noindent Like in Section~\ref{cspa.ssect-algB} let $\Aa=\Cc_c(\tsG,\mu,\sigma)$. Given $M$ a closed invariant subset of $\Gamma^{(0)}$, let $\Jj_M$ denotes the ideal of functions $f\in\Aa$ such that $f(M,\gamma)=0$ for all $\gamma\in \Gamma_M$. In addition, $\Aa_M$ denotes the $\ast$-algebra $\Cc_c(\Gamma_m,\mu_M,\sigma_M)$ where $\sigma_M$ represents the restriction of the cocycle $\sigma$ on $\Gamma_M$ and, similarly, $\mu_M$ denotes the restriction of the Haar system $\mu$ to $\Gamma_M$. Then the following result holds

\begin{lemma}
\label{cspa.lem-extension}
Let $i_M: \Jj_M\to \Aa$ be the injection map and let $p_M:\Aa\to \Aa_M$ be the evaluation map at $M$. Then:

(a) The sequence of topological $\ast$-algebras $0\to \Jj_M\stackrel{i_M}{\to} \Aa\stackrel{p_M}{\to} \Aa_M\to 0$ is exact.

(b) The maps $i_M$ and $p_M$ extends as $\ast$-homomorphisms to the completions of the previous algebras for both the reduced and the full norms. In both case $i_M$ is injective and $p_M\circ i_M=0$.
\end{lemma}

\noindent  {\bf Proof: } (i) The map $i_M$ is defined as $i_M(f)(N,\gamma)=f(N,\gamma)$, for $f\in \Jj_M$, which makes sense since $\Jj_M\subseteq \Aa$. The map $p_M$ is defined by $p_M(f)(\gamma)=f(M,\gamma)$ for any $\gamma \in \Gamma_M$. 
Proposition~\ref{cspa.prop-quotient}~(b) implies $\Ker(p_M)=\Ima(i_M)$.

\vspace{.1cm}

\noindent (ii) The map $p_M$ is onto which will prove the exactness at $\Aa_M$. Namely if $g\in\Aa_M$, it defines a continuous function with compact support on $\Gamma_M$, which can be seen as a continuous function with compact support on $\{M\}\times \Gamma_M\subseteq \tsG$. Thanks to the Tietze extension theorem, there is $f\in\Aa$ coinciding with $g$ on $\{M\}\times \Gamma_M$. Then, by definition $p_M(f)=g$. Together with (i) this proves (a).

\vspace{.1cm}

\noindent (iii) Given a bounded, $\ast$-representation $\pi$ of $\Aa$, the map $f\in\Jj_M\mapsto \pi\circ i_M(f)$ gives a bounded, $\ast$-representation of $\Jj_M$. In particular, if $N\in\isG$ and if $x\in N$, let $\pi_{N,x}$ denote the left regular representation of $\tsG$ on the fiber space $L^2(\tsG^{N,x},\mu^{N,x})$. Then the reduce norm $\|g\|_{red}$ of $g\in \Aa$ is obtained by taking the supremum over $(N,x)$ of $\|\pi_{N,x}(g)\|$. It then follows that $\|i_M(f)\|_{red}=\sup_{N,x}\|\pi_{N,x}(i(f))\|\leq \|f\|_{red}$. Similarly, taking the supremum over all bounded, $\ast$-representations $\pi$ of $\Aa$, this argument gives $\|i_M(f)\|=\sup_{\pi}\|\pi(i_M(f))\|\leq \|f\|$. Hence $i_M$ extends as a $\ast$-homomorphism to the completions of $\Jj_M$ and $\Aa$ {\em w.r.t.} both the reduced and the full norm. By construction $i_M(f)=f$ for $f\in\Jj_M$, therefore the same identity holds in both completions. Hence $i_M$ is injective and therefore it is an isometry in both completions.

\vspace{.1cm}

\noindent (iv) A similar argument holds for $p_M$. If $\rho$ is a bounded, $\ast$-representation of $\Aa_M$, then $\rho\circ p_M$ is also a bounded, $\ast$-representation of $\Aa$ as can be checked by inspection, using the invariance of the closed set $M$. As a result the same argument shows that $p_M$ extends to both the reduced and the full completions of $\Aa$ and $\Aa_M$ as a $\ast$-homomorphism. The relation $p_M\circ i_M=0$ holds on $\Jj_M$ so it must hold on both completions as well by continuity.
\hfill $\Box$

\vspace{.2cm}

\noindent The following Lemma is using an argument due to Renault (see \cite{Re80}, p. 102) used in a different context.

\begin{lemma}
\label{cspa.lem-kernel}
Let $A$ denotes the full \Cs $A=\CG^\ast(\tsG,\mu,\sigma)$ and let $J_M$ be the closure of $\Jj_M$ in $A$ for $M\in\isG$. Similarly let $A_M$ be the full algebra $\CG^\ast(\Gamma_M,\mu_M,\sigma_M)$. If $\pi$ a $\ast$-representation of $A$ vanishing on $J_M$, then there is a $\ast$-representation $\rho$ of $A_M$ such that $\rho\circ p_M=\pi$.  
\end{lemma}

\noindent  {\bf Proof: }The restriction of $\pi$ on $\Aa$ is a bounded $\ast$-representation of $\Aa$. Since $\pi$ vanishes on $J_M$, it also vanishes on $\Jj_M$. Let then $\rho$ be defined, for $g\in\Aa_F$, by 
$\rho(g)=\pi(f)$ for any $f\in\Aa$ coinciding with $g$ on $\Gamma_M$, namely for any $f\in \Aa$ such that $f(M,\gamma)=g(\gamma)$ for $\gamma\in\Gamma_M$. This definition makes sense since changing $f$ into $f'$ implies $f-f'\in \Jj_M$ namely $\pi(f)=\pi(f')$. This construction gives $\rho\circ p_M=\pi$ when restricted to $\Aa$. By completion, $\rho$ extends to a $\ast$-representation of the full algebra $A_M$. The relation $\rho\circ p_M=\pi$, valid on $\Aa$, extends to $A$ since both $\pi$ and $p_M$ are $\ast$-homomorphisms.
\hfill $\Box$

\begin{coro}
\label{cspa.cor-exact}
The sequence of \Css $0\to J_M\stackrel{i_M}{\to} A \stackrel{p_M}{\to} A_M \to 0$ is exact. In particular, $A_M$ is $\ast$-isomorphic to the quotient algebra $A/J_M$. 
\end{coro}

\noindent  {\bf Proof: }Thanks to Lemma~\ref{cspa.lem-ideal}, it follows from Lemma~\ref{cspa.lem-kernel} that $\Ker(p_M)=\Ima(i_M)$ in this sequence. Hence it is exact at $A$. It remains to prove that $p_M:A\to A_M$ is surjective. Since $p_M$ is a $\ast$-homomorphism, the image $p_M(A)$ is closed (\cite{KR1}, Theorem 4.1.9). Since $\Aa$ is dense in $A$ and since $p_M(\Aa)=\Aa_M$ (thanks to Lemma~\ref{cspa.lem-extension}) is dense in $A_M$, it follows that $\Aa_M=p_M(\Aa)\subseteq {p_M(A)}\subseteq A_M$, so that $A_M=\overline{\Aa_M}=\overline{p_M(\Aa)}\subseteq p_M(A)\subseteq A_M$. Hence $p_M(A)=A_M$ and $p_M$ is surjective.
\hfill $\Box$

\vspace{.2cm}

\noindent The following result is using an argument that can be found in Rieffel (\cite{Rie89a}, Proposition 1.2) and the author credits Varela for it.

\begin{proposi}
\label{cspa.prop-fullusc}
For any 
$f\in \CG^\ast(\Gamma,\mu,\sigma)$, 
the map $M\in\isG\mapsto \|p_M(f)\|\in\RM_+$ is upper semi-continuous.
\end{proposi}

\noindent {\bf Proof: }Given $f\in A$, let $M\in\isG$ and let $\epsilon >0$. Using the exactness of the sequence in Corollary~\ref{cspa.cor-exact}, it follows that $p_M$ defines a $\ast$-isomorphism  from the quotient algebra $A/J_M$ into  $A_M$. By definition of the quotient norm, $\|p_M(f)\|=\inf\{\|f-g\|\,;\, g\in J_M\}$. Now, $J_M$ is the closure of $\Jj_M$. Moreover, thanks to Proposition~\ref{cspa.prop-quotient}, $\Jj_M=\Ii_M\Aa$. Hence there is a finite sum $g=\sum_{n=1}^k h_n g_n$ with $h_n\in \Ii_M$ and $g_n\in \Aa$, such that
$$\|p_M(f)\|\geq \|f-g\|-\epsilon\,.
$$

\noindent Let $\tth\in\Cc(\isG)$ be chosen such that $0\leq \tth(N)\leq 1$ and such that $\tth(N)=1$ in a neighborhood $U$ of $M$. In particular $\|\tth f\|\leq \|f\|$. Choosing $U$ and the support of $\tth$ small enough, it is even possible to assume that $\|\tth g\| \leq \epsilon$, since $\tth(M)h_n(M)=0$ and since $h_n$ is continuous for $1\leq n\leq k$. Consequently,
\begin{eqnarray*}
\|p_M(f)\| & \geq & \|f-g\|-\epsilon\geq 
           \|\tth(f-g)\|-\epsilon\geq \|\tth f\|-2\epsilon\\
& = & \|f-(1-\tth)f\|-2\epsilon
\end{eqnarray*}

\noindent Since $1-\tth$ vanishes at $N\in U$, the element $(1-\tth)f\in J_N$ for $N\in U$, so that $\|f-(1-\tth)f\|\geq \|p_N(f)\|$, leading to 
$$\|p_M(f)\|\geq \|p_N(f)\|-2\epsilon\,,
   \hspace{2cm}
     \forall N\in U\,. 
$$
\hfill $\Box$

\begin{rem}
\label{cspa.rem-tautfield2}
{\em Theorem~\ref{cspa.th-tautfield_2} and Proposition~\ref{cspa.prop-fullusc} imply that if $(\Gamma,T,p)$ is a continuous field of handy groupoids with continuous field of $2$-cocycle $\sigma$, then $t\mapsto\|f_t\|$ is upper semi-continuous for $f\in\CG^\ast(\Gamma,\mu,\sigma)$ where $f_t$ is the restriction of $f$ onto the closed subgroupoid $\Gamma_t:=p^{-1}(\{t\})$.}
\hfill $\Box$
\end{rem}

 \subsection{The lower semi-continuity for the reduced $C^\ast$-algebra}
 \label{cspa.ssect-reduced}

\noindent The lower semi-continuity of the reduced norm for the tautological groupoid is proved in this section.

\vspace{.1cm}

\noindent Using the previous notations, let $f\in\Aa$. Then $p_M(f)\in\Aa_M$ is the function $\gamma\in\Gamma_M\mapsto f(M,\gamma)$.

\begin{proposi}
\label{cspa.prop-redlsc}
Given $f\in\CG_{red}^\ast(\tsG,\mu,\sigma)$ the map $M\in\isG\mapsto \|p_M(f)\|_{red}$ is lower semi-continuous.
\end{proposi}

\noindent {\bf Proof: }(i) The left-regular representation of $\Aa_M$ is acting on the Hilbert space $\Hh_{M,x}=L^2(\Gamma^x,\mu^x)$ whenever $x\in M$. The corresponding operator is given by the following formula (see Section~\ref{cspa.ssect-red}, eq.~\eqref{cspa.eq-redrep}), where $\phi\in\Hh_{M,x}$
$$\pi_{M,x}(f)\phi(\gamma) =
   \int_{\Gamma^x} f(M,\gamma^{-1}\circ \eta)
    \sigma_M(\gamma,\gamma^{-1}\circ \eta)\;
     \phi(\eta)\;d\mu^x(\eta)\,.
$$

\noindent In particular, if $\phi'\in\Hh_{F,x}$ then the inner product is given by
$$\langle \phi|\pi_{M,x}(f)\phi'\rangle_{\Hh_{M,x}}=
   \int_{\Gamma^x\times \Gamma^x}
    \overline{\phi(\gamma)}\,
     f(M,\gamma^{-1}\circ \eta)
      \sigma_M(\gamma,\gamma^{-1}\circ \eta)\;
       \phi'(\eta)\;d\mu^x(\gamma)\,d\mu^x(\eta)\,. 
$$

\noindent (ii) At this point, it is worth reminding that $\Cc_c(\Gamma^x)\subseteq \Hh_{M,x}$ and the former is dense in the latter. Since $\tsG$ is normal, using the Tietze Extension Theorem, any continuous function with compact support defined on $\Gamma^x$, seen as the subset $\{M\}\times \Gamma^x\subseteq \tsG$, can be extended as a continuous function with compact support in $\tsG$. For such a function, say $\varphi$, its restriction $\varphi_{M,x}$ to $\{M\}\times \Gamma^x$ is an element of $\Hh_{M,x}$, while $\varphi$ itself belongs to $\Aa$. Therefore if $\varphi,\psi\in\Aa$ the following formula holds
$$\langle \varphi_{M,x}|\pi_{M,x}(f)\psi_{M,x}\rangle_{\Hh_{M,x}}=
   \int_{\Gamma^x\times \Gamma^x}
    \overline{\varphi(M,\gamma)}\,
     f(M,\gamma^{-1}\circ \eta)
      \sigma_M(\gamma,\gamma^{-1}\circ \eta)\;
       \psi(M,\eta)\;d\mu^x(\gamma)\,d\mu^x(\eta)\,. 
$$

\noindent From the definition of a Haar system, in particular the axiom (H2), it follows immediately that the map $(M,x)\in\tsGO\mapsto  \langle \varphi_{M,x}|\pi_{M,x}(f)\psi_{M,x}\rangle_{\Hh_{M,x}}$ is continuous whenever $f,\varphi,\psi\in\Aa$. Similarly, the map $(M,x)\in\tsGO\mapsto  \|\varphi_{M,x}\|_{\Hh_{M,x}}\in\RM_+$ is also continuous. Using the density result it follows that
$$\|\pi_{M,x}(f)\|=
  \sup_{\varphi,\psi\in\Aa} 
   \left|
   \frac{
    \langle \varphi_{M,x}|\pi_{M,x}(f)\psi_{M,x}\rangle_{\Hh_{M,x}}}
      {\|\varphi_{M,x}\|_{\Hh_{M,x}}\,\|\psi_{M,x}\|_{\Hh_{M,x}}}\right|\,.
$$

\noindent In addition, the reduced $C^\ast$-norm of $p_M(f)$ is given by

$$\|p_M(f)\|_{red}=\sup_{x\in M}\|\pi_{M,x}(f)\|\,.
$$

\noindent Consequently the map $M\in\isG\mapsto \|p_M(f)\|_{red}$ is lower semi-continuous, as it is the supremum over a family of continuous functions.

\vspace{.1cm}

\noindent (iii) If now $f\in\CG_{red}^\ast(\tsG,\mu,\sigma)$, it can be approximated by a sequence $(f_n)_{n\in\NM}$ of elements of $\Aa$ in the reduced norm. Hence, given $\epsilon >0$, there is $n_0\in\NM$ such that for $n\geq n_0$,  $\|f-f_n\|_{red}<\epsilon/3$. Since $p_M$ extends as a $\ast$-homomorphism from $\Aa$ to the reduced \CS, it follows immediately that $\|p_M(f)-p_M(f_n)\|_{red}<\epsilon/3$ for $n\geq n_0$. Now choosing $n\geq n_0$, by lower semi-continuity, there is a neighborhood $U$ of $M$ in $\isG$ such that for $N\in U$
$$\|p_M(f_n)\|_{red}\leq \|p_N(f_n)\|_{red} +\epsilon/3\,,
   \hspace{1cm}\Rightarrow\hspace{1cm}
    \|p_M(f)\|_{red}\leq \|p_N(f)\|_{red} +\epsilon\,.
$$
\noindent This finishes the proof.
\hfill $\Box$

\begin{rem}
\label{cspa.rem-tautfield3}
{\em Theorem~\ref{cspa.th-tautfield_2} and Proposition~\ref{cspa.prop-redlsc} imply that if $(\Gamma,T,p)$ is a continuous field of handy groupoids with continuous field of $2$-cocycle $\sigma$, then $t\mapsto\|f_t\|_{red}$ is lower semi-continuous for $f\in\CG_{red}^\ast(\tsG,\mu,\sigma)$ where $f_t$ is the restriction of $f$ onto the closed subgroupoid $\Gamma_t:=p^{-1}(\{t\})$.}
\hfill $\Box$
\end{rem}

 \subsection{Proof of Theorem~\ref{cspa.theo-contSp}}
 \label{cspa.ssect-Prtheo}

\noindent Due to Theorem~\ref{cspa.th-tautfield} and Proposition~\ref{cspa.prop-cfcoc}, $(\ts(\Gamma),\isG,p_\Gamma)$ is a continuous field of handy groupoids with continuous field of $2$-cocycle $\hat{\sigma}$. The map $i:\Cc_c(\Gamma,\mu,\sigma)\to\Cc_c(\tsG,\mu,\hat{\sigma})$ defined by $i(f)(M,\gamma):=f(\gamma)$ extends to an injective $\ast$-homomorphism from $\CG_{red}^\ast(\Gamma,\mu,\sigma)$ into $\CG_{red}^\ast(\tsG,\mu,\hat{\sigma})$. Thus, $\spec(f)=\spec(i(f))$ follows and so it suffices to show continuity of the map 
$$\Sigma_f:\isG\to\ks(\CM)\,,\; M\mapsto\spec(p_M(f))\,,
$$ 
for $f\in \CG_{red}^\ast(\tsG,\mu,\hat{\sigma})$ normal.

\vspace{.1cm}

\noindent The continuity of the norm $t\mapsto \|\phi(f)\|$ of a normal element $f\in\CG_{red}^\ast(\tsG,\mu,\hat{\sigma})$ for any continuous $\phi:\CM\to\CM$ follows directly from Propositions~\ref{cspa.prop-fullusc} and Proposition~\ref{cspa.prop-redlsc} whenever the reduced and full \Css agree. A sufficient condition for this is that the groupoid is amenable, see \cite{AR00,AR01} and \cite{BH14}, Corollary~4.3, for groupoid \Css with normalized $2$-cocycle. However, there is a subtlety depending upon whether or not the reduced \Cs is unital or not unital.

\vspace{.1cm}

\noindent (i) If the reduced \Cs is unital, and if $f$ is self-adjoint, the continuity of the spectrum follows by Theorem~\ref{cspa.th-spcont} proved in \cite{BB16}. For normal elements, see \cite{Di69} and a more precise statement in \cite[Thm.~2.7.9]{Bec16}, implying the continuity of the spectra.

\vspace{.1cm}

\noindent (ii) If the reduced \Cs is not unital, the usual way consists in adding a unit (in the Abelian case it consists in replacing the locally compact spectrum by its Alexandrov compactification adding one point at infinity). By convention then, the spectrum $\spec(p_M(f))$ contains always $0$. Then the continuity of the spectrum follows by \cite[Thm.~2.7.9]{Bec16} but only functions with $\phi(0)=0$ can be used \cite[Prop.~10.3.3]{Di69}.
\hfill $\Box$

\begin{rem}
\label{cspa.rem-unit}
{\em If $\Gamma$ is a handy {\'e}tale groupoid, the reduced \Cs is unital, where the characteristic function $\chi_{\Gamma^{(0)}}$ of the unit set is the unit. Note that $\chi_{\Gamma^{(0)}}\in\Cc_c(\Gamma,\mu,\sigma)$ since the unit set is compact and open \cite{Re09}.
}
\hfill $\Box$
\end{rem}

\vspace{.5cm}
\input{spcontv18.tex}

\end{document}